\documentclass{amsproc}
\usepackage{amsfonts}

\theoremstyle{plain}

\newtheorem{corollary}{Corollary}

\newtheorem{definition}{Definition}

\newtheorem{lemma}{Lemma}

\newtheorem{proposition}{Proposition}
\newtheorem{remark}{Remark}

\newtheorem{theorem}{Theorem}
\numberwithin{equation}{section}
\input{tcilatex}
\usepackage{graphicx}

\begin{document}
\title[Reverses of the Triangle Inequality]{Reverses of the Triangle
Inequality in Banach Spaces}
\author[S.S. Dragomir]{Sever S. Dragomir}
\address{School of Computer Science and Mathematics\\
Victoria University of Technology\\
PO Box 14428\\
Melbourne VIC 8001, Australia.}
\email{sever.dragomir@vu.edu.au}
\urladdr{http://rgmia.vu.edu.au/SSDragomirWeb.html}
\subjclass[2000]{Primary 46B05, 46C05; Secondary 26D15, 26D10}

\begin{abstract}
Recent reverses for the discrete generalised triangle inequality and its
continuous version for vector-valued integrals in Banach spaces are
surveyed. New results are also obtained. Particular instances of interest in
Hilbert spaces and for complex numbers and functions are pointed out as well.
\end{abstract}

\maketitle

\section{Introduction}

The \textit{generalised triangle inequality}, namely

\begin{equation*}
\left\Vert \sum_{i=1}^{n}x_{i}\right\Vert \leq \sum_{i=1}^{n}\left\Vert
x_{i}\right\Vert ,
\end{equation*}%
provided $\left( X,\left\Vert .\right\Vert \right) $ is a normed linear
space over the real or complex filed $\mathbb{K=R}$, $\mathbb{C}$ and $%
x_{i},i\in \left\{ 1,...,n\right\} $ are vectors in $X$ plays a fundamental
role in establishing various analytic and geometric properties of such
spaces.

With no less importance, the \textit{continuous} version of it, i.e.,

\begin{equation}
\left\Vert \int_{a}^{b}f\left( t\right) dt\right\Vert \leq
\int_{a}^{b}\left\Vert f\left( t\right) \right\Vert dt,  \label{C}
\end{equation}%
where $f:\left[ a,b\right] \subset \mathbb{R}\rightarrow X$ is a Bochner
measurable function on the compact interval $\left[ a,b\right] $ with values
in the Banach space $X$ and $\left\Vert f\left( .\right) \right\Vert $ is
Lebesgue integrable on $\left[ a,b\right] ,$ is crucial in the Analysis of
vector-valued functions with countless applications in Functional Analysis,
Operator Theory, Differential Equations, Semigroups Theory and related
fields.

Surprisingly enough, the reverses of these, i.e., inequalities of the
following type%
\begin{equation*}
\sum_{i=1}^{n}\left\Vert x_{i}\right\Vert \leq C\left\Vert
\sum_{i=1}^{n}x_{i}\right\Vert ,\int_{a}^{b}\left\Vert f\left( t\right)
\right\Vert dt\leq C\left\Vert \int_{a}^{b}f\left( t\right) dt\right\Vert ,
\end{equation*}%
with $C\geq 1,$ which we call \textit{multiplicative reverses}, or%
\begin{equation*}
\sum_{i=1}^{n}\left\Vert x_{i}\right\Vert \leq \left\Vert
\sum_{i=1}^{n}x_{i}\right\Vert +M,\int_{a}^{b}\left\Vert f\left( t\right)
\right\Vert dt\leq \left\Vert \int_{a}^{b}f\left( t\right) dt\right\Vert +M,
\end{equation*}%
with $M\geq 0,$ which we call \textit{additive reverses}, under suitable
assumptions for the involved vectors or functions, are far less known in the
literature.

It is worth mentioning though, the following reverse of the generalised
triangle inequality for complex numbers%
\begin{equation*}
\cos \theta \sum_{k=1}^{n}\left\vert z_{k}\right\vert \leq \left\vert
\sum_{k=1}^{n}z_{k}\right\vert ,
\end{equation*}%
provided the complex numbers $z_{k},$ $k\in \left\{ 1,\dots ,n\right\} $
satisfy the assumption%
\begin{equation*}
a-\theta \leq \arg \left( z_{k}\right) \leq a+\theta ,\ \ \text{for any \ }%
k\in \left\{ 1,\dots ,n\right\} ,
\end{equation*}%
where $a\in \mathbb{R}$ and $\theta \in \left( 0,\frac{\pi }{2}\right) $ was
first discovered by M. Petrovich in 1917, \cite{P} (see \cite[p. 492]{MPF})
and subsequently was rediscovered by other authors, including J. Karamata 
\cite[p. 300 -- 301]{K}, H.S. Wilf \cite{W}, and in an equivalent form by M.
Marden \cite{MA}. Marden and Wilf have outlined in their work the important
fact that reverses of the generalised triangle inequality may be
successfully applied to the location problem for the roots of complex
polynomials.

In 1966, J.B. Diaz and F.T. Metcalf \cite{DM} proved the following reverse
of the triangle inequality in the more general case of inner product spaces:

\begin{theorem}[Diaz-Metcalf, 1966]
\label{0.ta}Let $a$ be a unit vector in the inner product space $\left(
H;\left\langle \cdot ,\cdot \right\rangle \right) $ over the real or complex
number field $\mathbb{K}$. Suppose that the vectors $x_{i}\in H\backslash
\left\{ 0\right\} ,$ $i\in \left\{ 1,\dots ,n\right\} $ satisfy%
\begin{equation*}
0\leq r\leq \frac{\func{Re}\left\langle x_{i},a\right\rangle }{\left\Vert
x_{i}\right\Vert },\ \ \ \ \ i\in \left\{ 1,\dots ,n\right\} .
\end{equation*}%
Then%
\begin{equation*}
r\sum_{i=1}^{n}\left\Vert x_{i}\right\Vert \leq \left\Vert
\sum_{i=1}^{n}x_{i}\right\Vert ,
\end{equation*}%
where equality holds if and only if%
\begin{equation*}
\sum_{i=1}^{n}x_{i}=r\left( \sum_{i=1}^{n}\left\Vert x_{i}\right\Vert
\right) a.
\end{equation*}
\end{theorem}

A generalisation of this result for orthonormal families is incorporated in
the following result \cite{DM}.

\begin{theorem}[Diaz-Metcalf, 1966]
\label{0.tb}Let $a_{1},\dots ,a_{n}$ be orthonormal vectors in $H.$ Suppose
the vectors $x_{1},\dots ,x_{n}\in H\backslash \left\{ 0\right\} $ satisfy%
\begin{equation*}
0\leq r_{k}\leq \frac{\func{Re}\left\langle x_{i},a_{k}\right\rangle }{%
\left\Vert x_{i}\right\Vert },\ \ \ \ \ i\in \left\{ 1,\dots ,n\right\} ,\
k\in \left\{ 1,\dots ,m\right\} .
\end{equation*}%
Then%
\begin{equation*}
\left( \sum_{k=1}^{m}r_{k}^{2}\right) ^{\frac{1}{2}}\sum_{i=1}^{n}\left\Vert
x_{i}\right\Vert \leq \left\Vert \sum_{i=1}^{n}x_{i}\right\Vert ,
\end{equation*}%
where equality holds if and only if%
\begin{equation*}
\sum_{i=1}^{n}x_{i}=\left( \sum_{i=1}^{n}\left\Vert x_{i}\right\Vert \right)
\sum_{k=1}^{m}r_{k}a_{k}.
\end{equation*}
\end{theorem}

Similar results valid for semi-inner products may be found in \cite{KA}, 
\cite{KUR} and \cite{M}.

Now, for the scalar continuous case.

It appears, see \cite[p. 492]{MPF}, that the first reverse inequality for (%
\ref{C}) in the case of complex valued functions was obtained by J. Karamata
in his book from 1949, \cite{K}. It can be stated as%
\begin{equation*}
\cos \theta \int_{a}^{b}\left\vert f\left( x\right) \right\vert dx\leq
\left\vert \int_{a}^{b}f\left( x\right) dx\right\vert
\end{equation*}%
provided%
\begin{equation*}
-\theta \leq \arg f\left( x\right) \leq \theta ,\ \ x\in \left[ a,b\right]
\end{equation*}%
for given $\theta \in \left( 0,\frac{\pi }{2}\right) .$

This result has recently been extended by the author for the case of Bochner
integrable functions with values in a Hilbert space $H.$ If by $L\left( %
\left[ a,b\right] ;H\right) ,$ we denote the space of Bochner integrable
functions with values in a Hilbert space $H,$ i.e., we recall that $f\in
L\left( \left[ a,b\right] ;H\right) $ if and only if $f:\left[ a,b\right]
\rightarrow H$ is Bochner measurable on $\left[ a,b\right] $ and the
Lebesgue integral $\int_{a}^{b}\left\Vert f\left( t\right) \right\Vert dt$
is finite, then%
\begin{equation}
\int_{a}^{b}\left\Vert f\left( t\right) \right\Vert dt\leq K\left\Vert
\int_{a}^{b}f\left( t\right) dt\right\Vert ,  \label{R1}
\end{equation}%
provided that $f$ satisfies the condition%
\begin{equation*}
\left\Vert f\left( t\right) \right\Vert \leq K\func{Re}\left\langle f\left(
t\right) ,e\right\rangle \ \ \ \text{for a.e. }t\in \left[ a,b\right] ,
\end{equation*}%
where $e\in H,$ $\left\Vert e\right\Vert =1$ and $K\geq 1$ are given. The
case of equality holds in (\ref{R1}) if and only if%
\begin{equation*}
\int_{a}^{b}f\left( t\right) dt=\frac{1}{K}\left( \int_{a}^{b}\left\Vert
f\left( t\right) \right\Vert dt\right) e.
\end{equation*}%
The aim of the present paper is to survey some of the recent results
concerning multiplicative and additive reverses for both the discrete and
continuous version of the triangle inequalities in Banach spaces. New
results and applications for the important case of Hilbert spaces and for
complex numbers and complex functions have been provided as well.

\section{Diaz-Metcalf Type Inequalities}

In \cite{DM}, Diaz and Metcalf established the following reverse of the
generalised triangle inequality in real or complex normed linear spaces.

\begin{theorem}[Diaz-Metcalf, 1966]
If $F:X\rightarrow \mathbb{K}$, $\mathbb{K}=\mathbb{R},\mathbb{C}$ is a
linear functional of a unit norm defined on the normed linear space $X$
endowed with the norm $\left\Vert \cdot \right\Vert $ and the vectors $%
x_{1},\dots ,x_{n}$ satisfy the condition%
\begin{equation}
0\leq r\leq \func{Re}F\left( x_{i}\right) ,\ \ \ \ \ i\in \left\{ 1,\dots
,n\right\} ;  \label{1.1.1}
\end{equation}%
then%
\begin{equation}
r\sum_{i=1}^{n}\left\Vert x_{i}\right\Vert \leq \left\Vert
\sum_{i=1}^{n}x_{i}\right\Vert ,  \label{1.1.2}
\end{equation}%
where equality holds if and only if both%
\begin{equation}
F\left( \sum_{i=1}^{n}x_{i}\right) =r\sum_{i=1}^{n}\left\Vert
x_{i}\right\Vert   \label{1.1.3}
\end{equation}%
and%
\begin{equation}
F\left( \sum_{i=1}^{n}x_{i}\right) =\left\Vert
\sum_{i=1}^{n}x_{i}\right\Vert .  \label{1.1.4}
\end{equation}
\end{theorem}

If $X=H$, $\left( H;\left\langle \cdot ,\cdot \right\rangle \right) $ is an
inner product space and $F\left( x\right) =\left\langle x,e\right\rangle ,$ $%
\left\Vert e\right\Vert =1,$ then the condition (\ref{1.1.1}) may be
replaced with the simpler assumption%
\begin{equation}
0\leq r\left\Vert x_{i}\right\Vert \leq \func{Re}\left\langle
x_{i},e\right\rangle ,\qquad i=1,\dots ,n,  \label{1.1.5}
\end{equation}%
which implies the reverse of the generalised triangle inequality (\ref{1.1.2}%
). In this case the equality holds in (\ref{1.1.2}) if and only if \cite{DM}%
\begin{equation}
\sum_{i=1}^{n}x_{i}=r\left( \sum_{i=1}^{n}\left\Vert x_{i}\right\Vert
\right) e.  \label{1.1.6}
\end{equation}

\begin{theorem}[Diaz-Metcalf, 1966]
Let $F_{1},\dots ,F_{m}$ be linear functionals on $X,$ each of unit norm. As
in \cite{DM}, let consider the real number $c$ defined by%
\begin{equation*}
c=\sup_{x\neq 0}\left[ \frac{\sum_{k=1}^{m}\left\vert F_{k}\left( x\right)
\right\vert ^{2}}{\left\Vert x\right\Vert ^{2}}\right] ;
\end{equation*}%
it then follows that $1\leq c\leq m.$ Suppose the vectors $x_{1},\dots ,x_{n}
$ whenever $x_{i}\neq 0,$ satisfy%
\begin{equation}
0\leq r_{k}\left\Vert x_{i}\right\Vert \leq \func{Re}F_{k}\left(
x_{i}\right) ,\qquad i=1,\dots ,n,\ k=1,\dots ,m.  \label{1.1.7}
\end{equation}%
Then one has the following reverse of the generalised triangle inequality 
\cite{DM}%
\begin{equation}
\left( \frac{\sum_{k=1}^{m}r_{k}^{2}}{c}\right) ^{\frac{1}{2}%
}\sum_{i=1}^{n}\left\Vert x_{i}\right\Vert \leq \left\Vert
\sum_{i=1}^{n}x_{i}\right\Vert ,  \label{1.1.8}
\end{equation}%
where equality holds if and only if both%
\begin{equation}
F_{k}\left( \sum_{i=1}^{n}x_{i}\right) =r_{k}\sum_{i=1}^{n}\left\Vert
x_{i}\right\Vert ,\qquad k=1,\dots ,m  \label{1.1.9}
\end{equation}%
and 
\begin{equation}
\sum_{k=1}^{m}\left[ F_{k}\left( \sum_{i=1}^{n}x_{i}\right) \right]
^{2}=c\left\Vert \sum_{i=1}^{n}x_{i}\right\Vert ^{2}.  \label{1.1.10}
\end{equation}
\end{theorem}

If $X=H,$ an inner product space, then, for $F_{k}\left( x\right)
=\left\langle x,e_{k}\right\rangle ,$ where $\left\{ e_{k}\right\} _{k=%
\overline{1,n}}$ is an orthonormal family in $H,$ i.e., $\left\langle
e_{i},e_{j}\right\rangle =\delta _{ij},$ $i,j\in \left\{ 1,\dots ,k\right\}
, $ $\delta _{ij}$ is Kronecker delta, the condition (\ref{1.1.7}) may be
replaced by%
\begin{equation}
0\leq r_{k}\left\Vert x_{i}\right\Vert \leq \func{Re}\left\langle
x_{i},e_{k}\right\rangle ,\qquad i=1,\dots ,n,\ k=1,\dots ,m;  \label{1.1.11}
\end{equation}%
implying the following reverse of the generalised triangle inequality%
\begin{equation}
\left( \sum_{k=1}^{m}r_{k}^{2}\right) ^{\frac{1}{2}}\sum_{i=1}^{n}\left\Vert
x_{i}\right\Vert \leq \left\Vert \sum_{i=1}^{n}x_{i}\right\Vert ,
\label{1.1.12}
\end{equation}%
where the equality holds if and only if%
\begin{equation}
\sum_{i=1}^{n}x_{i}=\left( \sum_{i=1}^{n}\left\Vert x_{i}\right\Vert \right)
\sum_{k=1}^{m}r_{k}e_{k}.  \label{1.1.13}
\end{equation}

The aim of the following sections is to present recent reverses of the
triangle inequality obtained by the author in \cite{DRA1} and \cite{DRA2}.
New results are established for the general case of normed spaces. Their
versions in inner product spaces are analyzed and applications for complex
numbers are given as well.

For various classical inequalities related to the triangle inequality, see
Chapter XVII of the book \cite{MPF} and the references therein.

\section{ Inequalities of Diaz-Metcalf Type for $m$ Functionals}

\subsection{The Case of Normed Spaces}

The following result may be stated \cite{DRA1}.

\begin{theorem}[Dragomir, 2004]
\label{t1.2.1}Let $\left( X,\left\Vert \cdot \right\Vert \right) $ be a
normed linear space over the real or complex number field $\mathbb{K}$ and $%
F_{k}:X\rightarrow \mathbb{K}$, $k\in \left\{ 1,\dots ,m\right\} $
continuous linear functionals on $X.$ If $x_{i}\in X\backslash \left\{
0\right\} ,$ $i\in \left\{ 1,\dots ,n\right\} $ are such that there exists
the constants $r_{k}\geq 0$, $k\in \left\{ 1,\dots ,m\right\} $ with $%
\sum_{k=1}^{m}r_{k}>0$ and%
\begin{equation}
\func{Re}F_{k}\left( x_{i}\right) \geq r_{k}\left\Vert x_{i}\right\Vert 
\text{ }  \label{1.2.1}
\end{equation}%
\ for each \ $i\in \left\{ 1,\dots ,n\right\} $ and $k\in \left\{ 1,\dots
,m\right\} ,$ then%
\begin{equation}
\sum_{i=1}^{n}\left\Vert x_{i}\right\Vert \leq \frac{\left\Vert
\sum_{k=1}^{m}F_{k}\right\Vert }{\sum_{k=1}^{m}r_{k}}\left\Vert
\sum_{i=1}^{n}x_{i}\right\Vert .  \label{1.2.2}
\end{equation}%
The case of equality holds in (\ref{1.2.2}) if both%
\begin{equation}
\left( \sum_{k=1}^{m}F_{k}\right) \left( \sum_{i=1}^{n}x_{i}\right) =\left(
\sum_{k=1}^{m}r_{k}\right) \sum_{i=1}^{n}\left\Vert x_{i}\right\Vert
\label{1.2.3}
\end{equation}%
and%
\begin{equation}
\left( \sum_{k=1}^{m}F_{k}\right) \left( \sum_{i=1}^{n}x_{i}\right)
=\left\Vert \sum_{k=1}^{m}F_{k}\right\Vert \left\Vert
\sum_{i=1}^{n}x_{i}\right\Vert .  \label{1.2.4}
\end{equation}
\end{theorem}

\begin{proof}
Utilising the hypothesis (\ref{1.2.1}) and the properties of the modulus, we
have%
\begin{align}
I& :=\left\vert \left( \sum_{k=1}^{m}F_{k}\right) \left(
\sum_{i=1}^{n}x_{i}\right) \right\vert \geq \left\vert \func{Re}\left[
\left( \sum_{k=1}^{m}F_{k}\right) \left( \sum_{i=1}^{n}x_{i}\right) \right]
\right\vert  \label{1.2.5} \\
& \geq \sum_{k=1}^{m}\func{Re}F_{k}\left( \sum_{i=1}^{n}x_{i}\right)
=\sum_{k=1}^{m}\sum_{i=1}^{n}\func{Re}F_{k}\left( x_{i}\right)  \notag \\
& \geq \left( \sum_{k=1}^{m}r_{k}\right) \sum_{i=1}^{n}\left\Vert
x_{i}\right\Vert .  \notag
\end{align}%
On the other hand, by the continuity property of $F_{k},$ $k\in \left\{
1,\dots ,m\right\} $ we obviously have%
\begin{equation}
I=\left\vert \left( \sum_{k=1}^{m}F_{k}\right) \left(
\sum_{i=1}^{n}x_{i}\right) \right\vert \leq \left\Vert
\sum_{k=1}^{m}F_{k}\right\Vert \left\Vert \sum_{i=1}^{n}x_{i}\right\Vert .
\label{1.2.6}
\end{equation}%
Making use of (\ref{1.2.5}) and (\ref{1.2.6}), we deduce the desired
inequality (\ref{1.2.2}).

Now, if (\ref{1.2.3}) and (\ref{1.2.4}) are valid, then, obviously, the case
of equality holds true in the inequality (\ref{1.2.2}).

Conversely, if the case of equality holds in (\ref{1.2.2}), then it must
hold in all the inequalities used to prove (\ref{1.2.2}). Therefore we have%
\begin{equation}
\func{Re}F_{k}\left( x_{i}\right) =r_{k}\left\Vert x_{i}\right\Vert \text{ }
\label{1.2.7}
\end{equation}%
\ for each \ $i\in \left\{ 1,\dots ,n\right\} $, $k\in \left\{ 1,\dots
,m\right\} ;$%
\begin{equation}
\sum_{k=1}^{m}\func{Im}F_{k}\left( \sum_{i=1}^{n}x_{i}\right) =0
\label{1.2.8}
\end{equation}%
and 
\begin{equation}
\sum_{k=1}^{m}\func{Re}F_{k}\left( \sum_{i=1}^{n}x_{i}\right) =\left\Vert
\sum_{k=1}^{m}F_{k}\right\Vert \left\Vert \sum_{i=1}^{n}x_{i}\right\Vert .
\label{1.2.9}
\end{equation}%
Note that, from (\ref{1.2.7}), by summation over $i$ and $k,$ we get%
\begin{equation}
\func{Re}\left[ \left( \sum_{k=1}^{m}F_{k}\right) \left(
\sum_{i=1}^{n}x_{i}\right) \right] =\left( \sum_{k=1}^{m}r_{k}\right)
\sum_{i=1}^{n}\left\Vert x_{i}\right\Vert .  \label{1.2.10}
\end{equation}%
Since (\ref{1.2.8}) and (\ref{1.2.10}) imply (\ref{1.2.3}), while (\ref%
{1.2.9}) and (\ref{1.2.10}) imply (\ref{1.2.4}) hence the theorem is proved.
\end{proof}

\begin{remark}
\label{r1.2.2}If the norms $\left\Vert F_{k}\right\Vert ,$ $k\in \left\{
1,\dots ,m\right\} $ are easier to find, then, from (\ref{1.2.2}), one may
get the (coarser) inequality that might be more useful in practice:%
\begin{equation}
\sum_{i=1}^{n}\left\Vert x_{i}\right\Vert \leq \frac{\sum_{k=1}^{m}\left%
\Vert F_{k}\right\Vert }{\sum_{k=1}^{m}r_{k}}\left\Vert
\sum_{i=1}^{n}x_{i}\right\Vert .  \label{1.2.11}
\end{equation}
\end{remark}

\subsection{The Case of Inner Product Spaces}

The case of inner product spaces, in which we may provide a simpler
condition for equality, is of interest in applications \cite{DRA1}.

\begin{theorem}[Dragomir, 2004]
\label{t1.2.3}Let $\left( H;\left\langle \cdot ,\cdot \right\rangle \right) $
be an inner product space over the real or complex number field $\mathbb{K}$%
, $e_{k},$ $x_{i}\in H\backslash \left\{ 0\right\} $, $k\in \left\{ 1,\dots
,m\right\} ,$ $i\in \left\{ 1,\dots ,n\right\} .$ If $r_{k}\geq 0,$ $k\in
\left\{ 1,\dots ,m\right\} $ with $\sum_{k=1}^{m}r_{k}>0$ satisfy%
\begin{equation}
\func{Re}\left\langle x_{i},e_{k}\right\rangle \geq r_{k}\left\Vert
x_{i}\right\Vert \text{ \ }  \label{1.2.12}
\end{equation}%
for each \ $i\in \left\{ 1,\dots ,n\right\} $ and $k\in \left\{ 1,\dots
,m\right\} ,$ then%
\begin{equation}
\sum_{i=1}^{n}\left\Vert x_{i}\right\Vert \leq \frac{\left\Vert
\sum_{k=1}^{m}e_{k}\right\Vert }{\sum_{k=1}^{m}r_{k}}\left\Vert
\sum_{i=1}^{n}x_{i}\right\Vert .  \label{1.2.13}
\end{equation}%
The case of equality holds in (\ref{1.2.13}) if and only if%
\begin{equation}
\sum_{i=1}^{n}x_{i}=\frac{\sum_{k=1}^{m}r_{k}}{\left\Vert
\sum_{k=1}^{m}e_{k}\right\Vert ^{2}}\left( \sum_{i=1}^{n}\left\Vert
x_{i}\right\Vert \right) \sum_{k=1}^{m}e_{k}.  \label{1.2.14}
\end{equation}
\end{theorem}

\begin{proof}
By the properties of inner product and by (\ref{1.2.12}), we have%
\begin{align}
& \left\vert \left\langle
\sum_{i=1}^{n}x_{i},\sum_{k=1}^{m}e_{k}\right\rangle \right\vert
\label{1.2.15} \\
& \geq \left\vert \sum_{k=1}^{m}\func{Re}\left\langle
\sum_{i=1}^{n}x_{i},e_{k}\right\rangle \right\vert \geq \sum_{k=1}^{m}\func{%
Re}\left\langle \sum_{i=1}^{n}x_{i},e_{k}\right\rangle  \notag \\
& =\sum_{k=1}^{m}\sum_{i=1}^{n}\func{Re}\left\langle
x_{i},e_{k}\right\rangle \geq \left( \sum_{k=1}^{m}r_{k}\right)
\sum_{i=1}^{n}\left\Vert x_{i}\right\Vert >0.  \notag
\end{align}%
Observe also that, by (\ref{1.2.15}), $\sum_{k=1}^{m}e_{k}\neq 0.$

On utlising Schwarz's inequality in the inner product space $\left(
H;\left\langle \cdot ,\cdot \right\rangle \right) $ for $%
\sum_{i=1}^{n}x_{i}, $ $\sum_{k=1}^{m}e_{k},$ we have%
\begin{equation}
\left\Vert \sum_{i=1}^{n}x_{i}\right\Vert \left\Vert
\sum_{k=1}^{m}e_{k}\right\Vert \geq \left\vert \left\langle
\sum_{i=1}^{n}x_{i},\sum_{k=1}^{m}e_{k}\right\rangle \right\vert .
\label{1.2.16}
\end{equation}%
Making use of (\ref{1.2.15}) and (\ref{1.2.16}), we can conclude that (\ref%
{1.2.13}) holds.

Now, if (\ref{1.2.14}) holds true, then, by taking the norm, we have%
\begin{align*}
\left\Vert \sum_{i=1}^{n}x_{i}\right\Vert & =\frac{\left(
\sum_{k=1}^{m}r_{k}\right) \sum_{i=1}^{n}\left\Vert x_{i}\right\Vert }{%
\left\Vert \sum_{k=1}^{m}e_{k}\right\Vert ^{2}}\left\Vert
\sum_{k=1}^{m}e_{k}\right\Vert  \\
& =\frac{\left( \sum_{k=1}^{m}r_{k}\right) }{\left\Vert
\sum_{k=1}^{m}e_{k}\right\Vert }\sum_{i=1}^{n}\left\Vert x_{i}\right\Vert ,
\end{align*}%
i.e., the case of equality holds in (\ref{1.2.13}).

Conversely, if the case of equality holds in (\ref{1.2.13}), then it must
hold in all the inequalities used to prove (\ref{1.2.13}). Therefore, we have%
\begin{equation}
\func{Re}\left\langle x_{i},e_{k}\right\rangle =r_{k}\left\Vert
x_{i}\right\Vert \text{ \ }  \label{1.2.17}
\end{equation}%
for each \ $i\in \left\{ 1,\dots ,n\right\} $ and $k\in \left\{ 1,\dots
,m\right\} ,$%
\begin{equation}
\left\Vert \sum_{i=1}^{n}x_{i}\right\Vert \left\Vert
\sum_{k=1}^{m}e_{k}\right\Vert =\left\vert \left\langle
\sum_{i=1}^{n}x_{i},\sum_{k=1}^{m}e_{k}\right\rangle \right\vert 
\label{1.2.18}
\end{equation}%
and%
\begin{equation}
\func{Im}\left\langle \sum_{i=1}^{n}x_{i},\sum_{k=1}^{m}e_{k}\right\rangle
=0.  \label{1.2.19}
\end{equation}%
From (\ref{1.2.17}), on summing over $i$ and $k,$ we get%
\begin{equation}
\func{Re}\left\langle \sum_{i=1}^{n}x_{i},\sum_{k=1}^{m}e_{k}\right\rangle
=\left( \sum_{k=1}^{m}r_{k}\right) \sum_{i=1}^{n}\left\Vert x_{i}\right\Vert
.  \label{1.2.20}
\end{equation}%
By (\ref{1.2.19}) and (\ref{1.2.20}), we have%
\begin{equation}
\left\langle \sum_{i=1}^{n}x_{i},\sum_{k=1}^{m}e_{k}\right\rangle =\left(
\sum_{k=1}^{m}r_{k}\right) \sum_{i=1}^{n}\left\Vert x_{i}\right\Vert .
\label{1.2.21}
\end{equation}%
On the other hand, by the use of the following identity in inner product
spaces%
\begin{equation}
\left\Vert u-\frac{\left\langle u,v\right\rangle v}{\left\Vert v\right\Vert
^{2}}\right\Vert ^{2}=\frac{\left\Vert u\right\Vert ^{2}\left\Vert
v\right\Vert ^{2}-\left\vert \left\langle u,v\right\rangle \right\vert ^{2}}{%
\left\Vert v\right\Vert ^{2}},\quad v\neq 0,  \label{Id}
\end{equation}%
the relation (\ref{1.2.18}) holds if and only if%
\begin{equation}
\sum_{i=1}^{n}x_{i}=\frac{\left\langle
\sum_{i=1}^{n}x_{i},\sum_{k=1}^{m}e_{k}\right\rangle }{\left\Vert
\sum_{k=1}^{m}e_{k}\right\Vert ^{2}}\sum_{k=1}^{m}e_{k}.  \label{1.2.22}
\end{equation}%
Finally, on utilising (\ref{1.2.21}) and (\ref{1.2.22}), we deduce that the
condition (\ref{1.2.14}) is necessary for the equality case in (\ref{1.2.13}%
).
\end{proof}

Before we give a corollary of the above theorem, we need to state the
following lemma that has been basically obtained in \cite{SSD2}. For the
sake of completeness, we provide a short proof here as well.

\begin{lemma}[Dragomir, 2004]
\label{l1.2.4}Let $\left( H;\left\langle \cdot ,\cdot \right\rangle \right) $
be an inner product space over the real or complex number field $\mathbb{K}$
and $x,a\in H,$ $r>0$ such that:%
\begin{equation}
\left\Vert x-a\right\Vert \leq r<\left\Vert a\right\Vert .  \label{1.2.23}
\end{equation}%
Then we have the inequality%
\begin{equation}
\left\Vert x\right\Vert \left( \left\Vert a\right\Vert ^{2}-r^{2}\right) ^{%
\frac{1}{2}}\leq \func{Re}\left\langle x,a\right\rangle  \label{1.2.24}
\end{equation}%
or, equivalently%
\begin{equation}
\left\Vert x\right\Vert ^{2}\left\Vert a\right\Vert ^{2}-\left[ \func{Re}%
\left\langle x,a\right\rangle \right] ^{2}\leq r^{2}\left\Vert x\right\Vert
^{2}.  \label{1.2.25}
\end{equation}%
The case of equality holds in (\ref{1.2.24}) (or in (\ref{1.2.25})) if and
only if%
\begin{equation}
\left\Vert x-a\right\Vert =r\text{ \ and \ }\left\Vert x\right\Vert
^{2}+r^{2}=\left\Vert a\right\Vert ^{2}.  \label{1.2.26}
\end{equation}
\end{lemma}

\begin{proof}
From the first part of (\ref{1.2.23}), we have%
\begin{equation}
\left\Vert x\right\Vert ^{2}+\left\Vert a\right\Vert ^{2}-r^{2}\leq 2\func{Re%
}\left\langle x,a\right\rangle .  \label{1.2.27}
\end{equation}%
By the second part of (\ref{1.2.23}) we have $\left( \left\Vert a\right\Vert
^{2}-r^{2}\right) ^{\frac{1}{2}}>0,$ therefore, by (\ref{1.2.27}), we may
state that%
\begin{equation}
0<\frac{\left\Vert x\right\Vert ^{2}}{\left( \left\Vert a\right\Vert
^{2}-r^{2}\right) ^{\frac{1}{2}}}+\left( \left\Vert a\right\Vert
^{2}-r^{2}\right) ^{\frac{1}{2}}\leq \frac{2\func{Re}\left\langle
x,a\right\rangle }{\left( \left\Vert a\right\Vert ^{2}-r^{2}\right) ^{\frac{1%
}{2}}}.  \label{1.2.28}
\end{equation}%
Utilising the elementary inequality%
\begin{equation*}
\frac{1}{\alpha }q+\alpha p\geq 2\sqrt{pq},\quad \alpha >0,\ p>0,\ q\geq 0;
\end{equation*}%
with equality if and only if $\alpha =\sqrt{\frac{q}{p}},$ we may state (for 
$\alpha =\left( \left\Vert a\right\Vert ^{2}-r^{2}\right) ^{\frac{1}{2}},$ $%
p=1,$ $q=\left\Vert x\right\Vert ^{2}$) that%
\begin{equation}
2\left\Vert x\right\Vert \leq \frac{\left\Vert x\right\Vert ^{2}}{\left(
\left\Vert a\right\Vert ^{2}-r^{2}\right) ^{\frac{1}{2}}}+\left( \left\Vert
a\right\Vert ^{2}-r^{2}\right) ^{\frac{1}{2}}.  \label{1.2.29}
\end{equation}%
The inequality (\ref{1.2.24}) follows now by (\ref{1.2.28}) and (\ref{1.2.29}%
).

From the above argument, it is clear that the equality holds in (\ref{1.2.24}%
) if and only if it holds in (\ref{1.2.28}) and (\ref{1.2.29}). However, the
equality holds in (\ref{1.2.28}) if and only if $\left\Vert x-a\right\Vert
=r $ and in (\ref{1.2.29}) if and only if $\left( \left\Vert a\right\Vert
^{2}-r^{2}\right) ^{\frac{1}{2}}=\left\Vert x\right\Vert .$

The proof is thus completed.
\end{proof}

We may now state the following corollary \cite{DRA1}.

\begin{corollary}
\label{c1.2.5}Let $\left( H;\left\langle \cdot ,\cdot \right\rangle \right) $
be an inner product space over the real or complex number field $\mathbb{K}$%
, $e_{k},$ $x_{i}\in H\backslash \left\{ 0\right\} $, $k\in \left\{ 1,\dots
,m\right\} ,$ $i\in \left\{ 1,\dots ,n\right\} .$ If $\rho _{k}\geq 0,$ $%
k\in \left\{ 1,\dots ,m\right\} $ with%
\begin{equation}
\left\Vert x_{i}-e_{k}\right\Vert \leq \rho _{k}<\left\Vert e_{k}\right\Vert 
\text{\ }  \label{1.2.30}
\end{equation}%
for each \ $i\in \left\{ 1,\dots ,n\right\} $ and $k\in \left\{ 1,\dots
,m\right\} ,$ then%
\begin{equation}
\sum_{i=1}^{n}\left\Vert x_{i}\right\Vert \leq \frac{\left\Vert
\sum_{k=1}^{m}e_{k}\right\Vert }{\sum_{k=1}^{m}\left( \left\Vert
e_{k}\right\Vert ^{2}-\rho _{k}^{2}\right) ^{\frac{1}{2}}}\left\Vert
\sum_{i=1}^{n}x_{i}\right\Vert .  \label{1.2.31}
\end{equation}%
The case of equality holds in (\ref{1.2.31}) if and only if%
\begin{equation*}
\sum_{i=1}^{n}x_{i}=\frac{\sum_{k=1}^{m}\left( \left\Vert e_{k}\right\Vert
^{2}-\rho _{k}^{2}\right) ^{\frac{1}{2}}}{\left\Vert
\sum_{k=1}^{m}e_{k}\right\Vert ^{2}}\left( \sum_{i=1}^{n}\left\Vert
x_{i}\right\Vert \right) \sum_{k=1}^{m}e_{k}.
\end{equation*}
\end{corollary}

\begin{proof}
Utilising Lemma \ref{l1.2.4}, we have from (\ref{1.2.30}) that%
\begin{equation*}
\left\Vert x_{i}\right\Vert \left( \left\Vert e_{k}\right\Vert ^{2}-\rho
_{k}^{2}\right) ^{\frac{1}{2}}\leq \func{Re}\left\langle
x_{i},e_{k}\right\rangle
\end{equation*}%
for each $k\in \left\{ 1,\dots ,m\right\} $ and $i\in \left\{ 1,\dots
,n\right\} .$

Applying Theorem \ref{t1.2.3} for 
\begin{equation*}
r_{k}:=\left( \left\Vert e_{k}\right\Vert ^{2}-\rho _{k}^{2}\right) ^{\frac{1%
}{2}},\quad k\in \left\{ 1,\dots ,m\right\} ,
\end{equation*}%
we deduce the desired result.
\end{proof}

\begin{remark}
\label{r1.2.6}If $\left\{ e_{k}\right\} _{k\in \left\{ 1,\dots ,m\right\} }$
are orthogonal, then (\ref{1.2.31}) becomes%
\begin{equation}
\sum_{i=1}^{n}\left\Vert x_{i}\right\Vert \leq \frac{\left(
\sum_{k=1}^{m}\left\Vert e_{k}\right\Vert ^{2}\right) ^{\frac{1}{2}}}{%
\sum_{k=1}^{m}\left( \left\Vert e_{k}\right\Vert ^{2}-\rho _{k}^{2}\right) ^{%
\frac{1}{2}}}\left\Vert \sum_{i=1}^{n}x_{i}\right\Vert  \label{1.2.32}
\end{equation}%
with equality if and only if%
\begin{equation*}
\sum_{i=1}^{n}x_{i}=\frac{\sum_{k=1}^{m}\left( \left\Vert e_{k}\right\Vert
^{2}-\rho _{k}^{2}\right) ^{\frac{1}{2}}}{\sum_{k=1}^{m}\left\Vert
e_{k}\right\Vert ^{2}}\left( \sum_{i=1}^{n}\left\Vert x_{i}\right\Vert
\right) \sum_{k=1}^{m}e_{k}.
\end{equation*}%
Moreover, if $\left\{ e_{k}\right\} _{k\in \left\{ 1,\dots ,m\right\} }$ is
assumed to be orthonormal and%
\begin{equation*}
\left\Vert x_{i}-e_{k}\right\Vert \leq \rho _{k}\text{ \ for }k\in \left\{
1,\dots ,m\right\} ,\text{\ }i\in \left\{ 1,\dots ,n\right\}
\end{equation*}%
where $\rho _{k}\in \lbrack 0,1)$ for $k\in \left\{ 1,\dots ,m\right\} ,$
then%
\begin{equation}
\sum_{i=1}^{n}\left\Vert x_{i}\right\Vert \leq \frac{\sqrt{m}}{%
\sum_{k=1}^{m}\left( 1-\rho _{k}^{2}\right) ^{\frac{1}{2}}}\left\Vert
\sum_{i=1}^{n}x_{i}\right\Vert  \label{1.2.33}
\end{equation}%
with equality if and only if%
\begin{equation*}
\sum_{i=1}^{n}x_{i}=\frac{\sum_{k=1}^{m}\left( 1-\rho _{k}^{2}\right) ^{%
\frac{1}{2}}}{m}\left( \sum_{i=1}^{n}\left\Vert x_{i}\right\Vert \right)
\sum_{k=1}^{m}e_{k}.
\end{equation*}
\end{remark}

The following lemma may be stated as well \cite{SSD1}.

\begin{lemma}[Dragomir, 2004]
\label{l1.2.7}Let $\left( H;\left\langle \cdot ,\cdot \right\rangle \right) $
be an inner product space over the real or complex number field $\mathbb{K}$%
, $x,y\in H$ and $M\geq m>0.$ If%
\begin{equation}
\func{Re}\left\langle My-x,x-my\right\rangle \geq 0  \label{1.2.34}
\end{equation}%
or, equivalently,%
\begin{equation}
\left\Vert x-\frac{m+M}{2}y\right\Vert \leq \frac{1}{2}\left( M-m\right)
\left\Vert y\right\Vert ,  \label{1.2.35}
\end{equation}%
then%
\begin{equation}
\left\Vert x\right\Vert \left\Vert y\right\Vert \leq \frac{1}{2}\cdot \frac{%
M+m}{\sqrt{mM}}\func{Re}\left\langle x,y\right\rangle .  \label{1.2.36}
\end{equation}%
The equality holds in (\ref{1.2.36}) if and only if the case of equality
holds in (\ref{1.2.34}) and%
\begin{equation}
\left\Vert x\right\Vert =\sqrt{mM}\left\Vert y\right\Vert .  \label{1.2.37}
\end{equation}
\end{lemma}

\begin{proof}
Obviously,%
\begin{equation*}
\func{Re}\left\langle My-x,x-my\right\rangle =\left( M+m\right) \func{Re}%
\left\langle x,y\right\rangle -\left\Vert x\right\Vert ^{2}-mM\left\Vert
y\right\Vert ^{2}.
\end{equation*}%
Then (\ref{1.2.34}) is clearly equivalent to%
\begin{equation}
\frac{\left\Vert x\right\Vert ^{2}}{\sqrt{mM}}+\sqrt{mM}\left\Vert
y\right\Vert ^{2}\leq \frac{M+m}{\sqrt{mM}}\func{Re}\left\langle
x,y\right\rangle .  \label{1.2.38}
\end{equation}%
Since, obviously,%
\begin{equation}
2\left\Vert x\right\Vert \left\Vert y\right\Vert \leq \frac{\left\Vert
x\right\Vert ^{2}}{\sqrt{mM}}+\sqrt{mM}\left\Vert y\right\Vert ^{2},
\label{1.2.39}
\end{equation}%
with equality iff $\left\Vert x\right\Vert =\sqrt{mM}\left\Vert y\right\Vert
,$ hence (\ref{1.2.38}) and (\ref{1.2.39}) imply (\ref{1.2.36}).

The case of equality is obvious and we omit the details.
\end{proof}

Finally, we may state the following corollary of Theorem \ref{t1.2.3}, see 
\cite{DRA1}.

\begin{corollary}
\label{c1.2.8}Let $\left( H;\left\langle \cdot ,\cdot \right\rangle \right) $
be an inner product space over the real or complex number field $\mathbb{K}$%
, $e_{k},$ $x_{i}\in H\backslash \left\{ 0\right\} $, $k\in \left\{ 1,\dots
,m\right\} ,$ $i\in \left\{ 1,\dots ,n\right\} .$ If $M_{k}>\mu _{k}>0,$ $%
k\in \left\{ 1,\dots ,m\right\} $ are such that either%
\begin{equation}
\func{Re}\left\langle M_{k}e_{k}-x_{i},x_{i}-\mu _{k}e_{k}\right\rangle \geq
0  \label{1.2.40}
\end{equation}%
or, equivalently,%
\begin{equation*}
\left\Vert x_{i}-\frac{M_{k}+\mu _{k}}{2}e_{k}\right\Vert \leq \frac{1}{2}%
\left( M_{k}-\mu _{k}\right) \left\Vert e_{k}\right\Vert
\end{equation*}%
for each $k\in \left\{ 1,\dots ,m\right\} $ and $i\in \left\{ 1,\dots
,n\right\} ,$ then%
\begin{equation}
\sum_{i=1}^{n}\left\Vert x_{i}\right\Vert \leq \frac{\left\Vert
\sum_{k=1}^{m}e_{k}\right\Vert }{\sum_{k=1}^{m}\frac{2\cdot \sqrt{\mu
_{k}M_{k}}}{\mu _{k}+M_{k}}\left\Vert e_{k}\right\Vert }\left\Vert
\sum_{i=1}^{n}x_{i}\right\Vert .  \label{1.2.41}
\end{equation}%
The case of equality holds in (\ref{1.2.41}) if and only if%
\begin{equation*}
\sum_{i=1}^{n}x_{i}=\frac{\sum_{k=1}^{m}\frac{2\cdot \sqrt{\mu _{k}M_{k}}}{%
\mu _{k}+M_{k}}\left\Vert e_{k}\right\Vert }{\left\Vert
\sum_{k=1}^{m}e_{k}\right\Vert ^{2}}\sum_{i=1}^{n}\left\Vert
x_{i}\right\Vert \sum_{k=1}^{m}e_{k}.
\end{equation*}
\end{corollary}

\begin{proof}
Utilising Lemma \ref{l1.2.7}, by (\ref{1.2.40}) we deduce%
\begin{equation*}
\frac{2\cdot \sqrt{\mu _{k}M_{k}}}{\mu _{k}+M_{k}}\left\Vert
x_{i}\right\Vert \left\Vert e_{k}\right\Vert \leq \func{Re}\left\langle
x_{i},e_{k}\right\rangle
\end{equation*}%
for each $k\in \left\{ 1,\dots ,m\right\} $ and $i\in \left\{ 1,\dots
,n\right\} .$

Applying Theorem \ref{t1.2.3} for%
\begin{equation*}
r_{k}:=\frac{2\cdot \sqrt{\mu _{k}M_{k}}}{\mu _{k}+M_{k}}\left\Vert
e_{k}\right\Vert ,\quad k\in \left\{ 1,\dots ,m\right\} ,
\end{equation*}%
we deduce the desired result.
\end{proof}

\section{Diaz-Metcalf Inequality for Semi-Inner Products}

In 1961, G. Lumer \cite{L} introduced the following concept.

\begin{definition}
\label{d2.3.1}Let $X$ be a linear space over the real or complex number
field $\mathbb{K}$. The mapping $\left[ \cdot ,\cdot \right] :X\times
X\rightarrow \mathbb{K}$ is called a \textit{semi-inner product }on $X,$ if
the following properties are satisfied (see also \cite[p. 17]{SSD1}):

\begin{enumerate}
\item[$\left( i\right) $] $\left[ x+y,z\right] =\left[ x,z\right] +\left[ y,z%
\right] $ for all $x,y,z\in X;$

\item[$\left( ii\right) $] $\left[ \lambda x,y\right] =\lambda \left[ x,y%
\right] $ for all $x,y\in X$ and $\lambda \in \mathbb{K}$;

\item[$\left( iii\right) $] $\left[ x,x\right] \geq 0$ for all $x\in X$ and $%
\left[ x,x\right] =0$ implies $x=0$;

\item[$\left( iv\right) $] $\left\vert \left[ x,y\right] \right\vert
^{2}\leq \left[ x,x\right] \left[ y,y\right] $ for all $x,y\in X;$

\item[$\left( v\right) $] $\left[ x,\lambda y\right] =\bar{\lambda}\left[ x,y%
\right] $ for all $x,y\in X$ and $\lambda \in \mathbb{K}$.
\end{enumerate}
\end{definition}

It is well known that the mapping $X\ni x\longmapsto \left[ x,x\right] ^{%
\frac{1}{2}}\in \mathbb{R}$ is a norm on $X$ and for any $y\in X,$ the
functional $X\ni x\overset{\varphi _{y}}{\longmapsto }\left[ x,y\right] \in 
\mathbb{K}$ is a continuous linear functional on $X$ endowed with the norm $%
\left\Vert \cdot \right\Vert $ generated by $\left[ \cdot ,\cdot \right] .$
Moreover, one has $\left\Vert \varphi _{y}\right\Vert =\left\Vert
y\right\Vert $ (see for instance \cite[p. 17]{SSD1}).

Let $\left( X,\left\Vert \cdot \right\Vert \right) $ be a real or complex
normed space. If $J:X\rightarrow _{2}X^{\ast }$ is the \textit{normalised
duality mapping }defined on $X,$ i.e., we recall that (see for instance \cite%
[p. 1]{SSD1})%
\begin{equation*}
J\left( x\right) =\left\{ \varphi \in X^{\ast }|\varphi \left( x\right)
=\left\Vert \varphi \right\Vert \left\Vert x\right\Vert ,\ \left\Vert
\varphi \right\Vert =\left\Vert x\right\Vert \right\} ,\ \ \ x\in X,
\end{equation*}%
then we may state the following representation result (see for instance \cite%
[p. 18]{SSD1}):

\textit{Each semi-inner product }$\left[ \cdot ,\cdot \right] :X\times
X\rightarrow K$\textit{\ that generates the norm }$\left\Vert \cdot
\right\Vert $\textit{\ of the normed linear space }$\left( X,\left\Vert
\cdot \right\Vert \right) $\textit{\ over the real or complex number field }$%
K$\textit{, is of the form}%
\begin{equation*}
\left[ x,y\right] =\left\langle \tilde{J}\left( y\right) ,x\right\rangle 
\text{ \ for any \ }x,y\in X,
\end{equation*}%
\textit{where }$\tilde{J}$\textit{\ is a selection of the normalised duality
mapping and }$\left\langle \varphi ,x\right\rangle :=\varphi \left( x\right) 
$\textit{\ for }$\varphi \in X^{\ast }$\textit{\ and }$x\in X.$

Utilising the concept of semi-inner products, we can state the following
particular case of the Diaz-Metcalf inequality.

\begin{corollary}
\label{c2.3.1}Let $\left( X,\left\Vert \cdot \right\Vert \right) $ be a
normed linear space, $\left[ \cdot ,\cdot \right] :X\times X\rightarrow 
\mathbb{K}$ a semi-inner product generating the norm $\left\Vert \cdot
\right\Vert $ and $e\in X,$ $\left\Vert e\right\Vert =1.$ If $x_{i}\in X,$ $%
i\in \left\{ 1,\dots ,n\right\} $ and $r\geq 0$ such that%
\begin{equation}
r\left\Vert x_{i}\right\Vert \leq \func{Re}\left[ x_{i},e\right] \text{ \
for each \ }i\in \left\{ 1,\dots ,n\right\} ,  \label{2.3.1}
\end{equation}%
then we have the inequality%
\begin{equation}
r\sum_{i=1}^{n}\left\Vert x_{i}\right\Vert \leq \left\Vert
\sum_{i=1}^{n}x_{i}\right\Vert .  \label{2.3.2}
\end{equation}%
The case of equality holds in (\ref{2.3.2}) if and only if both%
\begin{equation}
\left[ \sum_{i=1}^{n}x_{i},e\right] =r\sum_{i=1}^{n}\left\Vert
x_{i}\right\Vert  \label{2.3.3}
\end{equation}%
and 
\begin{equation}
\left[ \sum_{i=1}^{n}x_{i},e\right] =\left\Vert
\sum_{i=1}^{n}x_{i}\right\Vert .  \label{2.3.4}
\end{equation}
\end{corollary}

The proof is obvious from the Diaz-Metcalf theorem \cite[Theorem 3]{DM}
applied for the continuous linear functional $F_{e}\left( x\right) =\left[
x,e\right] ,$ $x\in X.$

Before we provide a simpler necessary and sufficient condition of equality
in (\ref{2.3.2}), we need to recall the concept of strictly convex normed
spaces and a classical characterisation of these spaces.

\begin{definition}
\label{d2.3.2}A normed linear space $\left( X,\left\Vert \cdot \right\Vert
\right) $ is said to be strictly convex if for every $x,y$ from $X$ with $%
x\neq y$ and $\left\Vert x\right\Vert =\left\Vert y\right\Vert =1,$ we have $%
\left\Vert \lambda x+\left( 1-\lambda \right) y\right\Vert <1$ for all $%
\lambda \in \left( 0,1\right) .$
\end{definition}

The following characterisation of strictly convex spaces is useful in what
follows (see \cite{B}, \cite{GS}, or \cite[p. 21]{SSD1}).

\begin{theorem}
\label{t2.3.1}Let $\left( X,\left\Vert \cdot \right\Vert \right) $ be a
normed linear space over $\mathbb{K}$ and $\left[ \cdot ,\cdot \right] $ a
semi-inner product generating its norm. The following statements are
equivalent:

\begin{enumerate}
\item[$\left( i\right) $] $\left( X,\left\Vert \cdot \right\Vert \right) $
is strictly convex;

\item[$\left( ii\right) $] For every $x,y\in X,$ $x,y\neq 0$ with $\left[ x,y%
\right] =\left\Vert x\right\Vert \left\Vert y\right\Vert ,$ there exists a $%
\lambda >0$ such that $x=\lambda y.$
\end{enumerate}
\end{theorem}

The following result may be stated.

\begin{corollary}
\label{c2.3.2}Let $\left( X,\left\Vert \cdot \right\Vert \right) $ be a
strictly convex normed linear space, $\left[ \cdot ,\cdot \right] $ a
semi-inner product generating the norm and $e,$ $x_{i}$ $\left( i\in \left\{
1,\dots ,n\right\} \right) $ as in Corollary \ref{c2.3.1}. Then the case of
equality holds in (\ref{2.3.2}) if and only if%
\begin{equation}
\sum_{i=1}^{n}x_{i}=r\left( \sum_{i=1}^{n}\left\Vert x_{i}\right\Vert
\right) e.  \label{2.3.5}
\end{equation}
\end{corollary}

\begin{proof}
If (\ref{2.3.5}) holds true, then, obviously%
\begin{equation*}
\left\Vert \sum_{i=1}^{n}x_{i}\right\Vert =r\left( \sum_{i=1}^{n}\left\Vert
x_{i}\right\Vert \right) \left\Vert e\right\Vert =r\sum_{i=1}^{n}\left\Vert
x_{i}\right\Vert ,
\end{equation*}%
which is the equality case in (\ref{2.3.2}).

Conversely, if the equality holds in (\ref{2.3.2}), then by Corollary \ref%
{c2.3.1}, we have that (\ref{2.3.3}) and (\ref{2.3.4}) hold true. Utilising
Theorem \ref{t2.3.1}, we conclude that there exists a $\mu >0$ such that%
\begin{equation}
\sum_{i=1}^{n}x_{i}=\mu e.  \label{2.3.6}
\end{equation}%
Inserting this in (\ref{2.3.3}) we get%
\begin{equation*}
\mu \left\Vert e\right\Vert ^{2}=r\sum_{i=1}^{n}\left\Vert x_{i}\right\Vert
\end{equation*}%
giving%
\begin{equation}
\mu =r\sum_{i=1}^{n}\left\Vert x_{i}\right\Vert .  \label{2.3.7}
\end{equation}%
Finally, by (\ref{2.3.6}) and (\ref{2.3.7}) we deduce (\ref{2.3.5}) and the
corollary is proved.
\end{proof}

\section{Other Multiplicative Reverses for $m$ Functionals}

Assume that $F_{k},$ $k\in \left\{ 1,\dots ,m\right\} $ are bounded linear
functionals defined on the normed linear space $X.$

For $p\in \lbrack 1,\infty ),$ define%
\begin{equation}
c_{p}:=\sup_{x\neq 0}\left[ \frac{\sum_{k=1}^{m}\left\vert F_{k}\left(
x\right) \right\vert ^{p}}{\left\Vert x\right\Vert ^{p}}\right] ^{\frac{1}{p}%
}  \tag{$c_{p}$}  \label{cp}
\end{equation}%
and for $p=\infty ,$%
\begin{equation}
c_{\infty }:=\sup_{x\neq 0}\left[ \max_{1\leq k\leq m}\left\{ \frac{%
\left\vert F_{k}\left( x\right) \right\vert }{\left\Vert x\right\Vert }%
\right\} \right] .  \tag{$c_{\infty }$}  \label{cinf}
\end{equation}%
Then, by the fact that $\left\vert F_{k}\left( x\right) \right\vert \leq
\left\Vert F_{k}\right\Vert \left\Vert x\right\Vert $ for any $x\in X,$
where $\left\Vert F_{k}\right\Vert $ is the norm of the functional $F_{k},$
we have that%
\begin{equation*}
c_{p}\leq \left( \sum_{k=1}^{m}\left\Vert F_{k}\right\Vert ^{p}\right) ^{%
\frac{1}{p}},\quad p\geq 1
\end{equation*}%
and%
\begin{equation*}
c_{\infty }\leq \max_{1\leq k\leq m}\left\Vert F_{k}\right\Vert .
\end{equation*}

We may now state and prove a new reverse inequality for the generalised
triangle inequality in normed linear spaces.

\begin{theorem}
\label{ta3.1}Let $x_{i},$ $r_{k},$ $F_{k},$ $k\in \left\{ 1,\dots ,m\right\} 
$, $i\in \left\{ 1,\dots ,n\right\} $ be as in the hypothesis of Theorem \ref%
{t1.2.1}. Then we have the inequalities%
\begin{equation}
\left( 1\leq \right) \frac{\sum_{i=1}^{n}\left\Vert x_{i}\right\Vert }{%
\left\Vert \sum_{i=1}^{n}x_{i}\right\Vert }\leq \frac{c_{\infty }}{%
\max\limits_{1\leq k\leq m}\left\{ r_{k}\right\} }\left( \leq \frac{%
\max\limits_{1\leq k\leq m}\left\Vert F_{k}\right\Vert }{\max\limits_{1\leq
k\leq m}\left\{ r_{k}\right\} }\right) .  \label{a.3.1}
\end{equation}%
The case of equality holds in (\ref{a.3.1}) if and only if%
\begin{equation}
\func{Re}\left[ F_{k}\left( \sum_{i=1}^{n}x_{i}\right) \right]
=r_{k}\sum_{i=1}^{n}\left\Vert x_{i}\right\Vert \quad \text{for each \ }k\in
\left\{ 1,\dots ,m\right\}  \label{a.3.1.a}
\end{equation}%
and%
\begin{equation}
\max\limits_{1\leq k\leq m}\func{Re}\left[ F_{k}\left(
\sum_{i=1}^{n}x_{i}\right) \right] =c_{\infty }\left\Vert
\sum_{i=1}^{n}x_{i}\right\Vert .  \label{a.3.1.b}
\end{equation}
\end{theorem}

\begin{proof}
Since, by the definition of $c_{\infty },$ we have%
\begin{equation*}
c_{\infty }\left\Vert x\right\Vert \geq \max\limits_{1\leq k\leq
m}\left\vert F_{k}\left( x\right) \right\vert ,\quad \text{for any \ }x\in X,
\end{equation*}%
then we can state, for $x=\sum_{i=1}^{n}x_{i},$ that%
\begin{align}
c_{\infty }\left\Vert \sum_{i=1}^{n}x_{i}\right\Vert & \geq
\max\limits_{1\leq k\leq m}\left\vert F_{k}\left( \sum_{i=1}^{n}x_{i}\right)
\right\vert \geq \max\limits_{1\leq k\leq m}\left[ \left\vert \func{Re}%
F_{k}\left( \sum_{i=1}^{n}x_{i}\right) \right\vert \right]  \label{a.3.2} \\
& \geq \max\limits_{1\leq k\leq m}\left[ \func{Re}\sum_{i=1}^{n}F_{k}\left(
x_{i}\right) \right] =\max\limits_{1\leq k\leq m}\left[ \sum_{i=1}^{n}\func{%
Re}F_{k}\left( x_{i}\right) \right] .  \notag
\end{align}%
Utilising the hypothesis (\ref{1.2.1}) we obviously have%
\begin{equation*}
\max\limits_{1\leq k\leq m}\left[ \sum_{i=1}^{n}\func{Re}F_{k}\left(
x_{i}\right) \right] \geq \max\limits_{1\leq k\leq m}\left\{ r_{k}\right\}
\cdot \sum_{i=1}^{n}\left\Vert x_{i}\right\Vert .
\end{equation*}%
Also, $\sum_{i=1}^{n}x_{i}\neq 0,$ because, by the initial assumptions, not
all $r_{k}$ and $x_{i}$ with $k\in \left\{ 1,\dots ,m\right\} $ and $i\in
\left\{ 1,\dots ,n\right\} $ are allowed to be zero. Hence the desired
inequality (\ref{a.3.1}) is obtained.

Now, if (\ref{a.3.1.a}) is valid, then, taking the maximum over $k\in
\left\{ 1,\dots ,m\right\} $ in this equality we get%
\begin{equation*}
\max\limits_{1\leq k\leq m}\func{Re}\left[ F_{k}\left(
\sum_{i=1}^{n}x_{i}\right) \right] =\max\limits_{1\leq k\leq m}\left\{
r_{k}\right\} \left\Vert \sum_{i=1}^{n}x_{i}\right\Vert ,
\end{equation*}%
which, together with (\ref{a.3.1.b}) provides the equality case in (\ref%
{a.3.1}).

Now, if the equality holds in (\ref{a.3.1}), it must hold in all the
inequalities used to prove (\ref{a.3.1}), therefore, we have%
\begin{equation}
\func{Re}F_{k}\left( x_{i}\right) =r_{k}\left\Vert x_{i}\right\Vert \quad 
\text{for each \ }i\in \left\{ 1,\dots ,n\right\} \text{\ \ and \ }k\in
\left\{ 1,\dots ,m\right\}  \label{a.3.3}
\end{equation}%
and, from (\ref{a.3.2}),%
\begin{equation*}
c_{\infty }\left\Vert \sum_{i=1}^{n}x_{i}\right\Vert =\max\limits_{1\leq
k\leq m}\func{Re}\left[ F_{k}\left( \sum_{i=1}^{n}x_{i}\right) \right] ,
\end{equation*}%
which is (\ref{a.3.1.b}).

From (\ref{a.3.3}), on summing over $i\in \left\{ 1,\dots ,n\right\} ,$ we
get (\ref{a.3.1.a}), and the theorem is proved.
\end{proof}

The following result in normed spaces also holds.

\begin{theorem}
\label{at3.2}Let $x_{i},r_{k},F_{k},$ $k\in \left\{ 1,\dots ,m\right\} ,$ $%
i\in \left\{ 1,\dots ,n\right\} $ be as in the hypothesis of Theorem \ref%
{t1.2.1}. Then we have the inequality%
\begin{equation}
\left( 1\leq \right) \frac{\sum_{i=1}^{n}\left\Vert x_{i}\right\Vert }{%
\left\Vert \sum_{i=1}^{n}x_{i}\right\Vert }\leq \frac{c_{p}}{\left(
\sum_{k=1}^{m}r_{k}^{p}\right) ^{\frac{1}{p}}}\left( \leq \frac{%
\sum_{k=1}^{m}\left\Vert F_{k}\right\Vert ^{p}}{\sum_{k=1}^{m}r_{k}^{p}}%
\right) ^{\frac{1}{p}},  \label{a.3.4}
\end{equation}%
where $p\geq 1.$

The case of equality holds in (\ref{a.3.4}) if and only if%
\begin{equation}
\func{Re}\left[ F_{k}\left( \sum_{i=1}^{n}x_{i}\right) \right]
=r_{k}\sum_{i=1}^{n}\left\Vert x_{i}\right\Vert \quad \text{for each \ }k\in
\left\{ 1,\dots ,m\right\}  \label{a.3.5}
\end{equation}%
and 
\begin{equation}
\sum_{k=1}^{m}\left[ \func{Re}F_{k}\left( \sum_{i=1}^{n}x_{i}\right) \right]
^{p}=c_{p}^{p}\left\Vert \sum_{i=1}^{n}x_{i}\right\Vert ^{p}.  \label{a.3.6}
\end{equation}
\end{theorem}

\begin{proof}
By the definition of $c_{p},$ $p\geq 1,$ we have%
\begin{equation*}
c_{p}^{p}\left\Vert x\right\Vert ^{p}\geq \sum_{k=1}^{m}\left\vert
F_{k}\left( x\right) \right\vert ^{p}\quad \text{for any \ }x\in X,
\end{equation*}%
implying that%
\begin{align}
c_{p}^{p}\left\Vert \sum_{i=1}^{n}x_{i}\right\Vert ^{p}& \geq
\sum_{k=1}^{m}\left\vert F_{k}\left( \sum_{i=1}^{n}x_{i}\right) \right\vert
^{p}\geq \sum_{k=1}^{m}\left\vert \func{Re}F_{k}\left(
\sum_{i=1}^{n}x_{i}\right) \right\vert ^{p}  \label{a.3.7} \\
& \geq \sum_{k=1}^{m}\left[ \func{Re}F_{k}\left( \sum_{i=1}^{n}x_{i}\right) %
\right] ^{p}=\sum_{k=1}^{m}\left[ \sum_{i=1}^{n}\func{Re}F_{k}\left(
x_{i}\right) \right] ^{p}.  \notag
\end{align}%
Utilising the hypothesis (\ref{1.2.1}), we obviously have that%
\begin{equation}
\sum_{k=1}^{m}\left[ \sum_{i=1}^{n}\func{Re}F_{k}\left( x_{i}\right) \right]
^{p}\geq \sum_{k=1}^{m}\left[ \sum_{i=1}^{n}r_{k}\left\Vert x_{i}\right\Vert %
\right] ^{p}=\sum_{k=1}^{m}r_{k}^{p}\left( \sum_{i=1}^{n}\left\Vert
x_{i}\right\Vert \right) ^{p}.  \label{a.3.8}
\end{equation}%
Making use of (\ref{a.3.7})\ and (\ref{a.3.8}), we deduce%
\begin{equation*}
c_{p}^{p}\left\Vert \sum_{i=1}^{n}x_{i}\right\Vert ^{p}\geq \left(
\sum_{k=1}^{m}r_{k}^{p}\right) \left( \sum_{i=1}^{n}\left\Vert
x_{i}\right\Vert \right) ^{p},
\end{equation*}%
which implies the desired inequality (\ref{a.3.4}).

If (\ref{a.3.5}) holds true, then, taking the power $p$ and summing over $%
k\in \left\{ 1,\dots ,m\right\} ,$ we deduce%
\begin{equation*}
\sum_{k=1}^{m}\left[ \func{Re}\left[ F_{k}\left( \sum_{i=1}^{n}x_{i}\right) %
\right] \right] ^{p}=\sum_{k=1}^{m}r_{k}^{p}\left( \sum_{i=1}^{n}\left\Vert
x_{i}\right\Vert \right) ^{p},
\end{equation*}%
which, together with (\ref{a.3.6}) shows that the equality case holds true
in (\ref{a.3.4}).

Conversely, if the case of equality holds in (\ref{a.3.4}), then it must
hold in all inequalities needed to prove (\ref{a.3.4}), therefore, we must
have:%
\begin{equation}
\func{Re}F_{k}\left( x_{i}\right) =r_{k}\left\Vert x_{i}\right\Vert \quad 
\text{for each \ }i\in \left\{ 1,\dots ,n\right\} \text{\ \ and \ }k\in
\left\{ 1,\dots ,m\right\}  \label{a.3.9}
\end{equation}%
and, from (\ref{a.3.7}),%
\begin{equation*}
c_{p}^{p}\left\Vert \sum_{i=1}^{n}x_{i}\right\Vert ^{p}=\sum_{k=1}^{m}\left[ 
\func{Re}F_{k}\left( \sum_{i=1}^{n}x_{i}\right) \right] ^{p},
\end{equation*}%
which is exactly (\ref{a.3.6}).

From (\ref{a.3.9}), on summing over $i$ from $1$ to $n,$ we deduce (\ref%
{a.3.5}), and the theorem is proved.
\end{proof}

\section{An Additive Reverse for the Triangle Inequality}

\subsection{The Case of One Functional}

In the following we provide an alternative of the Diaz-Metcalf reverse of
the generalised triangle inequality \cite{DRA2}.

\begin{theorem}[Dragomir, 2004]
\label{t2.5.1}Let $\left( X,\left\Vert \cdot \right\Vert \right) $ be a
normed linear space over the real or complex number field $\mathbb{K}$ and $%
F:X\rightarrow \mathbb{K}$ a linear functional with the property that $%
\left\vert F\left( x\right) \right\vert \leq \left\Vert x\right\Vert $ for
any $x\in X$. If $x_{i}\in X,$ $k_{i}\geq 0,$ $i\in \left\{ 1,\dots
,n\right\} $ are such that%
\begin{equation}
\left( 0\leq \right) \left\Vert x_{i}\right\Vert -\func{Re}F\left(
x_{i}\right) \leq k_{i}\text{ \ for each \ }i\in \left\{ 1,\dots ,n\right\} ,
\label{2.5.1}
\end{equation}%
then we have the inequality%
\begin{equation}
\left( 0\leq \right) \sum_{i=1}^{n}\left\Vert x_{i}\right\Vert -\left\Vert
\sum_{i=1}^{n}x_{i}\right\Vert \leq \sum_{i=1}^{n}k_{i}.  \label{2.5.2}
\end{equation}%
The equality holds in (\ref{2.5.2}) if and only if both%
\begin{equation}
F\left( \sum_{i=1}^{n}x_{i}\right) =\left\Vert
\sum_{i=1}^{n}x_{i}\right\Vert \text{ \ and \ }F\left(
\sum_{i=1}^{n}x_{i}\right) =\sum_{i=1}^{n}\left\Vert x_{i}\right\Vert
-\sum_{i=1}^{n}k_{i}.  \label{2.5.3}
\end{equation}
\end{theorem}

\begin{proof}
If we sum in (\ref{2.5.1}) over $i$ from $1$ to $n,$ then we get%
\begin{equation}
\sum_{i=1}^{n}\left\Vert x_{i}\right\Vert \leq \func{Re}\left[ F\left(
\sum_{i=1}^{n}x_{i}\right) \right] +\sum_{i=1}^{n}k_{i}.  \label{2.5.4}
\end{equation}%
Taking into account that $\left\vert F\left( x\right) \right\vert \leq
\left\Vert x\right\Vert $ for each $x\in X,$ then we may state that%
\begin{align}
\func{Re}\left[ F\left( \sum_{i=1}^{n}x_{i}\right) \right] & \leq \left\vert 
\func{Re}F\left( \sum_{i=1}^{n}x_{i}\right) \right\vert  \label{2.5.5} \\
& \leq \left\vert F\left( \sum_{i=1}^{n}x_{i}\right) \right\vert \leq
\left\Vert \sum_{i=1}^{n}x_{i}\right\Vert .  \notag
\end{align}%
Now, making use of (\ref{2.5.4}) and (\ref{2.5.5}), we deduce (\ref{2.5.2}).

Obviously, if (\ref{2.5.3}) is valid, then the case of equality in (\ref%
{2.5.2}) holds true.

Conversely, if the equality holds in (\ref{2.5.2}), then it must hold in all
the inequalities used to prove (\ref{2.5.2}), therefore we have%
\begin{equation*}
\sum_{i=1}^{n}\left\Vert x_{i}\right\Vert =\func{Re}\left[ F\left(
\sum_{i=1}^{n}x_{i}\right) \right] +\sum_{i=1}^{n}k_{i}
\end{equation*}%
and%
\begin{equation*}
\func{Re}\left[ F\left( \sum_{i=1}^{n}x_{i}\right) \right] =\left\vert
F\left( \sum_{i=1}^{n}x_{i}\right) \right\vert =\left\Vert
\sum_{i=1}^{n}x_{i}\right\Vert ,
\end{equation*}%
which imply (\ref{2.5.3}).
\end{proof}

The following corollary may be stated \cite{DRA2}.

\begin{corollary}
\label{c2.5.1}Let $\left( X,\left\Vert \cdot \right\Vert \right) $ be a
normed linear space, $\left[ \cdot ,\cdot \right] :X\times X\rightarrow 
\mathbb{K}$ a semi-inner product generating the norm $\left\Vert \cdot
\right\Vert $ and $e\in X,$ $\left\Vert e\right\Vert =1.$ If $x_{i}\in X,$ $%
k_{i}\geq 0,$ $\ i\in \left\{ 1,\dots ,n\right\} $ are such that%
\begin{equation}
\left( 0\leq \right) \left\Vert x_{i}\right\Vert -\func{Re}\left[ x_{i},e%
\right] \leq k_{i}\text{ \ for each \ }i\in \left\{ 1,\dots ,n\right\} ,
\label{2.5.6}
\end{equation}%
then we have the inequality%
\begin{equation}
\left( 0\leq \right) \sum_{i=1}^{n}\left\Vert x_{i}\right\Vert -\left\Vert
\sum_{i=1}^{n}x_{i}\right\Vert \leq \sum_{i=1}^{n}k_{i}.  \label{2.5.7}
\end{equation}%
The equality holds in (\ref{2.5.7}) if and only if both%
\begin{equation}
\left[ \sum_{i=1}^{n}x_{i},e\right] =\left\Vert
\sum_{i=1}^{n}x_{i}\right\Vert \text{ \ \ and \ \ }\left[
\sum_{i=1}^{n}x_{i},e\right] =\sum_{i=1}^{n}\left\Vert x_{i}\right\Vert
-\sum_{i=1}^{n}k_{i}.  \label{2.5.8}
\end{equation}%
Moreover, if $\left( X,\left\Vert \cdot \right\Vert \right) $ is strictly
convex, then the case of equality holds in (\ref{2.5.7}) if and only if%
\begin{equation}
\sum_{i=1}^{n}\left\Vert x_{i}\right\Vert \geq \sum_{i=1}^{n}k_{i}
\label{2.5.9}
\end{equation}%
and 
\begin{equation}
\sum_{i=1}^{n}x_{i}=\left( \sum_{i=1}^{n}\left\Vert x_{i}\right\Vert
-\sum_{i=1}^{n}k_{i}\right) \cdot e.  \label{2.5.10}
\end{equation}
\end{corollary}

\begin{proof}
The first part of the corollary is obvious by Theorem \ref{t2.5.1} applied
for the continuous linear functional of unit norm $F_{e},$ $F_{e}\left(
x\right) =\left[ x,e\right] ,$ $x\in X.$ The second part may be shown on
utilising a similar argument to the one from the proof of Corollary \ref%
{c2.3.2}. We omit the details.
\end{proof}

\begin{remark}
\label{r2.5.1}If $X=H,$ $\left( H;\left\langle \cdot ,\cdot \right\rangle
\right) $ is an inner product space, then from Corollary \ref{c2.5.1} we
deduce the additive reverse inequality obtained in Theorem 7 of \cite{SSD3}.
For further similar results in inner product spaces, see \cite{SSD2} and 
\cite{SSD3}.
\end{remark}

\subsection{The Case of $m$ Functionals}

The following result generalising Theorem \ref{t2.5.1} may be stated \cite%
{DRA2}.

\begin{theorem}[Dragomir, 2004]
\label{t2.2.1}Let $\left( X,\left\Vert \cdot \right\Vert \right) $ be a
normed linear space over the real or complex number field $\mathbb{K}$. If $%
F_{k}$, $k\in \left\{ 1,\dots ,m\right\} $ are bounded linear functionals
defined on $X\ $and $x_{i}\in X,$ $M_{ik}\geq 0$ for $i\in \left\{ 1,\dots
,n\right\} $, $k\in \left\{ 1,\dots ,m\right\} $ are such that%
\begin{equation}
\left\Vert x_{i}\right\Vert -\func{Re}F_{k}\left( x_{i}\right) \leq M_{ik}
\label{2.2.1}
\end{equation}%
\ for each \ $i\in \left\{ 1,\dots ,n\right\} ,\ k\in \left\{ 1,\dots
,m\right\} ,$ then we have the inequality%
\begin{equation}
\sum_{i=1}^{n}\left\Vert x_{i}\right\Vert \leq \left\Vert \frac{1}{m}%
\sum_{k=1}^{m}F_{k}\right\Vert \left\Vert \sum_{i=1}^{n}x_{i}\right\Vert +%
\frac{1}{m}\sum_{k=1}^{m}\sum_{i=1}^{n}M_{ik}.  \label{2.2.2}
\end{equation}%
The case of equality holds in (\ref{2.2.2}) if both%
\begin{equation}
\frac{1}{m}\sum_{k=1}^{m}F_{k}\left( \sum_{i=1}^{n}x_{i}\right) =\left\Vert 
\frac{1}{m}\sum_{k=1}^{m}F_{k}\right\Vert \left\Vert
\sum_{i=1}^{n}x_{i}\right\Vert  \label{2.2.3}
\end{equation}%
and%
\begin{equation}
\frac{1}{m}\sum_{k=1}^{m}F_{k}\left( \sum_{i=1}^{n}x_{i}\right)
=\sum_{i=1}^{n}\left\Vert x_{i}\right\Vert -\frac{1}{m}\sum_{k=1}^{m}%
\sum_{j=1}^{n}M_{ik}.  \label{2.2.4}
\end{equation}
\end{theorem}

\begin{proof}
If we sum (\ref{2.2.1}) over $i$ from $1$ to $n,$ then we deduce%
\begin{equation*}
\sum_{i=1}^{n}\left\Vert x_{i}\right\Vert -\func{Re}F_{k}\left(
\sum_{i=1}^{n}x_{i}\right) \leq \sum_{i=1}^{n}M_{ik}
\end{equation*}%
for each $k\in \left\{ 1,\dots ,m\right\} .$

Summing these inequalities over $k$ from $1$ to $m,$ we deduce 
\begin{equation}
\sum_{i=1}^{n}\left\Vert x_{i}\right\Vert \leq \frac{1}{m}\sum_{k=1}^{m}%
\func{Re}F_{k}\left( \sum_{i=1}^{n}x_{i}\right) +\frac{1}{m}%
\sum_{k=1}^{m}\sum_{i=1}^{n}M_{ik}.  \label{2.2.5}
\end{equation}%
Utilising the continuity property of the functionals $F_{k}$ and the
properties of the modulus, we have%
\begin{align}
\sum_{k=1}^{m}\func{Re}F_{k}\left( \sum_{i=1}^{n}x_{i}\right) & \leq
\left\vert \sum_{k=1}^{m}\func{Re}F_{k}\left( \sum_{i=1}^{n}x_{i}\right)
\right\vert  \label{2.2.6} \\
& \leq \left\vert \sum_{k=1}^{m}F_{k}\left( \sum_{i=1}^{n}x_{i}\right)
\right\vert \leq \left\Vert \sum_{k=1}^{m}F_{k}\right\Vert \left\Vert
\sum_{i=1}^{n}x_{i}\right\Vert .  \notag
\end{align}%
Now, by (\ref{2.2.5}) and (\ref{2.2.6}), we deduce (\ref{2.2.2}).

Obviously, if (\ref{2.2.3}) and (\ref{2.2.4}) hold true, then the case of
equality is valid in (\ref{2.2.2}).

Conversely, if the case of equality holds in (\ref{2.2.2}), then it must
hold in all the inequalities used to prove (\ref{2.2.2}). Therefore we have%
\begin{equation*}
\sum_{i=1}^{n}\left\Vert x_{i}\right\Vert =\frac{1}{m}\sum_{k=1}^{m}\func{Re}%
F_{k}\left( \sum_{i=1}^{n}x_{i}\right) +\frac{1}{m}\sum_{k=1}^{m}%
\sum_{i=1}^{n}M_{ik},
\end{equation*}%
\begin{equation*}
\sum_{k=1}^{m}\func{Re}F_{k}\left( \sum_{i=1}^{n}x_{i}\right) =\left\Vert
\sum_{k=1}^{m}F_{k}\right\Vert \left\Vert \sum_{i=1}^{n}x_{i}\right\Vert
\end{equation*}%
and%
\begin{equation*}
\sum_{k=1}^{m}\func{Im}F_{k}\left( \sum_{i=1}^{n}x_{i}\right) =0.
\end{equation*}%
These imply that (\ref{2.2.3}) and (\ref{2.2.4}) hold true, and the theorem
is completely proved.
\end{proof}

\begin{remark}
\label{r2.2.1}If $F_{k},$ $k\in \left\{ 1,\dots ,m\right\} $ are of unit
norm, then, from (\ref{2.2.2}), we deduce the inequality%
\begin{equation}
\sum_{i=1}^{n}\left\Vert x_{i}\right\Vert \leq \left\Vert
\sum_{i=1}^{n}x_{i}\right\Vert +\frac{1}{m}\sum_{k=1}^{m}%
\sum_{i=1}^{n}M_{ik},  \label{2.2.7}
\end{equation}%
which is obviously coarser than (\ref{2.2.2}), but perhaps more useful for
applications.
\end{remark}

\subsection{The Case of Inner Product Spaces}

The case of inner product spaces, in which we may provide a simpler
condition of equality, is of interest in applications \cite{DRA2}.

\begin{theorem}[Dragomir, 2004]
\label{t2.2.2}Let $\left( X,\left\Vert \cdot \right\Vert \right) $ be an
inner product space over the real or complex number field $\mathbb{K}$, $%
e_{k},$ $x_{i}\in H\backslash \left\{ 0\right\} ,$ $k\in \left\{ 1,\dots
,m\right\} ,$ $i\in \left\{ 1,\dots ,n\right\} .$ If $M_{ik}\geq 0$ for $%
i\in \left\{ 1,\dots ,n\right\} ,$ $\left\{ 1,\dots ,n\right\} $ such that%
\begin{equation}
\left\Vert x_{i}\right\Vert -\func{Re}\left\langle x_{i},e_{k}\right\rangle
\leq M_{ik}\text{ }  \label{2.2.8}
\end{equation}%
\ for each \ $i\in \left\{ 1,\dots ,n\right\} ,\ k\in \left\{ 1,\dots
,m\right\} ,$ then we have the inequality%
\begin{equation}
\sum_{i=1}^{n}\left\Vert x_{i}\right\Vert \leq \left\Vert \frac{1}{m}%
\sum_{k=1}^{m}e_{k}\right\Vert \left\Vert \sum_{i=1}^{n}x_{i}\right\Vert +%
\frac{1}{m}\sum_{k=1}^{m}\sum_{i=1}^{n}M_{ik}.  \label{2.2.9}
\end{equation}%
The case of equality holds in (\ref{2.2.9}) if and only if%
\begin{equation}
\sum_{i=1}^{n}\left\Vert x_{i}\right\Vert \geq \frac{1}{m}%
\sum_{k=1}^{m}\sum_{i=1}^{n}M_{ik}  \label{2.2.9.a}
\end{equation}%
and%
\begin{equation}
\sum_{i=1}^{n}x_{i}=\frac{m\left( \sum_{i=1}^{n}\left\Vert x_{i}\right\Vert -%
\frac{1}{m}\sum_{k=1}^{m}\sum_{i=1}^{n}M_{ik}\right) }{\left\Vert
\sum_{k=1}^{m}e_{k}\right\Vert ^{2}}\sum_{k=1}^{m}e_{k}.  \label{2.2.9.1}
\end{equation}
\end{theorem}

\begin{proof}
As in the proof of Theorem \ref{t2.2.1}, we have%
\begin{equation}
\sum_{i=1}^{n}\left\Vert x_{i}\right\Vert \leq \func{Re}\left\langle \frac{1%
}{m}\sum_{k=1}^{m}e_{k},\sum_{i=1}^{n}x_{i}\right\rangle +\frac{1}{m}%
\sum_{k=1}^{m}\sum_{i=1}^{n}M_{ik},  \label{2.2.11}
\end{equation}%
and $\sum_{k=1}^{m}e_{k}\neq 0.$

On utilising the Schwarz inequality in the inner product space $\left(
H;\left\langle \cdot ,\cdot \right\rangle \right) $ for $%
\sum_{i=1}^{n}x_{i}, $ $\sum_{k=1}^{m}e_{k},$ we have%
\begin{align}
\left\Vert \sum_{i=1}^{n}x_{i}\right\Vert \left\Vert
\sum_{k=1}^{m}e_{k}\right\Vert & \geq \left\vert \left\langle
\sum_{i=1}^{n}x_{i},\sum_{k=1}^{m}e_{k}\right\rangle \right\vert
\label{2.2.12} \\
& \geq \left\vert \func{Re}\left\langle
\sum_{i=1}^{n}x_{i},\sum_{k=1}^{m}e_{k}\right\rangle \right\vert  \notag \\
& \geq \func{Re}\left\langle
\sum_{i=1}^{n}x_{i},\sum_{k=1}^{m}e_{k}\right\rangle .  \notag
\end{align}%
By (\ref{2.2.11}) and (\ref{2.2.12}) we deduce (\ref{2.2.9}).

Taking the norm in (\ref{2.2.9.1}) and using (\ref{2.2.9.a}), we have%
\begin{equation*}
\left\Vert \sum_{i=1}^{n}x_{i}\right\Vert =\frac{m\left(
\sum_{i=1}^{n}\left\Vert x_{i}\right\Vert -\frac{1}{m}\sum_{k=1}^{m}%
\sum_{i=1}^{n}M_{ik}\right) }{\left\Vert \sum_{k=1}^{m}e_{k}\right\Vert },
\end{equation*}%
showing that the equality holds in (\ref{2.2.9}).

Conversely, if the case of equality holds in (\ref{2.2.9}), then it must
hold in all the inequalities used to prove (\ref{2.2.9}). Therefore we have%
\begin{equation}
\left\Vert x_{i}\right\Vert =\func{Re}\left\langle x_{i},e_{k}\right\rangle
+M_{ik}  \label{2.2.13}
\end{equation}%
\ for each \ $i\in \left\{ 1,\dots ,n\right\} ,\ k\in \left\{ 1,\dots
,m\right\} ,$%
\begin{equation}
\left\Vert \sum_{i=1}^{n}x_{i}\right\Vert \left\Vert
\sum_{k=1}^{m}e_{k}\right\Vert =\left\vert \left\langle
\sum_{i=1}^{n}x_{i},\sum_{k=1}^{m}e_{k}\right\rangle \right\vert 
\label{2.2.14}
\end{equation}%
and%
\begin{equation}
\func{Im}\left\langle \sum_{i=1}^{n}x_{i},\sum_{k=1}^{m}e_{k}\right\rangle
=0.  \label{2.2.15}
\end{equation}%
From (\ref{2.2.13}), on summing over $i$ and $k,$ we get%
\begin{equation}
\func{Re}\left\langle \sum_{i=1}^{n}x_{i},\sum_{k=1}^{m}e_{k}\right\rangle
=m\sum_{i=1}^{n}\left\Vert x_{i}\right\Vert
-\sum_{k=1}^{m}\sum_{i=1}^{n}M_{ik}.  \label{2.2.16}
\end{equation}%
On the other hand, by the use of the identity (\ref{Id}), the relation (\ref%
{2.2.14}) holds if and only if%
\begin{equation*}
\sum_{i=1}^{n}x_{i}=\frac{\left\langle
\sum_{i=1}^{n}x_{i},\sum_{k=1}^{m}e_{k}\right\rangle }{\left\Vert
\sum_{k=1}^{m}e_{k}\right\Vert ^{2}}\sum_{k=1}^{m}e_{k},
\end{equation*}%
giving, from (\ref{2.2.15}) and (\ref{2.2.16}), that%
\begin{equation*}
\sum_{i=1}^{n}x_{i}=\frac{m\sum_{i=1}^{n}\left\Vert x_{i}\right\Vert
-\sum_{k=1}^{m}\sum_{i=1}^{n}M_{ik}}{\left\Vert
\sum_{k=1}^{m}e_{k}\right\Vert ^{2}}\sum_{k=1}^{m}e_{k}.
\end{equation*}%
If the inequality holds in (\ref{2.2.9}), then obviously (\ref{2.2.9.a}) is
valid, and the theorem is proved.
\end{proof}

\begin{remark}
\label{r2.2.2}If in the above theorem the vectors $\left\{ e_{k}\right\} _{k=%
\overline{1,m}}$ are assumed to be orthogonal, then (\ref{2.2.9}) becomes:%
\begin{equation}
\sum_{i=1}^{n}\left\Vert x_{i}\right\Vert \leq \frac{1}{m}\left(
\sum_{k=1}^{m}\left\Vert e_{k}\right\Vert ^{2}\right) ^{\frac{1}{2}%
}\left\Vert \sum_{i=1}^{n}x_{i}\right\Vert +\frac{1}{m}\sum_{k=1}^{m}%
\sum_{i=1}^{n}M_{ik}.  \label{2.2.16.a}
\end{equation}%
Moreover, if $\left\{ e_{k}\right\} _{k=\overline{1,m}}$ is an orthonormal
family, then (\ref{2.2.16.a}) becomes 
\begin{equation}
\sum_{i=1}^{n}\left\Vert x_{i}\right\Vert \leq \frac{\sqrt{m}}{m}\left\Vert
\sum_{i=1}^{n}x_{i}\right\Vert +\frac{1}{m}\sum_{k=1}^{m}%
\sum_{i=1}^{n}M_{ik},  \label{2.2.16.b}
\end{equation}%
which has been obtained in \cite{SSD3}.
\end{remark}

Before we provide some natural consequences of Theorem \ref{t2.2.2}, we need
some preliminary results concerning another reverse of Schwarz's inequality
in inner product spaces (see for instance \cite[p. 27]{SSD2}).

\begin{lemma}[Dragomir, 2004]
\label{l2.2.1}Let $\left( X,\left\Vert \cdot \right\Vert \right) $ be an
inner product space over the real or complex number field $\mathbb{K}$ and $%
x,a\in H,$ $r>0.$ If $\left\Vert x-a\right\Vert \leq r,$ then we have the
inequality%
\begin{equation}
\left\Vert x\right\Vert \left\Vert a\right\Vert -\func{Re}\left\langle
x,a\right\rangle \leq \frac{1}{2}r^{2}.  \label{2.2.17}
\end{equation}%
The case of equality holds in (\ref{2.2.17}) if and only if%
\begin{equation}
\left\Vert x-a\right\Vert =r\text{ \ and \ }\left\Vert x\right\Vert
=\left\Vert a\right\Vert .  \label{2.2.18}
\end{equation}
\end{lemma}

\begin{proof}
The condition $\left\Vert x-a\right\Vert \leq r$ is clearly equivalent to%
\begin{equation}
\left\Vert x\right\Vert ^{2}+\left\Vert a\right\Vert ^{2}\leq 2\func{Re}%
\left\langle x,a\right\rangle +r^{2}.  \label{2.2.19}
\end{equation}%
Since%
\begin{equation}
2\left\Vert x\right\Vert \left\Vert a\right\Vert \leq \left\Vert
x\right\Vert ^{2}+\left\Vert a\right\Vert ^{2},  \label{2.2.20}
\end{equation}%
with equality if and only if $\left\Vert x\right\Vert =\left\Vert
a\right\Vert ,$ hence by (\ref{2.2.19}) and (\ref{2.2.20}) we deduce (\ref%
{2.2.17}).

The case of equality is obvious.
\end{proof}

Utilising the above lemma we may state the following corollary of Theorem %
\ref{t2.2.2} \cite{DRA2}.

\begin{corollary}
\label{c2.2.1}Let $\left( H;\left\langle \cdot ,\cdot \right\rangle \right)
, $ $e_{k},$ $x_{i}$ be as in Theorem \ref{t2.2.2}. If $r_{ik}>0,$ $i\in
\left\{ 1,\dots ,n\right\} ,$ $k\in \left\{ 1,\dots ,m\right\} $ such that%
\begin{equation}
\left\Vert x_{i}-e_{k}\right\Vert \leq r_{ik}\text{ \ for each \ }i\in
\left\{ 1,\dots ,n\right\} \text{ and }k\in \left\{ 1,\dots ,m\right\} ,
\label{2.2.21}
\end{equation}%
then we have the inequality 
\begin{equation}
\sum_{i=1}^{n}\left\Vert x_{i}\right\Vert \leq \left\Vert \frac{1}{m}%
\sum_{k=1}^{m}e_{k}\right\Vert \left\Vert \sum_{i=1}^{n}x_{i}\right\Vert +%
\frac{1}{2m}\sum_{k=1}^{m}\sum_{i=1}^{n}r_{ik}^{2}.  \label{2.2.22}
\end{equation}%
The equality holds in (\ref{2.2.22}) if and only if%
\begin{equation*}
\sum_{i=1}^{n}\left\Vert x_{i}\right\Vert \geq \frac{1}{2m}%
\sum_{k=1}^{m}\sum_{i=1}^{n}r_{ik}^{2}
\end{equation*}%
and%
\begin{equation*}
\sum_{i=1}^{n}x_{i}=\frac{m\left( \sum_{i=1}^{n}\left\Vert x_{i}\right\Vert -%
\frac{1}{2m}\sum_{k=1}^{m}\sum_{i=1}^{n}r_{ik}^{2}\right) }{\left\Vert
\sum_{k=1}^{m}e_{k}\right\Vert ^{2}}\sum_{k=1}^{m}e_{k}.
\end{equation*}
\end{corollary}

The following lemma may provide another sufficient condition for (\ref{2.2.8}%
) to hold (see also \cite[p. 28]{SSD2}).

\begin{lemma}[Dragomir, 2004]
\label{l2.2.2}Let $\left( H;\left\langle \cdot ,\cdot \right\rangle \right) $
be an inner product space over the real or complex number field $\mathbb{K}$
and $x,y\in H,$ $M\geq m>0.$ If either%
\begin{equation}
\func{Re}\left\langle My-x,x-my\right\rangle \geq 0  \label{2.2.23}
\end{equation}%
or, equivalently,%
\begin{equation}
\left\Vert x-\frac{m+M}{2}y\right\Vert \leq \frac{1}{2}\left( M-m\right)
\left\Vert y\right\Vert ,  \label{2.2.24}
\end{equation}%
holds, then%
\begin{equation}
\left\Vert x\right\Vert \left\Vert y\right\Vert -\func{Re}\left\langle
x,y\right\rangle \leq \frac{1}{4}\cdot \frac{\left( M-m\right) ^{2}}{m+M}%
\left\Vert y\right\Vert ^{2}.  \label{2.2.25}
\end{equation}%
The case of equality holds in (\ref{2.2.25}) if and only if the equality
case is realised in (\ref{2.2.23}) and%
\begin{equation*}
\left\Vert x\right\Vert =\frac{M+m}{2}\left\Vert y\right\Vert .
\end{equation*}
\end{lemma}

The proof is obvious by Lemma \ref{l2.2.1} for $a=\frac{M+m}{2}y$ and $r=%
\frac{1}{2}\left( M-m\right) \left\Vert y\right\Vert .$

Finally, the following corollary of Theorem \ref{t2.2.2} may be stated \cite%
{DRA2}.

\begin{corollary}
\label{c2.2.2}Assume that $\left( H,\left\langle \cdot ,\cdot \right\rangle
\right) ,$ $e_{k},$ $x_{i}$ are as in Theorem \ref{t2.2.2}. If $M_{ik}\geq
m_{ik}>0$ satisfy the condition%
\begin{equation*}
\func{Re}\left\langle M_{k}e_{k}-x_{i},x_{i}-\mu _{k}e_{k}\right\rangle \geq
0\text{ }
\end{equation*}%
\ for each \ $i\in \left\{ 1,\dots ,n\right\} $ and $k\in \left\{ 1,\dots
,m\right\} ,$ then%
\begin{equation*}
\sum_{i=1}^{n}\left\Vert x_{i}\right\Vert \leq \left\Vert \frac{1}{m}%
\sum_{k=1}^{m}e_{k}\right\Vert \left\Vert \sum_{i=1}^{n}x_{i}\right\Vert +%
\frac{1}{4m}\sum_{k=1}^{m}\sum_{i=1}^{n}\frac{\left( M_{ik}-m_{ik}\right)
^{2}}{M_{ik}+m_{ik}}\left\Vert e_{k}\right\Vert ^{2}.
\end{equation*}
\end{corollary}

\section{Other Additive Reverses for $m$ Functionals}

A different approach in obtaining other additive reverses for the
generalised triangle inequality is incorporated in the following new result:

\begin{theorem}
\label{bt3.2}Let $\left( X,\left\Vert \cdot \right\Vert \right) $ be a
normed linear space over the real or complex number field $\mathbb{K}$.
Assume $F_{k\text{ }},$ $k\in \left\{ 1,\dots ,m\right\} ,$ are bounded
linear functionals on the normed linear space $X$ and $x_{i}\in X,$ $i\in
\left\{ 1,\dots ,n\right\} ,$ $M_{ik}\geq 0,$ $i\in \left\{ 1,\dots
,n\right\} ,$ $k\in \left\{ 1,\dots ,m\right\} $ are such that%
\begin{equation}
\left\Vert x_{i}\right\Vert -\func{Re}F_{k}\left( x_{i}\right) \leq
M_{ik}\quad  \label{b.3.5}
\end{equation}%
for each \ $i\in \left\{ 1,\dots ,n\right\} $\ \ and \ $k\in \left\{ 1,\dots
,m\right\} .$

\begin{enumerate}
\item[(i)] If $c_{\infty }$ is defined by (\ref{cinf}), then we have the
inequality%
\begin{equation}
\sum_{i=1}^{n}\left\Vert x_{i}\right\Vert \leq c_{\infty }\left\Vert
\sum_{i=1}^{n}x_{i}\right\Vert +\frac{1}{m}\sum_{k=1}^{m}%
\sum_{i=1}^{n}M_{ik}.  \label{b.3.6}
\end{equation}

\item[(ii)] If $c_{p}$ is defined by (\ref{cp}) for $p\geq 1,$ then we have
the inequality:%
\begin{equation}
\sum_{i=1}^{n}\left\Vert x_{i}\right\Vert \leq \frac{1}{m^{\frac{1}{p}}}%
c_{p}\left\Vert \sum_{i=1}^{n}x_{i}\right\Vert +\frac{1}{m}%
\sum_{k=1}^{m}\sum_{i=1}^{n}M_{ik}.  \label{b.3.8}
\end{equation}
\end{enumerate}
\end{theorem}

\begin{proof}
\noindent (i) Since%
\begin{equation*}
\max_{1\leq k\leq m}\left\Vert F_{k}\left( x\right) \right\Vert \leq
c_{\infty }\left\Vert x\right\Vert \quad \text{for any }x\in X,
\end{equation*}%
then we have%
\begin{equation*}
\sum_{k=1}^{m}\left\vert F_{k}\left( \sum_{i=1}^{n}x_{i}\right) \right\vert
\leq m\max_{1\leq k\leq m}\left\vert F_{k}\left( \sum_{i=1}^{n}x_{i}\right)
\right\vert \leq mc_{\infty }\left\Vert \sum_{i=1}^{n}x_{i}\right\Vert .
\end{equation*}%
Using (\ref{2.2.6}), we may state that%
\begin{equation*}
\frac{1}{m}\sum_{k=1}^{m}\func{Re}F_{k}\left( \sum_{i=1}^{n}x_{i}\right)
\leq c_{\infty }\left\Vert \sum_{i=1}^{n}x_{i}\right\Vert ,
\end{equation*}%
which, together with (\ref{2.2.5}) imply the desired inequality (\ref{b.3.6}%
).

\noindent (ii) Using the fact that, obviously%
\begin{equation*}
\left( \sum_{k=1}^{m}\left\vert F_{k}\left( x\right) \right\vert ^{p}\right)
^{\frac{1}{p}}\leq c_{p}\left\Vert x\right\Vert \quad \text{for any }x\in X,
\end{equation*}%
then, by H\"{o}lder's inequality for $p>1,\frac{1}{p}+\frac{1}{q}=1,$ we have%
\begin{eqnarray*}
\sum_{k=1}^{m}\left\vert F_{k}\left( \sum_{i=1}^{n}x_{i}\right) \right\vert
&\leq &m^{\frac{1}{q}}\left( \sum_{k=1}^{m}\left\vert F_{k}\left(
\sum_{i=1}^{n}x_{i}\right) \right\vert ^{p}\right) ^{\frac{1}{p}} \\
&\leq &c_{p}m^{\frac{1}{q}}\left\Vert \sum_{i=1}^{n}x_{i}\right\Vert ,
\end{eqnarray*}%
which, combined with (\ref{2.2.5}) and (\ref{2.2.6}) will give the desired
inequality (\ref{b.3.8}).

The case $p=1$ goes likewise and we omit the details.
\end{proof}

\begin{remark}
Since, obviously $c_{\infty }\leq \max\limits_{1\leq k\leq m}\left\Vert
F_{k}\right\Vert ,$ then from (\ref{b.3.6}) we have%
\begin{equation}
\sum_{i=1}^{n}\left\Vert x_{i}\right\Vert \leq \max\limits_{1\leq k\leq
m}\left\{ \left\Vert F_{k}\right\Vert \right\} \cdot \left\Vert
\sum_{i=1}^{n}x_{i}\right\Vert +\frac{1}{m}\sum_{k=1}^{m}%
\sum_{i=1}^{n}M_{ik}.  \label{b.3.9}
\end{equation}%
Finally, since $c_{p}\leq \left( \sum_{k=1}^{m}\left\Vert F_{k}\right\Vert
^{p}\right) ^{\frac{1}{p}},p\geq 1,$ hence by (\ref{b.3.8}) we have%
\begin{equation}
\sum_{i=1}^{n}\left\Vert x_{i}\right\Vert \leq \left( \frac{%
\sum_{k=1}^{m}\left\Vert F_{k}\right\Vert ^{p}}{m}\right) ^{\frac{1}{p}%
}\left\Vert \sum_{i=1}^{n}x_{i}\right\Vert +\frac{1}{m}\sum_{k=1}^{m}%
\sum_{i=1}^{n}M_{ik}.  \label{b.3.11}
\end{equation}
\end{remark}

The following corollary for semi-inner products may be stated as well.

\begin{corollary}
\label{b.c.3.3}Let $\left( X,\left\Vert \cdot \right\Vert \right) $ be a
real or complex normed space and $\left[ \cdot ,\cdot \right] :X\times
X\rightarrow \mathbb{K}$ a semi-inner product generating the norm $%
\left\Vert \cdot \right\Vert .$ Assume $e_{k},$ $x_{i}\in H$ and $M_{ik}\geq
0,$ $i\in \left\{ 1,\dots ,n\right\} ,$ $k\in \left\{ 1,\dots ,m\right\} $
are such that%
\begin{equation}
\left\Vert x_{i}\right\Vert -\func{Re}\left[ x_{i},e_{k}\right] \leq M_{ik},
\label{b.3.12}
\end{equation}%
for any $i\in \left\{ 1,\dots ,n\right\} ,$ $k\in \left\{ 1,\dots ,m\right\}
.$

\begin{enumerate}
\item[(i)] If%
\begin{equation*}
d_{\infty }:=\sup_{x\neq 0}\left\{ \frac{\max_{1\leq k\leq n}\left\vert %
\left[ x,e_{k}\right] \right\vert }{\left\Vert x\right\Vert }\right\} \left(
\leq \max_{1\leq k\leq n}\left\Vert e_{k}\right\Vert \right) ,
\end{equation*}%
then%
\begin{align}
\sum_{i=1}^{n}\left\Vert x_{i}\right\Vert & \leq d_{\infty }\left\Vert
\sum_{i=1}^{n}x_{i}\right\Vert +\frac{1}{m}\sum_{k=1}^{m}\sum_{i=1}^{n}M_{ik}
\label{b.3.14} \\
& \left( \leq \max_{1\leq k\leq n}\left\Vert e_{k}\right\Vert \cdot
\left\Vert \sum_{i=1}^{n}x_{i}\right\Vert +\frac{1}{m}\sum_{k=1}^{m}%
\sum_{i=1}^{n}M_{ik}\right) ;  \notag
\end{align}

\item[(ii)] If%
\begin{equation*}
d_{p}:=\sup_{x\neq 0}\left\{ \frac{\sum_{k=1}^{m}\left\vert \left[ x,e_{k}%
\right] \right\vert ^{p}}{\left\Vert x\right\Vert ^{p}}\right\} ^{\frac{1}{p}%
}\left( \leq \left( \sum_{k=1}^{m}\left\Vert e_{k}\right\Vert ^{p}\right) ^{%
\frac{1}{p}}\right) ,
\end{equation*}%
where $p\geq 1,$ then%
\begin{align}
\sum_{i=1}^{n}\left\Vert x_{i}\right\Vert & \leq \frac{1}{m^{\frac{1}{p}}}%
d_{p}\left\Vert \sum_{i=1}^{n}x_{i}\right\Vert +\frac{1}{m}%
\sum_{k=1}^{m}\sum_{i=1}^{n}M_{ik}  \label{b.3.15} \\
& \left( \leq \left( \frac{\sum_{k=1}^{m}\left\Vert e_{k}\right\Vert ^{p}}{m}%
\right) ^{\frac{1}{p}}\left\Vert \sum_{i=1}^{n}x_{i}\right\Vert +\frac{1}{m}%
\sum_{k=1}^{m}\sum_{i=1}^{n}M_{ik}\right) .  \notag
\end{align}
\end{enumerate}
\end{corollary}

\section{Applications for Complex Numbers}

Let $\mathbb{C}$ be the field of complex numbers. If $z=\func{Re}z+i\func{Im}%
z,$ then by $\left\vert \cdot \right\vert _{p}:\mathbb{C}\rightarrow \lbrack
0,\infty ),$ $p\in \left[ 1,\infty \right] $ we define the $p-$\textit{%
modulus }of $z$ as%
\begin{equation*}
\left\vert z\right\vert _{p}:=\left\{ 
\begin{array}{ll}
\max \left\{ \left\vert \func{Re}z\right\vert ,\left\vert \func{Im}%
z\right\vert \right\} & \text{if \ }p=\infty , \\ 
&  \\ 
\left( \left\vert \func{Re}z\right\vert ^{p}+\left\vert \func{Im}%
z\right\vert ^{p}\right) ^{\frac{1}{p}} & \text{if \ }p\in \lbrack 1,\infty
),%
\end{array}%
\right.
\end{equation*}%
where $\left\vert a\right\vert ,$ $a\in \mathbb{R}$ is the usual modulus of
the real number $a.$

For $p=2,$ we recapture the usual modulus of a complex number, i.e.,%
\begin{equation*}
\left\vert z\right\vert _{2}=\sqrt{\left\vert \func{Re}z\right\vert
^{2}+\left\vert \func{Im}z\right\vert ^{2}}=\left\vert z\right\vert ,\quad
z\in \mathbb{C}.
\end{equation*}%
It is well known that $\left( \mathbb{C},\left\vert \cdot \right\vert
_{p}\right) ,$ $p\in \left[ 1,\infty \right] $ is a Banach space over the%
\textit{\ real number field} $\mathbb{R}$.

Consider the Banach space $\left( \mathbb{C},\left\vert \cdot \right\vert
_{1}\right) $ and $F:\mathbb{C\rightarrow C}$, $F\left( z\right) =az$ with $%
a\in \mathbb{C}$, $a\neq 0.$ Obviously, $F$ is linear on $\mathbb{C}$. For $%
z\neq 0,$ we have%
\begin{equation*}
\frac{\left\vert F\left( z\right) \right\vert }{\left\vert z\right\vert _{1}}%
=\frac{\left\vert a\right\vert \left\vert z\right\vert }{\left\vert
z\right\vert _{1}}=\frac{\left\vert a\right\vert \sqrt{\left\vert \func{Re}%
z\right\vert ^{2}+\left\vert \func{Im}z\right\vert ^{2}}}{\left\vert \func{Re%
}z\right\vert +\left\vert \func{Im}z\right\vert }\leq \left\vert
a\right\vert .
\end{equation*}%
Since, for $z_{0}=1,$ we have $\left\vert F\left( z_{0}\right) \right\vert
=\left\vert a\right\vert $ and $\left\vert z_{0}\right\vert _{1}=1,$ hence%
\begin{equation*}
\left\Vert F\right\Vert _{1}:=\sup_{z\neq 0}\frac{\left\vert F\left(
z\right) \right\vert }{\left\vert z\right\vert _{1}}=\left\vert a\right\vert
,
\end{equation*}%
showing that $F$ is a bounded linear functional on $\left( \mathbb{C}%
,\left\vert \cdot \right\vert _{1}\right) $ and $\left\Vert F\right\Vert
_{1}=\left\vert a\right\vert .$

We can apply Theorem \ref{t1.2.1} to state the following reverse of the
generalised triangle inequality\ for complex numbers \cite{DRA1}.

\begin{proposition}
\label{p1.4.1}Let $a_{k},$ $x_{j}\in \mathbb{C}$, $k\in \left\{ 1,\dots
,m\right\} $ and $j\in \left\{ 1,\dots ,n\right\} .$ If there exist the
constants $r_{k}\geq 0,$ $k\in \left\{ 1,\dots ,m\right\} $ with $%
\sum_{k=1}^{m}r_{k}>0$ and%
\begin{equation}
r_{k}\left[ \left\vert \func{Re}x_{j}\right\vert +\left\vert \func{Im}%
x_{j}\right\vert \right] \leq \func{Re}a_{k}\cdot \func{Re}x_{j}-\func{Im}%
a_{k}\cdot \func{Im}x_{j}  \label{1.4.1}
\end{equation}%
for each $j\in \left\{ 1,\dots ,n\right\} $ and $k\in \left\{ 1,\dots
,m\right\} ,$ then%
\begin{equation}
\sum_{j=1}^{n}\left[ \left\vert \func{Re}x_{j}\right\vert +\left\vert \func{%
Im}x_{j}\right\vert \right] \leq \frac{\left\vert
\sum_{k=1}^{m}a_{k}\right\vert }{\sum_{k=1}^{m}r_{k}}\left[ \left\vert
\sum_{j=1}^{n}\func{Re}x_{j}\right\vert +\left\vert \sum_{j=1}^{n}\func{Im}%
x_{j}\right\vert \right] .  \label{1.4.2}
\end{equation}%
The case of equality holds in (\ref{1.4.2}) if both%
\begin{align*}
& \func{Re}\left( \sum_{k=1}^{m}a_{k}\right) \func{Re}\left(
\sum_{j=1}^{n}x_{j}\right) -\func{Im}\left( \sum_{k=1}^{m}a_{k}\right) \func{%
Im}\left( \sum_{j=1}^{n}x_{j}\right) \\
& =\left( \sum_{k=1}^{m}r_{k}\right) \sum_{j=1}^{n}\left[ \left\vert \func{Re%
}x_{j}\right\vert +\left\vert \func{Im}x_{j}\right\vert \right] \\
& =\left\vert \sum_{k=1}^{m}a_{k}\right\vert \left[ \left\vert \sum_{j=1}^{n}%
\func{Re}x_{j}\right\vert +\left\vert \sum_{j=1}^{n}\func{Im}%
x_{j}\right\vert \right] .
\end{align*}
\end{proposition}

The proof follows by Theorem \ref{t1.2.1} applied for the Banach space $%
\left( \mathbb{C},\left\vert \cdot \right\vert _{1}\right) $ and $%
F_{k}\left( z\right) =a_{k}z,$ $k\in \left\{ 1,\dots ,m\right\} $ on taking
into account that:%
\begin{equation*}
\left\Vert \sum_{k=1}^{m}F_{k}\right\Vert _{1}=\left\vert
\sum_{k=1}^{m}a_{k}\right\vert .
\end{equation*}%
Now, consider the Banach space $\left( \mathbb{C},\left\vert \cdot
\right\vert _{\infty }\right) .$ If $F\left( z\right) =dz,$ then for $z\neq
0 $ we have%
\begin{equation*}
\frac{\left\vert F\left( z\right) \right\vert }{\left\vert z\right\vert
_{\infty }}=\frac{\left\vert d\right\vert \left\vert z\right\vert }{%
\left\vert z\right\vert _{\infty }}=\frac{\left\vert d\right\vert \sqrt{%
\left\vert \func{Re}z\right\vert ^{2}+\left\vert \func{Im}z\right\vert ^{2}}%
}{\max \left\{ \left\vert \func{Re}z\right\vert ,\left\vert \func{Im}%
z\right\vert \right\} }\leq \sqrt{2}\left\vert d\right\vert .
\end{equation*}%
Since, for $z_{0}=1+i,$ we have $\left\vert F\left( z_{0}\right) \right\vert
=\sqrt{2}\left\vert d\right\vert ,$ $\left\vert z_{0}\right\vert _{\infty
}=1,$ hence%
\begin{equation*}
\left\Vert F\right\Vert _{\infty }:=\sup_{z\neq 0}\frac{\left\vert F\left(
z\right) \right\vert }{\left\vert z\right\vert _{\infty }}=\sqrt{2}%
\left\vert d\right\vert ,
\end{equation*}%
showing that $F$ is a bounded linear functional on $\left( \mathbb{C}%
,\left\vert \cdot \right\vert _{\infty }\right) $ and $\left\Vert
F\right\Vert _{\infty }=\sqrt{2}\left\vert d\right\vert .$

If we apply Theorem \ref{t1.2.1}, then we can state the following reverse of
the generalised triangle inequality for complex numbers \cite{DRA1}.

\begin{proposition}
\label{p1.4.2}Let $a_{k},$ $x_{j}\in \mathbb{C}$, $k\in \left\{ 1,\dots
,m\right\} $ and $j\in \left\{ 1,\dots ,n\right\} .$ If there exist the
constants $r_{k}\geq 0,$ $k\in \left\{ 1,\dots ,m\right\} $ with $%
\sum_{k=1}^{m}r_{k}>0$ and%
\begin{equation*}
r_{k}\max \left\{ \left\vert \func{Re}x_{j}\right\vert ,\left\vert \func{Im}%
x_{j}\right\vert \right\} \leq \func{Re}a_{k}\cdot \func{Re}x_{j}-\func{Im}%
a_{k}\cdot \func{Im}x_{j}
\end{equation*}%
for each $j\in \left\{ 1,\dots ,n\right\} $ and $k\in \left\{ 1,\dots
,m\right\} ,$ then%
\begin{multline}
\sum_{j=1}^{n}\max \left\{ \left\vert \func{Re}x_{j}\right\vert ,\left\vert 
\func{Im}x_{j}\right\vert \right\}  \label{1.4.3} \\
\leq \sqrt{2}\cdot \frac{\left\vert \sum_{k=1}^{m}a_{k}\right\vert }{%
\sum_{k=1}^{m}r_{k}}\max \left\{ \left\vert \sum_{j=1}^{n}\func{Re}%
x_{j}\right\vert ,\left\vert \sum_{j=1}^{n}\func{Im}x_{j}\right\vert
\right\} .
\end{multline}%
The case of equality holds in (\ref{1.4.3}) if both%
\begin{align*}
& \func{Re}\left( \sum_{k=1}^{m}a_{k}\right) \func{Re}\left(
\sum_{j=1}^{n}x_{j}\right) -\func{Im}\left( \sum_{k=1}^{m}a_{k}\right) \func{%
Im}\left( \sum_{j=1}^{n}x_{j}\right) \\
& =\left( \sum_{k=1}^{m}r_{k}\right) \sum_{j=1}^{n}\max \left\{ \left\vert 
\func{Re}x_{j}\right\vert ,\left\vert \func{Im}x_{j}\right\vert \right\} \\
& =\sqrt{2}\left\vert \sum_{k=1}^{m}a_{k}\right\vert \max \left\{ \left\vert
\sum_{j=1}^{n}\func{Re}x_{j}\right\vert ,\left\vert \sum_{j=1}^{n}\func{Im}%
x_{j}\right\vert \right\} .
\end{align*}
\end{proposition}

Finally, consider the Banach space $\left( \mathbb{C},\left\vert \cdot
\right\vert _{2p}\right) $ with $p\geq 1.$

Let $F:\mathbb{C\rightarrow C}$, $F\left( z\right) =cz.$ By H\"{o}lder's
inequality, we have%
\begin{equation*}
\frac{\left\vert F\left( z\right) \right\vert }{\left\vert z\right\vert _{2p}%
}=\frac{\left\vert c\right\vert \sqrt{\left\vert \func{Re}z\right\vert
^{2}+\left\vert \func{Im}z\right\vert ^{2}}}{\left( \left\vert \func{Re}%
z\right\vert ^{2p}+\left\vert \func{Im}z\right\vert ^{2p}\right) ^{\frac{1}{%
2p}}}\leq 2^{\frac{1}{2}-\frac{1}{2p}}\left\vert c\right\vert .
\end{equation*}%
Since, for $z_{0}=1+i$ we have $\left\vert F\left( z_{0}\right) \right\vert
=2^{\frac{1}{2}}\left\vert c\right\vert ,$ $\left\vert z_{0}\right\vert
_{2p}=2^{\frac{1}{2p}}$ $\left( p\geq 1\right) ,$ hence%
\begin{equation*}
\left\Vert F\right\Vert _{2p}:=\sup_{z\neq 0}\frac{\left\vert F\left(
z\right) \right\vert }{\left\vert z\right\vert _{2p}}=2^{\frac{1}{2}-\frac{1%
}{2p}}\left\vert c\right\vert ,
\end{equation*}%
showing that $F$ is a bounded linear functional on $\left( \mathbb{C}%
,\left\vert \cdot \right\vert _{2p}\right) ,$ $p\geq 1$ and $\left\Vert
F\right\Vert _{2p}=2^{\frac{1}{2}-\frac{1}{2p}}\left\vert c\right\vert .$

If we apply Theorem \ref{t1.2.1}, then we can state the following
proposition \cite{DRA1}.

\begin{proposition}
\label{p1.4.3}Let $a_{k},$ $x_{j}\in \mathbb{C}$, $k\in \left\{ 1,\dots
,m\right\} $ and $j\in \left\{ 1,\dots ,n\right\} .$ If there exist the
constants $r_{k}\geq 0,$ $k\in \left\{ 1,\dots ,m\right\} $ with $%
\sum_{k=1}^{m}r_{k}>0$ and%
\begin{equation*}
r_{k}\left[ \left\vert \func{Re}x_{j}\right\vert ^{2p}+\left\vert \func{Im}%
x_{j}\right\vert ^{2p}\right] ^{\frac{1}{2p}}\leq \func{Re}a_{k}\cdot \func{%
Re}x_{j}-\func{Im}a_{k}\cdot \func{Im}x_{j}
\end{equation*}%
for each $j\in \left\{ 1,\dots ,n\right\} $ and $k\in \left\{ 1,\dots
,m\right\} ,$ then%
\begin{multline}
\sum_{j=1}^{n}\left[ \left\vert \func{Re}x_{j}\right\vert ^{2p}+\left\vert 
\func{Im}x_{j}\right\vert ^{2p}\right] ^{\frac{1}{2p}}  \label{1.4.4} \\
\leq 2^{\frac{1}{2}-\frac{1}{2p}}\frac{\left\vert
\sum_{k=1}^{m}a_{k}\right\vert }{\sum_{k=1}^{m}r_{k}}\left[ \left\vert
\sum_{j=1}^{n}\func{Re}x_{j}\right\vert ^{2p}+\left\vert \sum_{j=1}^{n}\func{%
Im}x_{j}\right\vert ^{2p}\right] ^{\frac{1}{2p}}.
\end{multline}%
The case of equality holds in (\ref{1.4.4}) if both:%
\begin{align*}
& \func{Re}\left( \sum_{k=1}^{m}a_{k}\right) \func{Re}\left(
\sum_{j=1}^{n}x_{j}\right) -\func{Im}\left( \sum_{k=1}^{m}a_{k}\right) \func{%
Im}\left( \sum_{j=1}^{n}x_{j}\right) \\
& =\left( \sum_{k=1}^{m}r_{k}\right) \sum_{j=1}^{n}\left[ \left\vert \func{Re%
}x_{j}\right\vert ^{2p}+\left\vert \func{Im}x_{j}\right\vert ^{2p}\right] ^{%
\frac{1}{2p}} \\
& =2^{\frac{1}{2}-\frac{1}{2p}}\left\vert \sum_{k=1}^{m}a_{k}\right\vert %
\left[ \left\vert \sum_{j=1}^{n}\func{Re}x_{j}\right\vert ^{2p}+\left\vert
\sum_{j=1}^{n}\func{Im}x_{j}\right\vert ^{2p}\right] ^{\frac{1}{2p}}.
\end{align*}
\end{proposition}

\begin{remark}
If in the above proposition we choose $p=1,$ then we have the following
reverse of the generalised triangle inequality for complex numbers%
\begin{equation*}
\sum_{j=1}^{n}\left\vert x_{j}\right\vert \leq \frac{\left\vert
\sum_{k=1}^{m}a_{k}\right\vert }{\sum_{k=1}^{m}r_{k}}\left\vert
\sum_{j=1}^{n}x_{j}\right\vert
\end{equation*}%
provided $x_{j},a_{k},$ $j\in \left\{ 1,\dots ,n\right\} $, $k\in \left\{
1,\dots ,m\right\} $ satisfy the assumption%
\begin{equation*}
r_{k}\left\vert x_{j}\right\vert \leq \func{Re}a_{k}\cdot \func{Re}x_{j}-%
\func{Im}a_{k}\cdot \func{Im}x_{j}
\end{equation*}%
for each $j\in \left\{ 1,\dots ,n\right\} $, $k\in \left\{ 1,\dots
,m\right\} .$ Here $\left\vert \cdot \right\vert $ is the usual modulus of a
complex number and $r_{k}>0,$ $k\in \left\{ 1,\dots ,m\right\} $ are given.
\end{remark}

We can apply Theorem \ref{t2.2.1} to state the following reverse of the
generalised triangle inequality\ for complex numbers \cite{DRA2}.

\begin{proposition}
\label{p2.4.1}Let $a_{k},$ $x_{j}\in \mathbb{C}$, $k\in \left\{ 1,\dots
,m\right\} $ and $j\in \left\{ 1,\dots ,n\right\} .$ If there exist the
constants $M_{jk}\geq 0,$ $k\in \left\{ 1,\dots ,m\right\} ,$ $j\in \left\{
1,\dots ,n\right\} $ such that%
\begin{equation}
\left\vert \func{Re}x_{j}\right\vert +\left\vert \func{Im}x_{j}\right\vert
\leq \func{Re}a_{k}\cdot \func{Re}x_{j}-\func{Im}a_{k}\cdot \func{Im}%
x_{j}+M_{jk}  \label{2.4.1}
\end{equation}%
for each $j\in \left\{ 1,\dots ,n\right\} $ and $k\in \left\{ 1,\dots
,m\right\} ,$ then%
\begin{multline}
\quad \sum_{j=1}^{n}\left[ \left\vert \func{Re}x_{j}\right\vert +\left\vert 
\func{Im}x_{j}\right\vert \right]  \label{2.4.2} \\
\leq \frac{1}{m}\left\vert \sum_{k=1}^{m}a_{k}\right\vert \left[ \left\vert
\sum_{j=1}^{n}\func{Re}x_{j}\right\vert +\left\vert \sum_{j=1}^{n}\func{Im}%
x_{j}\right\vert \right] +\frac{1}{m}\sum_{k=1}^{m}\sum_{j=1}^{n}M_{jk}.\quad
\end{multline}
\end{proposition}

The proof follows by Theorem \ref{t2.2.1} applied for the Banach space $%
\left( \mathbb{C},\left\vert \cdot \right\vert _{1}\right) $ and $%
F_{k}\left( z\right) =a_{k}z,$ $k\in \left\{ 1,\dots ,m\right\} $ on taking
into account that:%
\begin{equation*}
\left\Vert \sum_{k=1}^{m}F_{k}\right\Vert _{1}=\left\vert
\sum_{k=1}^{m}a_{k}\right\vert .
\end{equation*}

If we apply Theorem \ref{t2.2.1} for the Banach space $\left( \mathbb{C}%
,\left\vert \cdot \right\vert _{\infty }\right) $, then we can state the
following reverse of the generalised triangle inequality for complex numbers 
\cite{DRA2}.

\begin{proposition}
\label{p2.4.2}Let $a_{k},$ $x_{j}\in \mathbb{C}$, $k\in \left\{ 1,\dots
,m\right\} $ and $j\in \left\{ 1,\dots ,n\right\} .$ If there exist the
constants $M_{jk}\geq 0,$ $k\in \left\{ 1,\dots ,m\right\} ,$ $j\in \left\{
1,\dots ,n\right\} $ such that%
\begin{equation*}
\max \left\{ \left\vert \func{Re}x_{j}\right\vert ,\left\vert \func{Im}%
x_{j}\right\vert \right\} \leq \func{Re}a_{k}\cdot \func{Re}x_{j}-\func{Im}%
a_{k}\cdot \func{Im}x_{j}+M_{jk}
\end{equation*}%
for each $j\in \left\{ 1,\dots ,n\right\} $ and $k\in \left\{ 1,\dots
,m\right\} ,$ then%
\begin{multline}
\sum_{j=1}^{n}\max \left\{ \left\vert \func{Re}x_{j}\right\vert ,\left\vert 
\func{Im}x_{j}\right\vert \right\}  \label{2.4.3} \\
\leq \frac{\sqrt{2}}{m}\left\vert \sum_{k=1}^{m}a_{k}\right\vert \max
\left\{ \left\vert \sum_{j=1}^{n}\func{Re}x_{j}\right\vert ,\left\vert
\sum_{j=1}^{n}\func{Im}x_{j}\right\vert \right\} +\frac{1}{m}%
\sum_{k=1}^{m}\sum_{j=1}^{n}M_{jk}.
\end{multline}
\end{proposition}

Finally, if we apply Theorem \ref{t2.2.1}, for the Banach space $\left( 
\mathbb{C},\left\vert \cdot \right\vert _{2p}\right) $ with $p\geq 1,$ then
we can state the following proposition \cite{DRA2}.

\begin{proposition}
\label{p2.4.3}Let $a_{k},$ $x_{j},$ $M_{jk}$ be as in Proposition \ref%
{p2.4.2}. If%
\begin{equation*}
\left[ \left\vert \func{Re}x_{j}\right\vert ^{2p}+\left\vert \func{Im}%
x_{j}\right\vert ^{2p}\right] ^{\frac{1}{2p}}\leq \func{Re}a_{k}\cdot \func{%
Re}x_{j}-\func{Im}a_{k}\cdot \func{Im}x_{j}+M_{jk}
\end{equation*}%
for each $j\in \left\{ 1,\dots ,n\right\} $ and $k\in \left\{ 1,\dots
,m\right\} ,$ then%
\begin{multline}
\sum_{j=1}^{n}\left[ \left\vert \func{Re}x_{j}\right\vert ^{2p}+\left\vert 
\func{Im}x_{j}\right\vert ^{2p}\right] ^{\frac{1}{2p}}  \label{2.4.4} \\
\leq \frac{2^{\frac{1}{2}-\frac{1}{2p}}}{m}\left\vert
\sum_{k=1}^{m}a_{k}\right\vert \left[ \left\vert \sum_{j=1}^{n}\func{Re}%
x_{j}\right\vert ^{2p}+\left\vert \sum_{j=1}^{n}\func{Im}x_{j}\right\vert
^{2p}\right] ^{\frac{1}{2p}}+\frac{1}{m}\sum_{k=1}^{m}\sum_{j=1}^{n}M_{jk}.
\end{multline}%
where $p\geq 1.$
\end{proposition}

\begin{remark}
If in the above proposition we choose $p=1,$ then we have the following
reverse of the generalised triangle inequality for complex numbers%
\begin{equation*}
\sum_{j=1}^{n}\left\vert x_{j}\right\vert \leq \left\vert \frac{1}{m}%
\sum_{k=1}^{m}a_{k}\right\vert \left\vert \sum_{j=1}^{n}x_{j}\right\vert +%
\frac{1}{m}\sum_{k=1}^{m}\sum_{j=1}^{n}M_{jk}
\end{equation*}%
provided $x_{j},a_{k},$ $j\in \left\{ 1,\dots ,n\right\} $, $k\in \left\{
1,\dots ,m\right\} $ satisfy the assumption%
\begin{equation*}
\left\vert x_{j}\right\vert \leq \func{Re}a_{k}\cdot \func{Re}x_{j}-\func{Im}%
a_{k}\cdot \func{Im}x_{j}+M_{jk}
\end{equation*}%
for each $j\in \left\{ 1,\dots ,n\right\} $, $k\in \left\{ 1,\dots
,m\right\} .$ Here $\left\vert \cdot \right\vert $ is the usual modulus of a
complex number and $M_{jk}>0,j\in \left\{ 1,\dots ,n\right\} $, $k\in
\left\{ 1,\dots ,m\right\} $ are given.
\end{remark}

\section{Karamata Type Inequalities in Hilbert Spaces}

Let $f:\left[ a,b\right] \rightarrow \mathbb{K}$, $\mathbb{K}=\mathbb{C}$ or 
$\mathbb{R}$ be a Lebesgue integrable function. The following inequality,
which is the continuous version of the \textit{triangle inequality}%
\begin{equation}
\left\vert \int_{a}^{b}f\left( x\right) dx\right\vert \leq
\int_{a}^{b}\left\vert f\left( x\right) \right\vert dx,  \label{1.1}
\end{equation}%
plays a fundamental role in Mathematical Analysis and its applications.

It appears, see \cite[p. 492]{MPF}, that the first reverse inequality for (%
\ref{1.1}) was obtained by J. Karamata in his book from 1949, \cite{K}. It
can be stated as%
\begin{equation}
\cos \theta \int_{a}^{b}\left\vert f\left( x\right) \right\vert dx\leq
\left\vert \int_{a}^{b}f\left( x\right) dx\right\vert  \label{1.2}
\end{equation}%
provided%
\begin{equation*}
-\theta \leq \arg f\left( x\right) \leq \theta ,\ \ x\in \left[ a,b\right]
\end{equation*}%
for given $\theta \in \left( 0,\frac{\pi }{2}\right) .$

This result has recently been extended by the author for the case of Bochner
integrable functions with values in a Hilbert space $H$ (see also \cite{SSD}%
):

\begin{theorem}[Dragomir, 2004]
If  $f\in L\left( \left[ a,b\right] ;H\right) $ (this means that $f:\left[
a,b\right] \rightarrow H$ is Bochner measurable on $\left[ a,b\right] $ and
the Lebesgue integral $\int_{a}^{b}\left\Vert f\left( t\right) \right\Vert dt
$ is finite), then%
\begin{equation}
\int_{a}^{b}\left\Vert f\left( t\right) \right\Vert dt\leq K\left\Vert
\int_{a}^{b}f\left( t\right) dt\right\Vert ,  \label{1.3}
\end{equation}%
provided that $f$ satisfies the condition%
\begin{equation}
\left\Vert f\left( t\right) \right\Vert \leq K\func{Re}\left\langle f\left(
t\right) ,e\right\rangle \ \ \ \text{for a.e. }t\in \left[ a,b\right] ,
\label{1.4}
\end{equation}%
where $e\in H,$ $\left\Vert e\right\Vert =1$ and $K\geq 1$ are given. 

The case of equality holds in (\ref{1.4}) if and only if%
\begin{equation}
\int_{a}^{b}f\left( t\right) dt=\frac{1}{K}\left( \int_{a}^{b}\left\Vert
f\left( t\right) \right\Vert dt\right) e.  \label{1.5}
\end{equation}
\end{theorem}

As some natural consequences of the above results, we have noticed in \cite%
{SSD} that, if $\rho \in \lbrack 0,1)$ and $f\in L\left( \left[ a,b\right]
;H\right) $ are such that%
\begin{equation}
\left\Vert f\left( t\right) -e\right\Vert \leq \rho \ \ \ \text{for a.e. }%
t\in \left[ a,b\right] ,  \label{1.6}
\end{equation}%
then%
\begin{equation}
\sqrt{1-\rho ^{2}}\int_{a}^{b}\left\Vert f\left( t\right) \right\Vert dt\leq
\left\Vert \int_{a}^{b}f\left( t\right) dt\right\Vert   \label{1.7}
\end{equation}%
with equality if and only if%
\begin{equation*}
\int_{a}^{b}f\left( t\right) dt=\sqrt{1-\rho ^{2}}\left(
\int_{a}^{b}\left\Vert f\left( t\right) \right\Vert dt\right) \cdot e.
\end{equation*}%
Also, for $e$ as above and if $M\geq m>0,$ $f\in L\left( \left[ a,b\right]
;H\right) $ such that either%
\begin{equation}
\func{Re}\left\langle Me-f\left( t\right) ,f\left( t\right) -me\right\rangle
\geq 0  \label{1.8}
\end{equation}%
or, equivalently,%
\begin{equation}
\left\Vert f\left( t\right) -\frac{M+m}{2}e\right\Vert \leq \frac{1}{2}%
\left( M-m\right)   \label{1.9}
\end{equation}%
for a.e. $t\in \left[ a,b\right] ,$ then%
\begin{equation}
\int_{a}^{b}\left\Vert f\left( t\right) \right\Vert dt\leq \frac{M+m}{2\sqrt{%
mM}}\left\Vert \int_{a}^{b}f\left( t\right) dt\right\Vert ,  \label{1.10}
\end{equation}%
with equality if and only if%
\begin{equation*}
\int_{a}^{b}f\left( t\right) dt=\frac{2\sqrt{mM}}{M+m}\left(
\int_{a}^{b}\left\Vert f\left( t\right) \right\Vert dt\right) \cdot e.
\end{equation*}

The main aim of the following sections is to extend the integral
inequalities mentioned above for the case of Banach spaces. Applications for
Hilbert spaces and for complex-valued functions are given as well.

\section{Multiplicative Reverses of the Continuous Triangle Inequality\label%
{s2}}

\subsection{The Case of One Functional}

Let $\left( X,\left\Vert \cdot \right\Vert \right) $ be a Banach space over
the real or complex number field. Then one has the following reverse of the
continuous triangle inequality \cite{SSD3a}.

\begin{theorem}[Dragomir, 2004]
\label{t2.1}Let $F$ be a continuous linear functional of unit norm on $X.$
Suppose that the function $f:\left[ a,b\right] \rightarrow X$ is Bochner
integrable on $\left[ a,b\right] $ and there exists a $r\geq 0$ such that%
\begin{equation}
r\left\Vert f\left( t\right) \right\Vert \leq \func{Re}F\left( f\left(
t\right) \right) \ \ \ \text{for a.e. }t\in \left[ a,b\right] .  \label{2.1}
\end{equation}%
Then%
\begin{equation}
r\int_{a}^{b}\left\Vert f\left( t\right) \right\Vert dt\leq \left\Vert
\int_{a}^{b}f\left( t\right) dt\right\Vert ,  \label{2.2}
\end{equation}%
where equality holds in (\ref{2.2}) if and only if both%
\begin{equation}
F\left( \int_{a}^{b}f\left( t\right) dt\right) =r\int_{a}^{b}\left\Vert
f\left( t\right) \right\Vert dt  \label{2.3}
\end{equation}%
and%
\begin{equation}
F\left( \int_{a}^{b}f\left( t\right) dt\right) =\left\Vert
\int_{a}^{b}f\left( t\right) dt\right\Vert .  \label{2.4}
\end{equation}
\end{theorem}

\begin{proof}
Since the norm of $F$ is one, then%
\begin{equation*}
\left\vert F\left( x\right) \right\vert \leq \left\Vert x\right\Vert \ \ \ 
\text{for any \ }x\in X.
\end{equation*}%
Applying this inequality for the vector $\int_{a}^{b}f\left( t\right) dt,$
we get%
\begin{align}
\left\Vert \int_{a}^{b}f\left( t\right) dt\right\Vert & \geq \left\vert
F\left( \int_{a}^{b}f\left( t\right) dt\right) \right\vert  \label{2.5} \\
& \geq \left\vert \func{Re}F\left( \int_{a}^{b}f\left( t\right) dt\right)
\right\vert =\left\vert \int_{a}^{b}\func{Re}F\left( f\left( t\right)
\right) dt\right\vert .  \notag
\end{align}%
Now, by integration of (\ref{2.1}), we obtain%
\begin{equation}
\int_{a}^{b}\func{Re}F\left( f\left( t\right) \right) dt\geq
r\int_{a}^{b}\left\Vert f\left( t\right) \right\Vert dt,  \label{2.6}
\end{equation}%
and by (\ref{2.5}) and (\ref{2.6}) we deduce the desired inequality (\ref%
{2.1}).

Obviously, if (\ref{2.3}) and (\ref{2.4}) hold true, then the equality case
holds in (\ref{2.2}).

Conversely, if the case of equality holds in (\ref{2.2}), then it must hold
in all the inequalities used before in proving this inequality. Therefore,
we must have%
\begin{equation}
r\left\Vert f\left( t\right) \right\Vert =\func{Re}F\left( f\left( t\right)
\right) \ \ \ \text{for a.e. }t\in \left[ a,b\right] ,  \label{2.7}
\end{equation}%
\begin{equation}
\func{Im}F\left( \int_{a}^{b}f\left( t\right) dt\right) =0  \label{2.8}
\end{equation}%
and%
\begin{equation}
\left\Vert \int_{a}^{b}f\left( t\right) dt\right\Vert =\func{Re}F\left(
\int_{a}^{b}f\left( t\right) dt\right) .  \label{2.9}
\end{equation}%
Integrating (\ref{2.7}) on $\left[ a,b\right] ,$ we get%
\begin{equation}
r\int_{a}^{b}\left\Vert f\left( t\right) \right\Vert dt=\func{Re}F\left(
\int_{a}^{b}f\left( t\right) dt\right) .  \label{2.10}
\end{equation}%
On utilising (\ref{2.10}) and (\ref{2.8}), we deduce (\ref{2.3}) while (\ref%
{2.9}) and (\ref{2.10}) would imply (\ref{2.4}), and the theorem is proved.
\end{proof}

\begin{corollary}
\label{c2.2}Let $\left( X,\left\Vert \cdot \right\Vert \right) $ be a Banach
space, $\left[ \cdot ,\cdot \right] :X\times X\rightarrow \mathbb{R}$ a
semi-inner product generating the norm $\left\Vert \cdot \right\Vert $ and $%
e\in X,$ $\left\Vert e\right\Vert =1.$ Suppose that the function $f:\left[
a,b\right] \rightarrow X$ is Bochner integrable on $\left[ a,b\right] $ and
there exists a $r\geq 0$ such that%
\begin{equation}
r\left\Vert f\left( t\right) \right\Vert \leq \func{Re}\left[ f\left(
t\right) ,e\right] \ \ \ \text{for a.e. }t\in \left[ a,b\right] .
\label{2.11}
\end{equation}%
Then%
\begin{equation}
r\int_{a}^{b}\left\Vert f\left( t\right) \right\Vert dt\leq \left\Vert
\int_{a}^{b}f\left( t\right) dt\right\Vert  \label{2.12}
\end{equation}%
where equality holds in (\ref{2.12}) if and only if both%
\begin{equation}
\left[ \int_{a}^{b}f\left( t\right) dt,e\right] =r\int_{a}^{b}\left\Vert
f\left( t\right) \right\Vert dt  \label{2.13}
\end{equation}%
and 
\begin{equation}
\left[ \int_{a}^{b}f\left( t\right) dt,e\right] =\left\Vert
\int_{a}^{b}f\left( t\right) dt\right\Vert .  \label{2.14}
\end{equation}
\end{corollary}

The proof follows from Theorem \ref{t2.1} for the continuous linear
functional $F\left( x\right) =\left[ x,e\right] ,$ $x\in X,$ and we omit the
details.

The following corollary of Theorem \ref{t2.1} may be stated \cite{SSD4a}.

\begin{corollary}
\label{c2.5}Let $\left( X,\left\Vert \cdot \right\Vert \right) $ be a
strictly convex Banach space, $\left[ \cdot ,\cdot \right] :X\times
X\rightarrow \mathbb{K}$ a semi-inner product generating the norm $%
\left\Vert \cdot \right\Vert $ and $e\in X,$ $\left\Vert e\right\Vert =1.$
If $f:\left[ a,b\right] \rightarrow X$ is Bochner integrable on $\left[ a,b%
\right] $ and there exists a $r\geq 0$ such that (\ref{2.11}) holds true,
then (\ref{2.12}) is valid. The case of equality holds in (\ref{2.12}) if
and only if%
\begin{equation}
\int_{a}^{b}f\left( t\right) dt=r\left( \int_{a}^{b}\left\Vert f\left(
t\right) \right\Vert dt\right) e.  \label{2.15}
\end{equation}
\end{corollary}

\begin{proof}
If (\ref{2.15}) holds true, then, obviously%
\begin{equation*}
\left\Vert \int_{a}^{b}f\left( t\right) dt\right\Vert =r\left(
\int_{a}^{b}\left\Vert f\left( t\right) \right\Vert dt\right) \left\Vert
e\right\Vert =r\int_{a}^{b}\left\Vert f\left( t\right) \right\Vert dt,
\end{equation*}%
which is the equality case in (\ref{2.12}).

Conversely, if the equality holds in (\ref{2.12}), then, by Corollary \ref%
{c2.2}, we must have (\ref{2.13}) and (\ref{2.14}). Utilising Theorem \ref%
{t2.3.1}, by (\ref{2.14}) we can conclude that there exists a $\mu >0$ such
that%
\begin{equation}
\int_{a}^{b}f\left( t\right) dt=\mu e.  \label{2.16}
\end{equation}%
Replacing this in (\ref{2.13}), we get 
\begin{equation*}
\mu \left\Vert e\right\Vert ^{2}=r\int_{a}^{b}\left\Vert f\left( t\right)
\right\Vert dt,
\end{equation*}%
giving%
\begin{equation}
\mu =r\int_{a}^{b}\left\Vert f\left( t\right) \right\Vert dt.  \label{2.17}
\end{equation}%
Utilising (\ref{2.16}) and (\ref{2.17}) we deduce (\ref{2.15}) and the proof
is completed.
\end{proof}

\subsection{The Case of $m$ Functionals}

The following result may be stated \cite{SSD4a}:

\begin{theorem}[Dragomir, 2004]
\label{t3.1}Let $\left( X,\left\Vert \cdot \right\Vert \right) $ be a Banach
space over the real or complex number field $\mathbb{K}$ and $%
F_{k}:X\rightarrow \mathbb{K}$, $k\in \left\{ 1,\dots ,m\right\} $
continuous linear functionals on $X.$ If $f:\left[ a,b\right] \rightarrow X$
is a Bochner integrable function on $\left[ a,b\right] $ and there exists $%
r_{k}\geq 0,$ $k\in \left\{ 1,\dots ,m\right\} $ with $\sum_{k=1}^{m}r_{k}>0$
and%
\begin{equation}
r_{k}\left\Vert f\left( t\right) \right\Vert \leq \func{Re}F_{k}\left(
f\left( t\right) \right)  \label{3.1}
\end{equation}%
for each $k\in \left\{ 1,\dots ,m\right\} $ and a.e. $t\in \left[ a,b\right]
,$ then%
\begin{equation}
\int_{a}^{b}\left\Vert f\left( t\right) \right\Vert dt\leq \frac{\left\Vert
\sum_{k=1}^{m}F_{k}\right\Vert }{\sum_{k=1}^{m}r_{k}}\left\Vert
\int_{a}^{b}f\left( t\right) dt\right\Vert .  \label{3.2}
\end{equation}%
The case of equality holds in (\ref{3.2}) if both%
\begin{equation}
\left( \sum_{k=1}^{m}F_{k}\right) \left( \int_{a}^{b}f\left( t\right)
dt\right) =\left( \sum_{k=1}^{m}r_{k}\right) \int_{a}^{b}\left\Vert f\left(
t\right) \right\Vert dt  \label{3.3}
\end{equation}%
and%
\begin{equation}
\left( \sum_{k=1}^{m}F_{k}\right) \left( \int_{a}^{b}f\left( t\right)
dt\right) =\left\Vert \sum_{k=1}^{m}F_{k}\right\Vert \left\Vert
\int_{a}^{b}f\left( t\right) dt\right\Vert .  \label{3.4}
\end{equation}
\end{theorem}

\begin{proof}
Utilising the hypothesis (\ref{3.1}), we have%
\begin{align}
I& :=\left\vert \sum_{k=1}^{m}F_{k}\left( \int_{a}^{b}f\left( t\right)
dt\right) \right\vert \geq \left\vert \func{Re}\left[ \sum_{k=1}^{m}F_{k}%
\left( \int_{a}^{b}f\left( t\right) dt\right) \right] \right\vert
\label{3.5} \\
& \geq \func{Re}\left[ \sum_{k=1}^{m}F_{k}\left( \int_{a}^{b}f\left(
t\right) dt\right) \right] =\sum_{k=1}^{m}\left( \int_{a}^{b}\func{Re}%
F_{k}f\left( t\right) dt\right)  \notag \\
& \geq \left( \sum_{k=1}^{m}r_{k}\right) \cdot \int_{a}^{b}\left\Vert
f\left( t\right) \right\Vert dt.  \notag
\end{align}%
On the other hand, by the continuity property of $F_{k},$ $k\in \left\{
1,\dots ,m\right\} ,$ we obviously have%
\begin{equation}
I=\left\vert \left( \sum_{k=1}^{m}F_{k}\right) \left( \int_{a}^{b}f\left(
t\right) dt\right) \right\vert \leq \left\Vert
\sum_{k=1}^{m}F_{k}\right\Vert \left\Vert \int_{a}^{b}f\left( t\right)
dt\right\Vert .  \label{3.6}
\end{equation}%
Making use of (\ref{3.5}) and (\ref{3.6}), we deduce (\ref{3.2}).

Now, obviously, if (\ref{3.3}) and (\ref{3.4}) are valid, then the case of
equality holds true in (\ref{3.2}).

Conversely, if the equality holds in the inequality (\ref{3.2}), then it
must hold in all the inequalities used to prove (\ref{3.2}), therefore we
have%
\begin{equation}
r_{k}\left\Vert f\left( t\right) \right\Vert =\func{Re}F_{k}\left( f\left(
t\right) \right)   \label{3.7}
\end{equation}%
\ for each \ $k\in \left\{ 1,\dots ,m\right\} $ \ \ and a.e. $t\in \left[ a,b%
\right] ,$%
\begin{equation}
\func{Im}\left( \sum_{k=1}^{m}F_{k}\right) \left( \int_{a}^{b}f\left(
t\right) dt\right) =0,  \label{3.8}
\end{equation}%
\begin{equation}
\func{Re}\left( \sum_{k=1}^{m}F_{k}\right) \left( \int_{a}^{b}f\left(
t\right) dt\right) =\left\Vert \sum_{k=1}^{m}F_{k}\right\Vert \left\Vert
\int_{a}^{b}f\left( t\right) dt\right\Vert .  \label{3.9}
\end{equation}

Note that, by (\ref{3.7}), on integrating on $\left[ a,b\right] $ and
summing over $k\in \left\{ 1,\dots ,m\right\} ,$ we get%
\begin{equation}
\func{Re}\left( \sum_{k=1}^{m}F_{k}\right) \left( \int_{a}^{b}f\left(
t\right) dt\right) =\left( \sum_{k=1}^{m}r_{k}\right) \int_{a}^{b}\left\Vert
f\left( t\right) \right\Vert dt.  \label{3.10}
\end{equation}%
Now, (\ref{3.8}) and (\ref{3.10}) imply (\ref{3.3}) while (\ref{3.8}) and (%
\ref{3.9}) imply (\ref{3.4}), therefore the theorem is proved.
\end{proof}

The following new results may be stated as well:

\begin{theorem}
\label{t.3.1.a}Let $\left( X,\left\Vert \cdot \right\Vert \right) $ be a
Banach space over the real or complex number field $\mathbb{K}$ and $%
F_{k}:X\rightarrow \mathbb{K}$, $k\in \left\{ 1,\dots ,m\right\} $
continuous linear functionals on $X.$ Also, assume that $f:\left[ a,b\right]
\rightarrow X$ is a Bochner integrable function on $\left[ a,b\right] $ and
there exists $r_{k}\geq 0,$ $k\in \left\{ 1,\dots ,m\right\} $ with $%
\sum_{k=1}^{m}r_{k}>0$ and%
\begin{equation*}
r_{k}\left\Vert f\left( t\right) \right\Vert \leq \func{Re}F_{k}\left(
f\left( t\right) \right)
\end{equation*}%
for each $k\in \left\{ 1,\dots ,m\right\} $ and a.e. $t\in \left[ a,b\right]
.$

(i) If $c_{\infty }$ is defined by (\ref{cinf}), then we have the inequality%
\begin{equation}
\left( 1\leq \right) \frac{\int_{a}^{b}\left\Vert f\left( t\right)
\right\Vert dt}{\left\Vert \int_{a}^{b}f\left( t\right) dt\right\Vert }\leq 
\frac{c_{\infty }}{\max_{1\leq k\leq m}\{r_{k}\}}\left( \leq \frac{%
\max_{1\leq k\leq m}\left\Vert F_{k}\right\Vert }{\max_{1\leq k\leq
m}\{r_{k}\}}\right)   \label{3.10.a}
\end{equation}%
with equality if and only if%
\begin{equation*}
\func{Re}\left( F_{k}\right) \left( \int_{a}^{b}f\left( t\right) dt\right)
=r_{k}\int_{a}^{b}\left\Vert f\left( t\right) \right\Vert dt\text{ }
\end{equation*}%
for each $k\in \left\{ 1,\dots ,m\right\} $ and%
\begin{equation*}
\max_{1\leq k\leq m}\left[ \func{Re}\left( F_{k}\right) \left(
\int_{a}^{b}f\left( t\right) dt\right) \right] =c_{\infty
}\int_{a}^{b}\left\Vert f\left( t\right) \right\Vert dt.
\end{equation*}

(ii) If $c_{p},p\geq 1,$ is defined by (\ref{cp}) , then we have the
inequality%
\begin{equation*}
\left( 1\leq \right) \frac{\int_{a}^{b}\left\Vert f\left( t\right)
\right\Vert dt}{\left\Vert \int_{a}^{b}f\left( t\right) dt\right\Vert }\leq 
\frac{c_{p}}{\left( \sum_{k=1}^{m}r_{k}^{p}\right) ^{\frac{1}{p}}}\left(
\leq \frac{\sum_{k=1}^{m}\left\Vert F_{k}\right\Vert ^{p}}{%
\sum_{k=1}^{m}r_{k}^{p}}\right) ^{\frac{1}{p}}
\end{equation*}%
with equality if and only if%
\begin{equation*}
\func{Re}\left( F_{k}\right) \left( \int_{a}^{b}f\left( t\right) dt\right)
=r_{k}\int_{a}^{b}\left\Vert f\left( t\right) \right\Vert dt\text{ }
\end{equation*}%
for each $k\in \left\{ 1,\dots ,m\right\} $ and%
\begin{equation*}
\sum_{k=1}^{m}\left[ \func{Re}F_{k}\left( \int_{a}^{b}f\left( t\right)
dt\right) \right] ^{p}=c_{p}^{p}\left\Vert \int_{a}^{b}f\left( t\right)
dt\right\Vert ^{p}
\end{equation*}%
where $p\geq 1.$
\end{theorem}

The proof is similar to the ones from Theorems \ref{ta3.1}, \ref{at3.2} and %
\ref{t3.1} and we omit the details.

The case of Hilbert spaces which provides a simpler condition for equality
is of interest for applications \cite{SSD4a}.

\begin{theorem}[Dragomir, 2004]
\label{t3.2}Let $\left( X,\left\Vert \cdot \right\Vert \right) $ be a
Hilbert space over the real or complex number field $\mathbb{K}$ and $%
e_{k}\in H\backslash \left\{ 0\right\} ,$ $k\in \left\{ 1,\dots ,m\right\} .$
If $f:\left[ a,b\right] \rightarrow H$ is a Bochner integrable function and $%
r_{k}\geq 0,$ $k\in \left\{ 1,\dots ,m\right\} $ and $\sum_{k=1}^{m}r_{k}>0$
satisfy%
\begin{equation}
r_{k}\left\Vert f\left( t\right) \right\Vert \leq \func{Re}\left\langle
f\left( t\right) ,e_{k}\right\rangle  \label{3.11}
\end{equation}%
for each $k\in \left\{ 1,\dots ,m\right\} $ and for a.e. $t\in \left[ a,b%
\right] ,$ then 
\begin{equation}
\int_{a}^{b}\left\Vert f\left( t\right) \right\Vert dt\leq \frac{\left\Vert
\sum_{k=1}^{m}e_{k}\right\Vert }{\sum_{k=1}^{m}r_{k}}\left\Vert
\int_{a}^{b}f\left( t\right) dt\right\Vert .  \label{3.12}
\end{equation}%
The case of equality holds in (\ref{3.12}) for $f\neq 0$ a.e. on $\left[ a,b%
\right] $ if and only if%
\begin{equation}
\int_{a}^{b}f\left( t\right) dt=\frac{\left( \sum_{k=1}^{m}r_{k}\right)
\int_{a}^{b}\left\Vert f\left( t\right) \right\Vert dt}{\left\Vert
\sum_{k=1}^{m}e_{k}\right\Vert ^{2}}\sum_{k=1}^{m}e_{k}.  \label{3.13}
\end{equation}
\end{theorem}

\begin{proof}
Utilising the hypothesis (\ref{3.11}) and the modulus properties, we have%
\begin{align}
\left\vert \left\langle \int_{a}^{b}f\left( t\right)
dt,\sum_{k=1}^{m}e_{k}\right\rangle \right\vert & \geq \left\vert
\sum_{k=1}^{m}\func{Re}\left\langle \int_{a}^{b}f\left( t\right)
dt,e_{k}\right\rangle \right\vert  \label{3.14} \\
& \geq \sum_{k=1}^{m}\func{Re}\left\langle \int_{a}^{b}f\left( t\right)
dt,e_{k}\right\rangle  \notag \\
& =\sum_{k=1}^{m}\int_{a}^{b}\func{Re}\left\langle f\left( t\right)
,e_{k}\right\rangle dt  \notag \\
& \geq \left( \sum_{k=1}^{m}r_{k}\right) \int_{a}^{b}\left\Vert f\left(
t\right) \right\Vert dt.  \notag
\end{align}%
By Schwarz's inequality in Hilbert spaces applied for $\int_{a}^{b}f\left(
t\right) dt$ and $\sum_{k=1}^{m}e_{k},$ we have%
\begin{equation}
\left\Vert \int_{a}^{b}f\left( t\right) dt\right\Vert \left\Vert
\sum_{k=1}^{m}e_{k}\right\Vert \geq \left\vert \left\langle
\int_{a}^{b}f\left( t\right) dt,\sum_{k=1}^{m}e_{k}\right\rangle \right\vert
.  \label{3.15}
\end{equation}%
Making use of (\ref{3.14}) and (\ref{3.15}), we deduce (\ref{3.12}).

Now, if $f\neq 0$ a.e. on $\left[ a,b\right] ,$ then $\int_{a}^{b}\left\Vert
f\left( t\right) \right\Vert dt\neq 0$ and by (\ref{3.14}) $%
\sum_{k=1}^{m}e_{k}\neq 0.$ Obviously, if (\ref{3.13}) is valid, then taking
the norm we have%
\begin{align*}
\left\Vert \int_{a}^{b}f\left( t\right) dt\right\Vert & =\frac{\left(
\sum_{k=1}^{m}r_{k}\right) \int_{a}^{b}\left\Vert f\left( t\right)
\right\Vert dt}{\left\Vert \sum_{k=1}^{m}e_{k}\right\Vert ^{2}}\left\Vert
\sum_{k=1}^{m}e_{k}\right\Vert \\
& =\frac{\sum_{k=1}^{m}r_{k}}{\left\Vert \sum_{k=1}^{m}e_{k}\right\Vert }%
\int_{a}^{b}\left\Vert f\left( t\right) \right\Vert dt,
\end{align*}%
i.e., the case of equality holds true in (\ref{3.12}).

Conversely, if the equality case holds true in (\ref{3.12}), then it must
hold in all the inequalities used to prove (\ref{3.12}), therefore we have%
\begin{equation}
\func{Re}\left\langle f\left( t\right) ,e_{k}\right\rangle =r_{k}\left\Vert
f\left( t\right) \right\Vert  \label{3.16}
\end{equation}%
for each \ $k\in \left\{ 1,\dots ,m\right\} $ \ \ and a.e. $t\in \left[ a,b%
\right] ,$%
\begin{equation}
\left\Vert \int_{a}^{b}f\left( t\right) dt\right\Vert \left\Vert
\sum_{k=1}^{m}e_{k}\right\Vert =\left\vert \left\langle \int_{a}^{b}f\left(
t\right) dt,\sum_{k=1}^{m}e_{k}\right\rangle \right\vert ,  \label{3.17}
\end{equation}%
and 
\begin{equation}
\func{Im}\left\langle \int_{a}^{b}f\left( t\right)
dt,\sum_{k=1}^{m}e_{k}\right\rangle =0.  \label{3.18}
\end{equation}%
From (\ref{3.16}) on integrating on $\left[ a,b\right] $ and summing over $k$
from $1$ to $m,$ we get%
\begin{equation}
\func{Re}\left\langle \int_{a}^{b}f\left( t\right)
dt,\sum_{k=1}^{m}e_{k}\right\rangle =\left( \sum_{k=1}^{m}r_{k}\right)
\int_{a}^{b}\left\Vert f\left( t\right) \right\Vert dt,  \label{3.19}
\end{equation}%
and then, by (\ref{3.18}) and (\ref{3.19}), we have%
\begin{equation}
\left\langle \int_{a}^{b}f\left( t\right)
dt,\sum_{k=1}^{m}e_{k}\right\rangle =\left( \sum_{k=1}^{m}r_{k}\right)
\int_{a}^{b}\left\Vert f\left( t\right) \right\Vert dt.  \label{3.20}
\end{equation}

On the other hand, by the use of the identity (\ref{Id}), the relation (\ref%
{3.17}) holds true if and only if%
\begin{equation}
\int_{a}^{b}f\left( t\right) dt=\frac{\left\langle \int_{a}^{b}f\left(
t\right) dt,\sum_{k=1}^{m}e_{k}\right\rangle }{\left\Vert
\sum_{k=1}^{m}e_{k}\right\Vert }\sum_{k=1}^{m}e_{k}.  \label{3.21}
\end{equation}%
Finally, by (\ref{3.20}) and (\ref{3.21}) we deduce that (\ref{3.13}) is
also necessary for the equality case in (\ref{3.12}) and the theorem is
proved.
\end{proof}

\begin{remark}
\label{r3.3}If $\left\{ e_{k}\right\} _{k\in \left\{ 1,\dots ,m\right\} }$
are orthogonal, then (\ref{3.12}) can be replaced by%
\begin{equation}
\int_{a}^{b}\left\Vert f\left( t\right) \right\Vert dt\leq \frac{\left(
\sum_{k=1}^{m}\left\Vert e_{k}\right\Vert ^{2}\right) ^{\frac{1}{2}}}{%
\sum_{k=1}^{m}r_{k}}\left\Vert \int_{a}^{b}f\left( t\right) dt\right\Vert ,
\label{3.22}
\end{equation}%
with equality if and only if%
\begin{equation}
\int_{a}^{b}f\left( t\right) dt=\frac{\left( \sum_{k=1}^{m}r_{k}\right)
\int_{a}^{b}\left\Vert f\left( t\right) \right\Vert dt}{\sum_{k=1}^{m}\left%
\Vert e_{k}\right\Vert ^{2}}\sum_{k=1}^{m}e_{k}.  \label{3.23}
\end{equation}%
Moreover, if $\left\{ e_{k}\right\} _{k\in \left\{ 1,\dots ,m\right\} }$ are
orthonormal, then (\ref{3.22}) becomes%
\begin{equation}
\int_{a}^{b}\left\Vert f\left( t\right) \right\Vert dt\leq \frac{\sqrt{m}}{%
\sum_{k=1}^{m}r_{k}}\left\Vert \int_{a}^{b}f\left( t\right) dt\right\Vert ,
\label{3.24}
\end{equation}%
with equality if and only if%
\begin{equation}
\int_{a}^{b}f\left( t\right) dt=\frac{1}{m}\left( \sum_{k=1}^{m}r_{k}\right)
\left( \int_{a}^{b}\left\Vert f\left( t\right) \right\Vert dt\right)
\sum_{k=1}^{m}e_{k}.  \label{3.25}
\end{equation}
\end{remark}

The following corollary of Theorem \ref{t3.2} may be stated as well \cite%
{SSD4a}.

\begin{corollary}
\label{c3.7}Let $\left( H;\left\langle \cdot ,\cdot \right\rangle \right) $
be a Hilbert space over the real or complex number field $\mathbb{K}$ and $%
e_{k}\in H\backslash \left\{ 0\right\} ,$ $k\in \left\{ 1,\dots ,m\right\} .$
If \ $f:\left[ a,b\right] \rightarrow H$ is a Bochner integrable function on 
$\left[ a,b\right] $ and $\rho _{k}>0,$ $k\in \left\{ 1,\dots ,m\right\} $
with%
\begin{equation}
\left\Vert f\left( t\right) -e_{k}\right\Vert \leq \rho _{k}<\left\Vert
e_{k}\right\Vert  \label{3.32}
\end{equation}%
\ for each \ $k\in \left\{ 1,\dots ,m\right\} $ \ \ and a.e. $t\in \left[ a,b%
\right] ,$ then%
\begin{equation}
\int_{a}^{b}\left\Vert f\left( t\right) \right\Vert dt\leq \frac{\left\Vert
\sum_{k=1}^{m}e_{k}\right\Vert }{\sum_{k=1}^{m}\left( \left\Vert
e_{k}\right\Vert ^{2}-\rho _{k}^{2}\right) ^{\frac{1}{2}}}\left\Vert
\int_{a}^{b}f\left( t\right) dt\right\Vert .  \label{3.33}
\end{equation}%
The case of equality holds in (\ref{3.33}) if and only if%
\begin{equation}
\int_{a}^{b}f\left( t\right) dt=\frac{\sum_{k=1}^{m}\left( \left\Vert
e_{k}\right\Vert ^{2}-\rho _{k}^{2}\right) ^{\frac{1}{2}}}{\left\Vert
\sum_{k=1}^{m}e_{k}\right\Vert ^{2}}\left( \int_{a}^{b}\left\Vert f\left(
t\right) \right\Vert dt\right) \sum_{k=1}^{m}e_{k}.  \label{3.34}
\end{equation}
\end{corollary}

\begin{proof}
Utilising Lemma \ref{l1.2.4}, we have from (\ref{3.32}) that%
\begin{equation*}
\left\Vert f\left( t\right) \right\Vert \left( \left\Vert e_{k}\right\Vert
^{2}-\rho _{k}^{2}\right) ^{\frac{1}{2}}\leq \func{Re}\left\langle f\left(
t\right) ,e_{k}\right\rangle
\end{equation*}%
for any $k\in \left\{ 1,\dots ,m\right\} $ and a.e. $t\in \left[ a,b\right]
. $

Applying Theorem \ref{t3.2} for%
\begin{equation*}
r_{k}:=\left( \left\Vert e_{k}\right\Vert ^{2}-\rho _{k}^{2}\right) ^{\frac{1%
}{2}},\ \ \ k\in \left\{ 1,\dots ,m\right\} ,
\end{equation*}%
we deduce the desired result.
\end{proof}

\begin{remark}
\label{r3.8}If $\left\{ e_{k}\right\} _{k\in \left\{ 1,\dots ,m\right\} }$
are orthogonal, then (\ref{3.33}) becomes%
\begin{equation}
\int_{a}^{b}\left\Vert f\left( t\right) \right\Vert dt\leq \frac{\left(
\sum_{k=1}^{m}\left\Vert e_{k}\right\Vert ^{2}\right) ^{\frac{1}{2}}}{%
\sum_{k=1}^{m}\left( \left\Vert e_{k}\right\Vert ^{2}-\rho _{k}^{2}\right) ^{%
\frac{1}{2}}}\left\Vert \int_{a}^{b}f\left( t\right) dt\right\Vert ,
\label{3.35}
\end{equation}%
with equality if and only if%
\begin{equation}
\int_{a}^{b}f\left( t\right) dt=\frac{\sum_{k=1}^{m}\left( \left\Vert
e_{k}\right\Vert ^{2}-\rho _{k}^{2}\right) ^{\frac{1}{2}}}{%
\sum_{k=1}^{m}\left\Vert e_{k}\right\Vert ^{2}}\left( \int_{a}^{b}\left\Vert
f\left( t\right) \right\Vert dt\right) \sum_{k=1}^{m}e_{k}.  \label{3.36}
\end{equation}%
Moreover, if $\left\{ e_{k}\right\} _{k\in \left\{ 1,\dots ,m\right\} }$ is
assumed to be orthonormal and%
\begin{equation*}
\left\Vert f\left( t\right) -e_{k}\right\Vert \leq \rho _{k}\text{ \ \ \ for
a.e. }t\in \left[ a,b\right] ,
\end{equation*}%
where $\rho _{k}\in \lbrack 0,1),$ $k\in \left\{ 1,\dots ,m\right\} $, then%
\begin{equation}
\int_{a}^{b}\left\Vert f\left( t\right) \right\Vert dt\leq \frac{\sqrt{m}}{%
\sum_{k=1}^{m}\left( 1-\rho _{k}^{2}\right) ^{\frac{1}{2}}}\left\Vert
\int_{a}^{b}f\left( t\right) dt\right\Vert ,  \label{3.37}
\end{equation}%
with equality iff%
\begin{equation}
\int_{a}^{b}f\left( t\right) dt=\frac{\sum_{k=1}^{m}\left( 1-\rho
_{k}^{2}\right) ^{\frac{1}{2}}}{m}\left( \int_{a}^{b}\left\Vert f\left(
t\right) \right\Vert dt\right) \sum_{k=1}^{m}e_{k}.  \label{3.38}
\end{equation}
\end{remark}

Finally, we may state the following corollary of Theorem \ref{t3.2} \cite%
{SSD3a}.

\begin{corollary}
\label{c3.10}Let $\left( H;\left\langle \cdot ,\cdot \right\rangle \right) $
be a Hilbert space over the real or complex number field $\mathbb{K}$ and $%
e_{k}\in H\backslash \left\{ 0\right\} ,$ $k\in \left\{ 1,\dots ,m\right\} .$
If $f:\left[ a,b\right] \rightarrow H$ is a Bochner integrable function on $%
\left[ a,b\right] $ and $M_{k}\geq \mu _{k}>0,$ $k\in \left\{ 1,\dots
,m\right\} $ are such that either%
\begin{equation}
\func{Re}\left\langle M_{k}e_{k}-f\left( t\right) ,f\left( t\right) -\mu
_{k}e_{k}\right\rangle \geq 0  \label{3.45}
\end{equation}%
or, equivalently%
\begin{equation}
\left\Vert f\left( t\right) -\frac{M_{k}+\mu _{k}}{2}e_{k}\right\Vert \leq 
\frac{1}{2}\left( M_{k}-\mu _{k}\right) \left\Vert e_{k}\right\Vert
\label{3.45'}
\end{equation}%
for each $k\in \left\{ 1,\dots ,m\right\} $ and a.e. $t\in \left[ a,b\right]
,$ then%
\begin{equation}
\int_{a}^{b}\left\Vert f\left( t\right) \right\Vert dt\leq \frac{\left\Vert
\sum_{k=1}^{m}e_{k}\right\Vert }{\sum_{k=1}^{m}\frac{2\cdot \sqrt{\mu
_{k}M_{k}}}{\mu _{k}+M_{k}}\left\Vert e_{k}\right\Vert }\left\Vert
\int_{a}^{b}f\left( t\right) dt\right\Vert .  \label{3.46}
\end{equation}%
The case of equality holds if and only if%
\begin{equation*}
\int_{a}^{b}f\left( t\right) dt=\frac{\sum_{k=1}^{m}\frac{2\cdot \sqrt{\mu
_{k}M_{k}}}{\mu _{k}+M_{k}}\left\Vert e_{k}\right\Vert }{\left\Vert
\sum_{k=1}^{m}e_{k}\right\Vert ^{2}}\left( \int_{a}^{b}\left\Vert f\left(
t\right) \right\Vert dt\right) \cdot \sum_{k=1}^{m}e_{k}.
\end{equation*}
\end{corollary}

\begin{proof}
Utilising Lemma \ref{l1.2.7}, by (\ref{3.45}) we deduce%
\begin{equation*}
\left\Vert f\left( t\right) \right\Vert \frac{2\cdot \sqrt{\mu _{k}M_{k}}}{%
\mu _{k}+M_{k}}\left\Vert e_{k}\right\Vert \leq \func{Re}\left\langle
f\left( t\right) ,e_{k}\right\rangle
\end{equation*}%
for each $k\in \left\{ 1,\dots ,m\right\} $ and a.e. $t\in \left[ a,b\right]
.$

Applying Theorem \ref{t3.2} for%
\begin{equation*}
r_{k}:=\frac{2\cdot \sqrt{\mu _{k}M_{k}}}{\mu _{k}+M_{k}}\left\Vert
e_{k}\right\Vert ,\ \ \ k\in \left\{ 1,\dots ,m\right\}
\end{equation*}%
we deduce the desired result.
\end{proof}

\section{Additive Reverses of the Continuous Triangle Inequality\label{ARCs2}%
}

\subsection{The Case of One Functional}

The aim of this section is to provide a different approach to the problem of
reversing the continuous triangle inequality. Namely, we are interested in
finding upper bounds for the positive difference%
\begin{equation*}
\int_{a}^{b}\left\Vert f\left( t\right) \right\Vert dt-\left\Vert
\int_{a}^{b}f\left( t\right) dt\right\Vert
\end{equation*}%
under various assumptions for the Bochner integrable function $f:\left[ a,b%
\right] \rightarrow X.$

In the following we provide an additive reverse for the continuous triangle
inequality that has been established in \cite{SSD4a}.

\begin{theorem}[Dragomir, 2004]
\label{ARCt2.1}Let $\left( X,\left\Vert \cdot \right\Vert \right) $ be a
Banach space over the real or complex number field $\mathbb{K}$ and $%
F:X\rightarrow \mathbb{K}$ be a continuous linear functional of unit norm on 
$X.$ Suppose that the function $f:\left[ a,b\right] \rightarrow X$ is
Bochner integrable on $\left[ a,b\right] $ and there exists a Lebesgue
integrable function $k:\left[ a,b\right] \rightarrow \lbrack 0,\infty )$
such that%
\begin{equation}
\left\Vert f\left( t\right) \right\Vert -\func{Re}F\left[ f\left( t\right) %
\right] \leq k\left( t\right) \   \label{ARC2.1}
\end{equation}%
for a.e. $t\in \left[ a,b\right] .$ Then we have the inequality%
\begin{equation}
\left( 0\leq \right) \int_{a}^{b}\left\Vert f\left( t\right) \right\Vert
dt-\left\Vert \int_{a}^{b}f\left( t\right) dt\right\Vert \leq
\int_{a}^{b}k\left( t\right) dt.  \label{ARC2.2}
\end{equation}%
The equality holds in (\ref{ARC2.2}) if and only if both%
\begin{equation}
F\left( \int_{a}^{b}f\left( t\right) dt\right) =\left\Vert
\int_{a}^{b}f\left( t\right) dt\right\Vert  \label{ARC2.3}
\end{equation}%
and%
\begin{equation}
F\left( \int_{a}^{b}f\left( t\right) dt\right) =\int_{a}^{b}\left\Vert
f\left( t\right) \right\Vert dt-\int_{a}^{b}k\left( t\right) dt.
\label{ARC2.3b}
\end{equation}
\end{theorem}

\begin{proof}
Since the norm of $F$ is unity, then%
\begin{equation*}
\left\vert F\left( x\right) \right\vert \leq \left\Vert x\right\Vert \ \ \ 
\text{for any \ }x\in X.
\end{equation*}%
Applying this inequality for the vector $\int_{a}^{b}f\left( t\right) dt,$
we get%
\begin{align}
\left\Vert \int_{a}^{b}f\left( t\right) dt\right\Vert & \geq \left\vert
F\left( \int_{a}^{b}f\left( t\right) dt\right) \right\vert  \label{ARC2.4} \\
& \geq \left\vert \func{Re}F\left( \int_{a}^{b}f\left( t\right) dt\right)
\right\vert  \notag \\
& =\left\vert \int_{a}^{b}\func{Re}F\left[ f\left( t\right) \right]
dt\right\vert \geq \int_{a}^{b}\func{Re}F\left[ f\left( t\right) \right] dt.
\notag
\end{align}%
Integrating (\ref{ARC2.1}), we have%
\begin{equation}
\int_{a}^{b}\left\Vert f\left( t\right) \right\Vert dt-\func{Re}F\left(
\int_{a}^{b}f\left( t\right) dt\right) \leq \int_{a}^{b}k\left( t\right) dt.
\label{ARC2.5}
\end{equation}%
Now, making use of (\ref{ARC2.4}) and (\ref{ARC2.5}), we deduce (\ref{ARC2.2}%
).

Obviously, if the equality hold in (\ref{ARC2.3}) and (\ref{ARC2.3b}), then
it holds in (\ref{ARC2.2}) as well. Conversely, if the equality holds in (%
\ref{ARC2.2}), then it must hold in all the inequalities used to prove (\ref%
{ARC2.2}). Therefore, we have%
\begin{equation*}
\int_{a}^{b}\left\Vert f\left( t\right) \right\Vert dt=\func{Re}\left[
F\left( \int_{a}^{b}f\left( t\right) dt\right) \right] +\int_{a}^{b}k\left(
t\right) dt.
\end{equation*}%
and%
\begin{equation*}
\func{Re}\left[ F\left( \int_{a}^{b}f\left( t\right) dt\right) \right]
=\left\vert F\left( \int_{a}^{b}f\left( t\right) dt\right) \right\vert
=\left\Vert \int_{a}^{b}f\left( t\right) dt\right\Vert
\end{equation*}%
which imply (\ref{ARC2.3}) and (\ref{ARC2.3b}).
\end{proof}

\begin{corollary}
\label{ARCc2.1}Let $\left( X,\left\Vert \cdot \right\Vert \right) $ be a
Banach space, $\left[ \cdot ,\cdot \right] :X\times X\rightarrow \mathbb{K}$
a semi-inner product which generates its norm. If $e\in X$ is such that $%
\left\Vert e\right\Vert =1,$ $f:\left[ a,b\right] \rightarrow X$ is Bochner
integrable on $\left[ a,b\right] $ and there exists a Lebesgue integrable
function $k:\left[ a,b\right] \rightarrow \lbrack 0,\infty )$ such that 
\begin{equation}
\left( 0\leq \right) \left\Vert f\left( t\right) \right\Vert -\func{Re}\left[
f\left( t\right) ,e\right] \leq k\left( t\right) ,  \label{ARC2.6}
\end{equation}%
for a.e. $t\in \left[ a,b\right] ,$ then 
\begin{equation}
\left( 0\leq \right) \int_{a}^{b}\left\Vert f\left( t\right) \right\Vert
dt-\left\Vert \int_{a}^{b}f\left( t\right) dt\right\Vert \leq
\int_{a}^{b}k\left( t\right) dt.  \label{ARC2.7}
\end{equation}%
where equality holds in (\ref{ARC2.7}) if and only if both%
\begin{align}
\left[ \int_{a}^{b}f\left( t\right) dt,e\right] & =\left\Vert
\int_{a}^{b}f\left( t\right) dt\right\Vert \text{ \ and}  \label{ARC2.8} \\
\left[ \int_{a}^{b}f\left( t\right) dt,e\right] & =\left\Vert
\int_{a}^{b}f\left( t\right) dt\right\Vert -\int_{a}^{b}k\left( t\right) dt.
\notag
\end{align}
\end{corollary}

The proof is obvious by Theorem \ref{ARCt2.1} applied for the continuous
linear functional of unit norm $F_{e}:X\rightarrow \mathbb{K}$, $F_{e}\left(
x\right) =\left[ x,e\right] .$

The following corollary may be stated.

\begin{corollary}
\label{ARCc2.2}Let $\left( X,\left\Vert \cdot \right\Vert \right) $ be a
strictly convex Banach space, and $\left[ \cdot ,\cdot \right] ,$ $e,$ $f,$ $%
k$ as in Corollary \ref{ARCc2.1}. Then the case of equality holds in (\ref%
{ARC2.7}) if and only if%
\begin{equation}
\int_{a}^{b}\left\Vert f\left( t\right) \right\Vert dt\geq
\int_{a}^{b}k\left( t\right) dt  \label{ARC2.9}
\end{equation}%
and%
\begin{equation}
\int_{a}^{b}f\left( t\right) dt=\left( \int_{a}^{b}\left\Vert f\left(
t\right) \right\Vert dt-\int_{a}^{b}k\left( t\right) dt\right) e.
\label{ARC2.10}
\end{equation}
\end{corollary}

\begin{proof}
Suppose that (\ref{ARC2.9}) and (\ref{ARC2.10}) are valid. Taking the norm
on (\ref{ARC2.10}) we have%
\begin{equation*}
\left\Vert \int_{a}^{b}f\left( t\right) dt\right\Vert =\left\vert
\int_{a}^{b}\left\Vert f\left( t\right) \right\Vert dt-\int_{a}^{b}k\left(
t\right) dt\right\vert \left\Vert e\right\Vert =\int_{a}^{b}\left\Vert
f\left( t\right) \right\Vert dt-\int_{a}^{b}k\left( t\right) dt,
\end{equation*}%
and the case of equality holds true in (\ref{ARC2.7}).

Now, if the equality case holds in (\ref{ARC2.7}), then obviously (\ref%
{ARC2.9}) is valid, and by Corollary \ref{ARCc2.1},%
\begin{equation*}
\left[ \int_{a}^{b}f\left( t\right) dt,e\right] =\left\Vert
\int_{a}^{b}f\left( t\right) dt\right\Vert \left\Vert e\right\Vert .
\end{equation*}%
Utilising Theorem \ref{t2.3.1}, we get%
\begin{equation}
\int_{a}^{b}f\left( t\right) dt=\lambda e\text{ \ with \ }\lambda >0.
\label{ARC2.11}
\end{equation}%
Replacing $\int_{a}^{b}f\left( t\right) dt$ with $\lambda e$ in the second
equation of (\ref{ARC2.8}) we deduce%
\begin{equation}
\lambda =\int_{a}^{b}\left\Vert f\left( t\right) \right\Vert
dt-\int_{a}^{b}k\left( t\right) dt,  \label{ARC2.12}
\end{equation}%
and by (\ref{ARC2.11}) and (\ref{ARC2.12}) we deduce (\ref{ARC2.10}).
\end{proof}

\begin{remark}
\label{ARCr2.1}If $X=H,$ $\left( H;\left\langle \cdot ,\cdot \right\rangle
\right) $ is a Hilbert space, then from Corollary \ref{ARCc2.2} we deduce
the additive reverse inequality obtained in\ \cite{SSD4}. For further
similar results in Hilbert spaces, see \cite{SSD4} and \cite{SSD5}.
\end{remark}

\subsection{The Case of $m$ Functionals}

The following result may be stated \cite{SSD4a}:

\begin{theorem}[Dragomir, 2004]
\label{ARCt3.1}Let $\left( X,\left\Vert \cdot \right\Vert \right) $ be a
Banach space over the real or complex number field $\mathbb{K}$ and $%
F_{k}:X\rightarrow \mathbb{K}$, $k\in \left\{ 1,\dots ,m\right\} $
continuous linear functionals on $X.$ If $f:\left[ a,b\right] \rightarrow X$
is a Bochner integrable function on $\left[ a,b\right] $ and $M_{k}:\left[
a,b\right] \rightarrow \lbrack 0,\infty ),\ k\in \left\{ 1,\dots ,m\right\} $
are Lebesgue integrable functions such that 
\begin{equation}
\left\Vert f\left( t\right) \right\Vert -\func{Re}F_{k}\left[ f\left(
t\right) \right] \leq M_{k}\left( t\right)  \label{ARC3.1}
\end{equation}%
for each $k\in \left\{ 1,\dots ,m\right\} $ and a.e. $t\in \left[ a,b\right]
,$ then%
\begin{equation}
\int_{a}^{b}\left\Vert f\left( t\right) \right\Vert dt\leq \left\Vert \frac{1%
}{m}\sum_{k=1}^{m}F_{k}\right\Vert \left\Vert \int_{a}^{b}f\left( t\right)
dt\right\Vert +\frac{1}{m}\sum_{k=1}^{m}\int_{a}^{b}M_{k}\left( t\right) dt.
\label{ARC3.2}
\end{equation}%
The case of equality holds in (\ref{ARC3.2}) if and only if both%
\begin{equation}
\frac{1}{m}\sum_{k=1}^{m}F_{k}\left( \int_{a}^{b}f\left( t\right) dt\right)
=\left\Vert \frac{1}{m}\sum_{k=1}^{m}F_{k}\right\Vert \left\Vert
\int_{a}^{b}f\left( t\right) dt\right\Vert  \label{ARC3.3}
\end{equation}%
and%
\begin{equation}
\frac{1}{m}\sum_{k=1}^{m}F_{k}\left( \int_{a}^{b}f\left( t\right) dt\right)
=\int_{a}^{b}\left\Vert f\left( t\right) \right\Vert dt-\frac{1}{m}%
\sum_{k=1}^{m}\int_{a}^{b}M_{k}\left( t\right) dt.  \label{ARC3.4}
\end{equation}
\end{theorem}

\begin{proof}
If we integrate on $\left[ a,b\right] $ and sum over $k$ from $1$ to $m,$ we
deduce%
\begin{equation}
\int_{a}^{b}\left\Vert f\left( t\right) \right\Vert dt\leq \frac{1}{m}%
\sum_{k=1}^{m}\func{Re}\left[ F_{k}\left( \int_{a}^{b}f\left( t\right)
dt\right) \right] +\frac{1}{m}\sum_{k=1}^{m}\int_{a}^{b}M_{k}\left( t\right)
dt.  \label{ARC3.5}
\end{equation}%
Utilising the continuity property of the functionals $F_{k}$ and the
properties of the modulus, we have:%
\begin{align}
\sum_{k=1}^{m}\func{Re}F_{k}\left( \int_{a}^{b}f\left( t\right) dt\right) &
\leq \left\vert \sum_{k=1}^{m}\func{Re}\left[ F_{k}\left(
\int_{a}^{b}f\left( t\right) dt\right) \right] \right\vert  \label{ARC3.6} \\
& \leq \left\vert \sum_{k=1}^{m}F_{k}\left( \int_{a}^{b}f\left( t\right)
dt\right) \right\vert  \notag \\
& \leq \left\Vert \sum_{k=1}^{m}F_{k}\right\Vert \left\Vert
\int_{a}^{b}f\left( t\right) dt\right\Vert .  \notag
\end{align}%
Now, by (\ref{ARC3.5}) and (\ref{ARC3.6}) we deduce (\ref{ARC3.2}).

Obviously, if (\ref{ARC3.3}) and (\ref{ARC3.4}) hold true, then the case of
equality is valid in (\ref{ARC3.2}).

Conversely, if the case of equality holds in (\ref{ARC3.2}), then it must
hold in all the inequalities used to prove (\ref{ARC3.2}). Therefore, we have%
\begin{equation*}
\int_{a}^{b}\left\Vert f\left( t\right) \right\Vert dt=\frac{1}{m}%
\sum_{k=1}^{m}\func{Re}\left[ F_{k}\left( \int_{a}^{b}f\left( t\right)
dt\right) \right] +\frac{1}{m}\sum_{k=1}^{m}\int_{a}^{b}M_{k}\left( t\right)
dt,
\end{equation*}%
\begin{equation*}
\sum_{k=1}^{m}\func{Re}\left[ F_{k}\left( \int_{a}^{b}f\left( t\right)
dt\right) \right] =\left\Vert \int_{a}^{b}f\left( t\right) dt\right\Vert
\left\Vert \sum_{k=1}^{m}F_{k}\right\Vert
\end{equation*}%
and%
\begin{equation*}
\sum_{k=1}^{m}\func{Im}\left[ F_{k}\left( \int_{a}^{b}f\left( t\right)
dt\right) \right] =0.
\end{equation*}%
These imply that (\ref{ARC3.3}) and (\ref{ARC3.4}) hold true, and the
theorem is completely proved.
\end{proof}

\begin{remark}
\label{ARCr3.1}If $F_{k},$ $k\in \left\{ 1,\dots ,m\right\} $ are of unit
norm, then, from (\ref{ARC3.2}) we deduce the inequality 
\begin{equation}
\int_{a}^{b}\left\Vert f\left( t\right) \right\Vert dt\leq \left\Vert
\int_{a}^{b}f\left( t\right) dt\right\Vert +\frac{1}{m}\sum_{k=1}^{m}%
\int_{a}^{b}M_{k}\left( t\right) dt,  \label{ARC3.7}
\end{equation}%
which is obviously coarser than (\ref{ARC3.2}) but, perhaps more useful for
applications.
\end{remark}

The following new result may be stated as well:

\begin{theorem}
\label{ARCt3.1.a}Let $\left( X,\left\Vert \cdot \right\Vert \right) $ be a
Banach space over the real or complex number field $\mathbb{K}$ and $%
F_{k}:X\rightarrow \mathbb{K}$, $k\in \left\{ 1,\dots ,m\right\} $
continuous linear functionals on $X.$ Assume also that $f:\left[ a,b\right]
\rightarrow X$ is a Bochner integrable function on $\left[ a,b\right] $ and $%
M_{k}:\left[ a,b\right] \rightarrow \lbrack 0,\infty ),\ k\in \left\{
1,\dots ,m\right\} $ are Lebesgue integrable functions such that 
\begin{equation}
\left\Vert f\left( t\right) \right\Vert -\func{Re}F_{k}\left[ f\left(
t\right) \right] \leq M_{k}\left( t\right)
\end{equation}%
for each $k\in \left\{ 1,\dots ,m\right\} $ and a.e. $t\in \left[ a,b\right]
.$

(i) If $c_{\infty }$ is defined by (\ref{cinf}), then we have the inequality%
\begin{equation}
\int_{a}^{b}\left\Vert f\left( t\right) \right\Vert dt\leq c_{\infty
}\left\Vert \int_{a}^{b}f\left( t\right) dt\right\Vert +\frac{1}{m}%
\sum_{k=1}^{m}\int_{a}^{b}M_{k}\left( t\right) dt.
\end{equation}

(ii) If $c_{p},p\geq 1,$ is defined by (\ref{cp}) , then we have the
inequality%
\begin{equation*}
\int_{a}^{b}\left\Vert f\left( t\right) \right\Vert dt\leq \frac{c_{p}}{%
m^{1/p}}\left\Vert \int_{a}^{b}f\left( t\right) dt\right\Vert +\frac{1}{m}%
\sum_{k=1}^{m}\int_{a}^{b}M_{k}\left( t\right) dt.
\end{equation*}
\end{theorem}

The proof is similar to the ones from Theorem \ref{bt3.2} and \ref{ARCt3.1}
and we omit the details.

The case of Hilbert spaces, in which one may provide a simpler condition for
equality, is of interest in applications \cite{SSD4a}.

\begin{theorem}[Dragomir, 2004]
\label{ARCt3.2}Let $\left( H,\left\langle \cdot ,\cdot \right\rangle \right) 
$ be a Hilbert space over the real or complex number field $\mathbb{K}$ and $%
e_{k}\in H,$ $k\in \left\{ 1,\dots ,m\right\} .$ If $f:\left[ a,b\right]
\rightarrow H$ is a Bochner integrable function on $\left[ a,b\right] ,$ $%
f\left( t\right) \neq 0$ for a.e. $t\in \left[ a,b\right] $ and $M_{k}:\left[
a,b\right] \rightarrow \lbrack 0,\infty ),\ k\in \left\{ 1,\dots ,m\right\} $
is a Lebesgue integrable function such that%
\begin{equation}
\left\Vert f\left( t\right) \right\Vert -\func{Re}\left\langle f\left(
t\right) ,e_{k}\right\rangle \leq M_{k}\left( t\right)  \label{ARC3.8}
\end{equation}%
for each $k\in \left\{ 1,\dots ,m\right\} $ and for a.e. $t\in \left[ a,b%
\right] ,$ then 
\begin{equation}
\int_{a}^{b}\left\Vert f\left( t\right) \right\Vert dt\leq \left\Vert \frac{1%
}{m}\sum_{k=1}^{m}e_{k}\right\Vert \left\Vert \int_{a}^{b}f\left( t\right)
dt\right\Vert +\frac{1}{m}\sum_{k=1}^{m}\int_{a}^{b}M_{k}\left( t\right) dt.
\label{ARC3.9}
\end{equation}%
The case of equality holds in (\ref{ARC3.9}) if and only if%
\begin{equation}
\int_{a}^{b}\left\Vert f\left( t\right) \right\Vert dt\geq \frac{1}{m}%
\sum_{k=1}^{m}\int_{a}^{b}M_{k}\left( t\right) dt  \label{ARC3.10}
\end{equation}%
and%
\begin{equation}
\int_{a}^{b}f\left( t\right) dt=\frac{m\left( \int_{a}^{b}\left\Vert f\left(
t\right) \right\Vert dt-\frac{1}{m}\sum_{k=1}^{m}\int_{a}^{b}M_{k}\left(
t\right) dt\right) }{\left\Vert \sum_{k=1}^{m}e_{k}\right\Vert ^{2}}%
\sum_{k=1}^{m}e_{k}.  \label{ARC3.11}
\end{equation}
\end{theorem}

\begin{proof}
As in the proof of Theorem \ref{ARCt3.1}, we have%
\begin{equation}
\int_{a}^{b}\left\Vert f\left( t\right) \right\Vert dt\leq \func{Re}%
\left\langle \frac{1}{m}\sum_{k=1}^{m}e_{k},\int_{a}^{b}f\left( t\right)
dt\right\rangle +\frac{1}{m}\sum_{k=1}^{m}\int_{a}^{b}M_{k}\left( t\right) dt
\label{ARC3.12}
\end{equation}%
and $\sum_{k=1}^{m}e_{k}\neq 0.$

On utilising Schwarz's inequality in Hilbert space $\left( H,\left\langle
\cdot ,\cdot \right\rangle \right) $ for $\int_{a}^{b}f\left( t\right) dt$
and $\sum_{k=1}^{m}e_{k},$ we have%
\begin{align}
\left\Vert \int_{a}^{b}f\left( t\right) dt\right\Vert \left\Vert
\sum_{k=1}^{m}e_{k}\right\Vert & \geq \left\vert \left\langle
\int_{a}^{b}f\left( t\right) dt,\sum_{k=1}^{m}e_{k}\right\rangle \right\vert
\label{ARC3.13} \\
& \geq \left\vert \func{Re}\left\langle \int_{a}^{b}f\left( t\right)
dt,\sum_{k=1}^{m}e_{k}\right\rangle \right\vert  \notag \\
& \geq \func{Re}\left\langle \int_{a}^{b}f\left( t\right)
dt,\sum_{k=1}^{m}e_{k}\right\rangle .  \notag
\end{align}%
By (\ref{ARC3.12}) and (\ref{ARC3.13}), we deduce (\ref{ARC3.9}).

Taking the norm on (\ref{ARC3.11}) and using (\ref{ARC3.10}), we have%
\begin{equation*}
\left\Vert \int_{a}^{b}f\left( t\right) dt\right\Vert =\frac{m\left(
\int_{a}^{b}\left\Vert f\left( t\right) \right\Vert dt-\frac{1}{m}%
\sum_{k=1}^{m}\int_{a}^{b}M_{k}\left( t\right) dt\right) }{\left\Vert
\sum_{k=1}^{m}e_{k}\right\Vert },
\end{equation*}%
showing that the equality holds in (\ref{ARC3.9}).

Conversely, if the equality case holds in (\ref{ARC3.9}), then it must hold
in all the inequalities used to prove (\ref{ARC3.9}). Therefore we have%
\begin{equation}
\left\Vert f\left( t\right) \right\Vert =\func{Re}\left\langle f\left(
t\right) ,e_{k}\right\rangle +M_{k}\left( t\right)   \label{ARC3.14}
\end{equation}%
for each $k\in \left\{ 1,\dots ,m\right\} $ and for a.e. $t\in \left[ a,b%
\right] ,$%
\begin{equation}
\left\Vert \int_{a}^{b}f\left( t\right) dt\right\Vert \left\Vert
\sum_{k=1}^{m}e_{k}\right\Vert =\left\vert \left\langle \int_{a}^{b}f\left(
t\right) dt,\sum_{k=1}^{m}e_{k}\right\rangle \right\vert   \label{ARC3.15}
\end{equation}%
and 
\begin{equation}
\func{Im}\left\langle \int_{a}^{b}f\left( t\right)
dt,\sum_{k=1}^{m}e_{k}\right\rangle =0.  \label{ARC3.16}
\end{equation}%
From (\ref{ARC3.14}) on integrating on $\left[ a,b\right] $ and summing over 
$k,$ we get%
\begin{equation}
\func{Re}\left\langle \int_{a}^{b}f\left( t\right)
dt,\sum_{k=1}^{m}e_{k}\right\rangle =m\int_{a}^{b}\left\Vert f\left(
t\right) \right\Vert dt-\sum_{k=1}^{m}\int_{a}^{b}M_{k}\left( t\right) dt.
\label{ARC3.17}
\end{equation}%
On the other hand, by the use of the identity (\ref{Id}), the relation (\ref%
{ARC3.15}) holds if and only if%
\begin{equation*}
\int_{a}^{b}f\left( t\right) dt=\frac{\left\langle \int_{a}^{b}f\left(
t\right) dt,\sum_{k=1}^{m}e_{k}\right\rangle }{\left\Vert
\sum_{k=1}^{m}e_{k}\right\Vert ^{2}}\sum_{k=1}^{m}e_{k},
\end{equation*}%
giving, from (\ref{ARC3.16}) and (\ref{ARC3.17}), that (\ref{ARC3.11}) holds
true.

If the equality holds in (\ref{ARC3.9}), then obviously (\ref{ARC3.10}) is
valid and the theorem is proved.
\end{proof}

\begin{remark}
\label{ARCr3.3}If in the above theorem, the vectors $\left\{ e_{k}\right\}
_{k\in \left\{ 1,\dots ,m\right\} }$ are assumed to be orthogonal, then (\ref%
{ARC3.9}) becomes%
\begin{equation}
\int_{a}^{b}\left\Vert f\left( t\right) \right\Vert dt\leq \frac{1}{m}\left(
\sum_{k=1}^{m}\left\Vert e_{k}\right\Vert ^{2}\right) ^{\frac{1}{2}%
}\left\Vert \int_{a}^{b}f\left( t\right) dt\right\Vert +\frac{1}{m}%
\sum_{k=1}^{m}\int_{a}^{b}M_{k}\left( t\right) dt.  \label{ARC3.18}
\end{equation}%
Moreover, if $\left\{ e_{k}\right\} _{k\in \left\{ 1,\dots ,m\right\} }$ is
an orthonormal family, then (\ref{ARC3.18}) becomes%
\begin{equation}
\int_{a}^{b}\left\Vert f\left( t\right) \right\Vert dt\leq \frac{1}{\sqrt{m}}%
\left\Vert \int_{a}^{b}f\left( t\right) dt\right\Vert +\frac{1}{m}%
\sum_{k=1}^{m}\int_{a}^{b}M_{k}\left( t\right) dt  \label{ARC3.19}
\end{equation}%
which has been obtained in \cite{SSD2}.
\end{remark}

The following corollaries are of interest.

\begin{corollary}
\label{ARCc3.1}Let $\left( H;\left\langle \cdot ,\cdot \right\rangle \right) 
$, $e_{k},$ $k\in \left\{ 1,\dots ,m\right\} $ and $f$ be as in Theorem \ref%
{ARCt3.2}. If $r_{k}:\left[ a,b\right] \rightarrow \lbrack 0,\infty ),$ $%
k\in \left\{ 1,\dots ,m\right\} $ are such that $r_{k}\in L^{2}\left[ a,b%
\right] , $ $k\in \left\{ 1,\dots ,m\right\} $ and%
\begin{equation}
\left\Vert f\left( t\right) -e_{k}\right\Vert \leq r_{k}\left( t\right) ,
\label{ARC3.20}
\end{equation}%
for each $k\in \left\{ 1,\dots ,m\right\} $ and a.e. $t\in \left[ a,b\right] 
$, then%
\begin{equation}
\int_{a}^{b}\left\Vert f\left( t\right) \right\Vert dt\leq \left\Vert \frac{1%
}{m}\sum_{k=1}^{m}e_{k}\right\Vert \left\Vert \int_{a}^{b}f\left( t\right)
dt\right\Vert +\frac{1}{2m}\sum_{k=1}^{m}\int_{a}^{b}r_{k}^{2}\left(
t\right) dt.  \label{ARC3.21}
\end{equation}%
The case of equality holds in (\ref{ARC3.21}) if and only if%
\begin{equation*}
\int_{a}^{b}\left\Vert f\left( t\right) \right\Vert dt\geq \frac{1}{2m}%
\sum_{k=1}^{m}\int_{a}^{b}r_{k}^{2}\left( t\right) dt
\end{equation*}%
and%
\begin{equation*}
\int_{a}^{b}f\left( t\right) dt=\frac{m\left( \int_{a}^{b}\left\Vert f\left(
t\right) \right\Vert dt-\frac{1}{2m}\sum_{k=1}^{m}\int_{a}^{b}r_{k}^{2}%
\left( t\right) dt\right) }{\left\Vert \sum_{k=1}^{m}e_{k}\right\Vert ^{2}}%
\sum_{k=1}^{m}e_{k}.
\end{equation*}
\end{corollary}

Finally, the following corollary may be stated.

\begin{corollary}
\label{ARCc3.2}Let $\left( H;\left\langle \cdot ,\cdot \right\rangle \right) 
$, $e_{k},$ $k\in \left\{ 1,\dots ,m\right\} $ and $f$ be as in Theorem \ref%
{ARCt3.2}. If $M_{k},\mu _{k}:\left[ a,b\right] \rightarrow \mathbb{R}$ are
such that $M_{k}\geq \mu _{k}>0$ a.e. on $\left[ a,b\right] ,$ $\frac{\left(
M_{k}-\mu _{k}\right) ^{2}}{M_{k}+\mu _{k}}\in L\left[ a,b\right] $ and%
\begin{equation*}
\func{Re}\left\langle M_{k}\left( t\right) e_{k}-f\left( t\right) ,f\left(
t\right) -\mu _{k}\left( t\right) e_{k}\right\rangle \geq 0
\end{equation*}%
for each $k\in \left\{ 1,\dots ,m\right\} $ and for a.e. $t\in \left[ a,b%
\right] ,$ then%
\begin{eqnarray*}
\int_{a}^{b}\left\Vert f\left( t\right) \right\Vert dt &\leq &\left\Vert 
\frac{1}{m}\sum_{k=1}^{m}e_{k}\right\Vert \left\Vert \int_{a}^{b}f\left(
t\right) dt\right\Vert \\
&&+\frac{1}{4m}\sum_{k=1}^{m}\left\Vert e_{k}\right\Vert ^{2}\int_{a}^{b}%
\frac{\left[ M_{k}\left( t\right) -\mu _{k}\left( t\right) \right] ^{2}}{%
M_{k}\left( t\right) +\mu _{k}\left( t\right) }dt.
\end{eqnarray*}
\end{corollary}

\section{Applications for Complex-Valued Functions}

We now give some examples of inequalities for complex-valued functions that
are Lebesgue integrable on using the general result obtained in Section \ref%
{s2}.

Consider the Banach space $\left( \mathbb{C},\left\vert \cdot \right\vert
_{1}\right) $ and $F:\mathbb{C\rightarrow C}$, $F\left( z\right) =ez$ with $%
e=\alpha +i\beta $ and $\left\vert e\right\vert ^{2}=\alpha ^{2}+\beta
^{2}=1 $, then $F$ is linear on $\mathbb{C}$. For $z\neq 0,$ we have%
\begin{equation*}
\frac{\left\vert F\left( z\right) \right\vert }{\left\vert z\right\vert _{1}}%
=\frac{\left\vert e\right\vert \left\vert z\right\vert }{\left\vert
z\right\vert _{1}}=\frac{\sqrt{\left\vert \func{Re}z\right\vert
^{2}+\left\vert \func{Im}z\right\vert ^{2}}}{\left\vert \func{Re}%
z\right\vert +\left\vert \func{Im}z\right\vert }\leq 1.
\end{equation*}%
Since, for $z_{0}=1,$ we have $\left\vert F\left( z_{0}\right) \right\vert
=1 $ and $\left\vert z_{0}\right\vert _{1}=1,$ hence%
\begin{equation*}
\left\Vert F\right\Vert _{1}:=\sup_{z\neq 0}\frac{\left\vert F\left(
z\right) \right\vert }{\left\vert z\right\vert _{1}}=1,
\end{equation*}%
showing that $F$ is a bounded linear functional on $\left( \mathbb{C}%
,\left\vert \cdot \right\vert _{1}\right) $.

Therefore we can apply Theorem \ref{t2.1} to state the following result\ for
complex-valued functions.

\begin{proposition}
\label{p4.1}Let $\alpha ,\beta \in \mathbb{R}$ with $\alpha ^{2}+\beta
^{2}=1,$ $f:\left[ a,b\right] \rightarrow \mathbb{C}$ be a Lebesgue
integrable function on $\left[ a,b\right] $ and $r\geq 0$ such that 
\begin{equation}
r\left[ \left\vert \func{Re}f\left( t\right) \right\vert +\left\vert \func{Im%
}f\left( t\right) \right\vert \right] \leq \alpha \func{Re}f\left( t\right)
-\beta \func{Im}f\left( t\right)  \label{4.1}
\end{equation}%
for a.e. $t\in \left[ a,b\right] .$ Then%
\begin{equation}
r\left[ \int_{a}^{b}\left\vert \func{Re}f\left( t\right) \right\vert
dt+\int_{a}^{b}\left\vert \func{Im}f\left( t\right) \right\vert dt\right]
\leq \left\vert \int_{a}^{b}\func{Re}f\left( t\right) dt\right\vert
+\left\vert \int_{a}^{b}\func{Im}f\left( t\right) dt\right\vert .
\label{4.2}
\end{equation}%
The equality holds in (\ref{4.2}) if and only if both%
\begin{equation*}
\alpha \int_{a}^{b}\func{Re}f\left( t\right) dt-\beta \int_{a}^{b}\func{Im}%
f\left( t\right) dt=r\left[ \int_{a}^{b}\left\vert \func{Re}f\left( t\right)
\right\vert dt+\int_{a}^{b}\left\vert \func{Im}f\left( t\right) \right\vert
dt\right]
\end{equation*}%
and%
\begin{equation*}
\alpha \int_{a}^{b}\func{Re}f\left( t\right) dt-\beta \int_{a}^{b}\func{Im}%
f\left( t\right) dt=\left\vert \int_{a}^{b}\func{Re}f\left( t\right)
dt\right\vert +\left\vert \int_{a}^{b}\func{Im}f\left( t\right)
dt\right\vert .
\end{equation*}
\end{proposition}

Now, consider the Banach space $\left( \mathbb{C},\left\vert \cdot
\right\vert _{\infty }\right) .$ If $F\left( z\right) =dz$ with $d=\gamma
+i\delta $ and $\left\vert d\right\vert =\frac{\sqrt{2}}{2},$ i.e., $\gamma
^{2}+\delta ^{2}=\frac{1}{2},$ then $F$ is linear on $\mathbb{C}$. For $%
z\neq 0$ we have%
\begin{equation*}
\frac{\left\vert F\left( z\right) \right\vert }{\left\vert z\right\vert
_{\infty }}=\frac{\left\vert d\right\vert \left\vert z\right\vert }{%
\left\vert z\right\vert _{\infty }}=\frac{\sqrt{2}}{2}\cdot \frac{\sqrt{%
\left\vert \func{Re}z\right\vert ^{2}+\left\vert \func{Im}z\right\vert ^{2}}%
}{\max \left\{ \left\vert \func{Re}z\right\vert ,\left\vert \func{Im}%
z\right\vert \right\} }\leq 1.
\end{equation*}%
Since, for $z_{0}=1+i,$ we have $\left\vert F\left( z_{0}\right) \right\vert
=1,$ $\left\vert z_{0}\right\vert _{\infty }=1,$ hence%
\begin{equation*}
\left\Vert F\right\Vert _{\infty }:=\sup_{z\neq 0}\frac{\left\vert F\left(
z\right) \right\vert }{\left\vert z\right\vert _{\infty }}=1,
\end{equation*}%
showing that $F$ is a bounded linear functional of unit norm on $\left( 
\mathbb{C},\left\vert \cdot \right\vert _{\infty }\right) $.

Therefore, we can apply Theorem \ref{t2.1}, to state the following result\
for complex-valued functions.

\begin{proposition}
\label{p4.2}Let $\gamma ,\delta \in \mathbb{R}$ with $\gamma ^{2}+\delta
^{2}=\frac{1}{2},$ $f:\left[ a,b\right] \rightarrow \mathbb{C}$ be a
Lebesgue integrable function on $\left[ a,b\right] $ and $r\geq 0$ such that%
\begin{equation*}
r\max \left\{ \left\vert \func{Re}f\left( t\right) \right\vert ,\left\vert 
\func{Im}f\left( t\right) \right\vert \right\} \leq \gamma \func{Re}f\left(
t\right) -\delta \func{Im}f\left( t\right)
\end{equation*}%
for a.e. $t\in \left[ a,b\right] .$ Then%
\begin{multline}
r\int_{a}^{b}\max \left\{ \left\vert \func{Re}f\left( t\right) \right\vert
,\left\vert \func{Im}f\left( t\right) \right\vert \right\} dt  \label{4.3} \\
\leq \max \left\{ \left\vert \int_{a}^{b}\func{Re}f\left( t\right)
dt\right\vert ,\left\vert \int_{a}^{b}\func{Im}f\left( t\right)
dt\right\vert \right\} .
\end{multline}%
The equality holds in (\ref{4.3}) if and only if both%
\begin{equation*}
\gamma \int_{a}^{b}\func{Re}f\left( t\right) dt-\delta \int_{a}^{b}\func{Im}%
f\left( t\right) dt=r\int_{a}^{b}\max \left\{ \left\vert \func{Re}f\left(
t\right) \right\vert ,\left\vert \func{Im}f\left( t\right) \right\vert
\right\} dt
\end{equation*}%
and%
\begin{equation*}
\gamma \int_{a}^{b}\func{Re}f\left( t\right) dt-\delta \int_{a}^{b}\func{Im}%
f\left( t\right) dt=\max \left\{ \left\vert \int_{a}^{b}\func{Re}f\left(
t\right) dt\right\vert ,\left\vert \int_{a}^{b}\func{Im}f\left( t\right)
dt\right\vert \right\} .
\end{equation*}
\end{proposition}

Now, consider the Banach space $\left( \mathbb{C},\left\vert \cdot
\right\vert _{2p}\right) $ with $p\geq 1.$ Let $F:\mathbb{C\rightarrow C}$, $%
F\left( z\right) =cz$ with $\left\vert c\right\vert =2^{\frac{1}{2p}-\frac{1%
}{2}}$ $\left( p\geq 1\right) .$ Obviously, $F$ is linear and by H\"{o}%
lder's inequality%
\begin{equation*}
\frac{\left\vert F\left( z\right) \right\vert }{\left\vert z\right\vert _{2p}%
}=\frac{2^{\frac{1}{2p}-\frac{1}{2}}\sqrt{\left\vert \func{Re}z\right\vert
^{2}+\left\vert \func{Im}z\right\vert ^{2}}}{\left( \left\vert \func{Re}%
z\right\vert ^{2p}+\left\vert \func{Im}z\right\vert ^{2p}\right) ^{\frac{1}{%
2p}}}\leq 1.
\end{equation*}%
Since, for $z_{0}=1+i$ we have $\left\vert F\left( z_{0}\right) \right\vert
=2^{\frac{1}{p}},$ $\left\vert z_{0}\right\vert _{2p}=2^{\frac{1}{2p}}$ $%
\left( p\geq 1\right) ,$ hence%
\begin{equation*}
\left\Vert F\right\Vert _{2p}:=\sup_{z\neq 0}\frac{\left\vert F\left(
z\right) \right\vert }{\left\vert z\right\vert _{2p}}=1,
\end{equation*}%
showing that $F$ is a bounded linear functional of unit norm on $\left( 
\mathbb{C},\left\vert \cdot \right\vert _{2p}\right) ,\left( p\geq 1\right)
. $ Therefore on using Theorem \ref{t2.1}, we may state the following result.

\begin{proposition}
\label{p4.3}Let $\varphi ,\phi \in \mathbb{R}$ with $\varphi ^{2}+\phi
^{2}=2^{\frac{1}{2p}-\frac{1}{2}}\ \left( p\geq 1\right) ,$ $f:\left[ a,b%
\right] \rightarrow \mathbb{C}$ be a Lebesgue integrable function on $\left[
a,b\right] $ and $r\geq 0$ such that%
\begin{equation*}
r\left[ \left\vert \func{Re}f\left( t\right) \right\vert ^{2p}+\left\vert 
\func{Im}f\left( t\right) \right\vert ^{2p}\right] ^{\frac{1}{2p}}\leq
\varphi \func{Re}f\left( t\right) -\phi \func{Im}f\left( t\right)
\end{equation*}%
for a.e. $t\in \left[ a,b\right] ,$ then%
\begin{multline}
r\int_{a}^{b}\left[ \left\vert \func{Re}f\left( t\right) \right\vert
^{2p}+\left\vert \func{Im}f\left( t\right) \right\vert ^{2p}\right] ^{\frac{1%
}{2p}}dt  \label{4.4} \\
\leq \left[ \left\vert \int_{a}^{b}\func{Re}f\left( t\right) dt\right\vert
^{2p}+\left\vert \int_{a}^{b}\func{Im}f\left( t\right) dt\right\vert ^{2p}%
\right] ^{\frac{1}{2p}},\ \ \left( p\geq 1\right)
\end{multline}%
where equality holds in (\ref{4.4}) if and only if both%
\begin{equation*}
\varphi \int_{a}^{b}\func{Re}f\left( t\right) dt-\phi \int_{a}^{b}\func{Im}%
f\left( t\right) dt=r\int_{a}^{b}\left[ \left\vert \func{Re}f\left( t\right)
\right\vert ^{2p}+\left\vert \func{Im}f\left( t\right) \right\vert ^{2p}%
\right] ^{\frac{1}{2p}}dt
\end{equation*}%
and%
\begin{equation*}
\varphi \int_{a}^{b}\func{Re}f\left( t\right) dt-\phi \int_{a}^{b}\func{Im}%
f\left( t\right) dt=\left[ \left\vert \int_{a}^{b}\func{Re}f\left( t\right)
dt\right\vert ^{2p}+\left\vert \int_{a}^{b}\func{Im}f\left( t\right)
dt\right\vert ^{2p}\right] ^{\frac{1}{2p}}.
\end{equation*}
\end{proposition}

\begin{remark}
\label{rnew}If $p=1$ above, and%
\begin{equation*}
r\left\vert f\left( t\right) \right\vert \leq \varphi \func{Re}f\left(
t\right) -\psi \func{Im}f\left( t\right) \text{ \ \ for a.e. }t\in \left[ a,b%
\right] ,
\end{equation*}%
provided $\varphi $, $\psi \in \mathbb{R}$ and $\varphi ^{2}+\psi
^{2}=1,r\geq 0,$ then we have a reverse of the classical continuous triangle
inequality for modulus:%
\begin{equation*}
r\int_{a}^{b}\left\vert f\left( t\right) \right\vert dt\leq \left\vert
\int_{a}^{b}f\left( t\right) dt\right\vert ,
\end{equation*}%
with equality iff%
\begin{equation*}
\varphi \int_{a}^{b}\func{Re}f\left( t\right) dt-\psi \int_{a}^{b}\func{Im}%
f\left( t\right) dt=r\int_{a}^{b}\left\vert f\left( t\right) \right\vert dt
\end{equation*}%
and%
\begin{equation*}
\varphi \int_{a}^{b}\func{Re}f\left( t\right) dt-\psi \int_{a}^{b}\func{Im}%
f\left( t\right) dt=\left\vert \int_{a}^{b}f\left( t\right) dt\right\vert .
\end{equation*}
\end{remark}

If we apply Theorem \ref{ARCt2.1}, then, in a similar manner we can prove
the following result\ for complex-valued functions.

\begin{proposition}
\label{ARCp4.1}Let $\alpha ,\beta \in \mathbb{R}$ with $\alpha ^{2}+\beta
^{2}=1,$ $f,k:\left[ a,b\right] \rightarrow \mathbb{C}$ Lebesgue integrable
functions such that 
\begin{equation*}
\left\vert \func{Re}f\left( t\right) \right\vert +\left\vert \func{Im}%
f\left( t\right) \right\vert \leq \alpha \func{Re}f\left( t\right) -\beta 
\func{Im}f\left( t\right) +k\left( t\right)
\end{equation*}%
for a.e. $t\in \left[ a,b\right] .$ Then%
\begin{multline*}
\left( 0\leq \right) \int_{a}^{b}\left\vert \func{Re}f\left( t\right)
\right\vert dt+\int_{a}^{b}\left\vert \func{Im}f\left( t\right) \right\vert
dt-\left[ \left\vert \int_{a}^{b}\func{Re}f\left( t\right) dt\right\vert
+\left\vert \int_{a}^{b}\func{Im}f\left( t\right) dt\right\vert \right] \\
\leq \int_{a}^{b}k\left( t\right) dt.
\end{multline*}
\end{proposition}

Applying Theorem \ref{ARCt2.1}, for $\left( \mathbb{C},\left\vert \cdot
\right\vert _{\infty }\right) $ we may state:

\begin{proposition}
\label{ARCp4.2}Let $\gamma ,\delta \in \mathbb{R}$ with $\gamma ^{2}+\delta
^{2}=\frac{1}{2},$ $f,k:\left[ a,b\right] \rightarrow \mathbb{C}$ Lebesgue
integrable functions on $\left[ a,b\right] $ such that%
\begin{equation*}
\max \left\{ \left\vert \func{Re}f\left( t\right) \right\vert ,\left\vert 
\func{Im}f\left( t\right) \right\vert \right\} \leq \gamma \func{Re}f\left(
t\right) -\delta \func{Im}f\left( t\right) +k\left( t\right)
\end{equation*}%
for a.e. $t\in \left[ a,b\right] .$ Then%
\begin{multline*}
\left( 0\leq \right) \int_{a}^{b}\max \left\{ \left\vert \func{Re}f\left(
t\right) \right\vert ,\left\vert \func{Im}f\left( t\right) \right\vert
\right\} dt-\max \left\{ \left\vert \int_{a}^{b}\func{Re}f\left( t\right)
dt\right\vert ,\left\vert \int_{a}^{b}\func{Im}f\left( t\right)
dt\right\vert \right\} \\
\leq \int_{a}^{b}k\left( t\right) dt.
\end{multline*}
\end{proposition}

Finally, utilising Theorem \ref{ARCt2.1}, for $\left( \mathbb{C},\left\vert
\cdot \right\vert _{2p}\right) $ with $p\geq 1,$ we may state that:

\begin{proposition}
\label{ARCp4.3}Let $\varphi ,\phi \in \mathbb{R}$ with $\varphi ^{2}+\phi
^{2}=2^{\frac{1}{2p}-\frac{1}{2}}\ \left( p\geq 1\right) ,$ $f,k:\left[ a,b%
\right] \rightarrow \mathbb{C}$ be Lebesgue integrable functions such that%
\begin{equation*}
\left[ \left\vert \func{Re}f\left( t\right) \right\vert ^{2p}+\left\vert 
\func{Im}f\left( t\right) \right\vert ^{2p}\right] ^{\frac{1}{2p}}\leq
\varphi \func{Re}f\left( t\right) -\phi \func{Im}f\left( t\right) +k\left(
t\right)
\end{equation*}%
for a.e. $t\in \left[ a,b\right] .$ Then%
\begin{multline*}
\left( 0\leq \right) \int_{a}^{b}\left[ \left\vert \func{Re}f\left( t\right)
\right\vert ^{2p}+\left\vert \func{Im}f\left( t\right) \right\vert ^{2p}%
\right] ^{\frac{1}{2p}}dt \\
-\left[ \left\vert \int_{a}^{b}\func{Re}f\left( t\right) dt\right\vert
^{2p}+\left\vert \int_{a}^{b}\func{Im}f\left( t\right) dt\right\vert ^{2p}%
\right] ^{\frac{1}{2p}}\leq \int_{a}^{b}k\left( t\right) dt.
\end{multline*}
\end{proposition}

\begin{remark}
\label{ARCrnew}If $p=1$ in the above proposition, then, from%
\begin{equation*}
\left\vert f\left( t\right) \right\vert \leq \varphi \func{Re}f\left(
t\right) -\psi \func{Im}f\left( t\right) +k\left( t\right) \text{ \ \ for
a.e. }t\in \left[ a,b\right] ,
\end{equation*}%
provided $\varphi ,\psi \in \mathbb{R}$ and $\varphi ^{2}+\psi ^{2}=1,$ we
have the additive reverse of the classical continuous triangle inequality 
\begin{equation*}
\left( 0\leq \right) \int_{a}^{b}\left\vert f\left( t\right) \right\vert
dt-\left\vert \int_{a}^{b}f\left( t\right) dt\right\vert \leq
\int_{a}^{b}k\left( t\right) dt.
\end{equation*}
\end{remark}

\end{document}